\documentclass[11pt, a4paper,onecolumn]{belarticle4}
\pdfoutput=1
\usepackage{revsymb, amsmath, amsfonts, amssymb, enumerate, fullpage, amsthm, graphicx, braket, relsize, bbm, mathrsfs, mathtools,dsfont,pifont}
\usepackage{tikzit}
\usepackage{circuitikz} 
\usepackage{graphicx,color}
\definecolor{darkbrown}{rgb}{0.4, 0.26, 0.13}
\definecolor{ao}{rgb}{0.0, 0.5, 0.0}
\definecolor{bleudefrance}{rgb}{0.19, 0.55, 0.91}
\usepackage{cancel}

\usepackage{paralist}
\usepackage{adjustbox}
\usepackage{subcaption}
\usepackage[normalem]{ulem}
\usepackage{thmtools, thm-restate}

\usepackage[linktocpage=true, colorlinks=true, linkcolor=blue, urlcolor=blue, citecolor=blue]{hyperref}

\newtheorem{theo}{Theorem}
\newtheorem{thm}[theo]{Theorem}
\newtheorem{prop}[theo]{Proposition} 
\newtheorem{lem}[theo]{Lemma}

 \newcommand{\RT}{\mathbf{RTCaus}}
 \newcommand{\RTKnow}{\mathbf{RTKnowCaus}}
  \newcommand{\Shannon}{\mathbf{Shannon}}

\newcommand{\PP}{\mathbb{P}}

\newcommand{\Id}{\mathtt{I}}
\newcommand{\F}{\mathtt{F}}
\newcommand{\Rzero}{\mathtt{R_0}}
\newcommand{\Rone}{\mathtt{R_1}}
\newcommand{\Imm}{\text{Im}}

\newcommand{\blk}{\color{black}}
\newcommand{\red}{\color{red}}
\definecolor{fuchsia}{RGB}{255,0,255}

\definecolor{orangy}{RGB}{213,94,0}

\newcommand{\rob}{\color{brown}}

\usepackage{cite}
\usepackage{tcolorbox}
\tcbuselibrary{skins,breakable}
\usetikzlibrary{shadings,shadows}

    {\endtcolorbox}

\begin{document}
\title{The resource theory of causal influence and knowledge of causal influence}

\author{Marina Maciel Ansanelli$^*$}
\affiliation{Perimeter Institute for Theoretical Physics, Waterloo, Ontario, Canada, N2L 2Y5}
\affiliation{Department of Physics and Astronomy, University of Waterloo, Waterloo, Ontario, Canada, N2L 3G1}

\author{Beata Zjawin$^*$}
\affiliation{International Centre for Theory of Quantum Technologies, University of  Gda{\'n}sk, 80-309 Gda{\'n}sk, Poland}

\author{David Schmid$^*$}
\affiliation{Perimeter Institute for Theoretical Physics, Waterloo, Ontario, Canada, N2L 2Y5}

\author{Y{\`i}l{\`e} Y{\=\i}ng} 
\affiliation{Perimeter Institute for Theoretical Physics, Waterloo, Ontario, Canada, N2L 2Y5}
\affiliation{Department of Physics and Astronomy, University of Waterloo, Waterloo, Ontario, Canada, N2L 3G1}

\author{John H.~Selby}
\affiliation{International Centre for Theory of Quantum Technologies, University of  Gda{\'n}sk, 80-309 Gda{\'n}sk, Poland}

\author{Ciarán M.~Gilligan-Lee}
\affiliation{Spotify, Dublin, Ireland}
\affiliation{Department of Physics and Astronomy, University College London, Gower Street, London, WC1E 6BT, UK.}

\author{Ana Bel\'en Sainz$^\dagger$}
\affiliation{International Centre for Theory of Quantum Technologies, University of  Gda{\'n}sk, 80-309 Gda{\'n}sk, Poland}
\affiliation{Basic Research Community for Physics e.V., Germany}

\author{Robert W.~Spekkens$^\dagger$}
\affiliation{Perimeter Institute for Theoretical Physics, Waterloo, Ontario, Canada, N2L 2Y5}
\affiliation{Department of Physics and Astronomy, University of Waterloo, Waterloo, Ontario, Canada, N2L 3G1}

\begingroup
\renewcommand{\thefootnote}{*}
\footnotetext{These authors contributed equally to this work.}
\renewcommand{\thefootnote}{$\dagger$}
\footnotetext{These authors share last authorship.}
\endgroup

\begin{abstract}
Understanding and quantifying causal relationships between variables is essential for reasoning about the physical world.  In this work, we  develop a resource-theoretic framework to do so. 
 Here, we focus on the simplest nontrivial setting--- two variables that are causally ordered, meaning that the first has the potential to influence the second, without hidden confounding. 
First, we introduce the 
resource theory that directly quantifies causal influence of a functional dependence in this setting  and show that the problem of deciding convertibility of resources and identifying a complete set of monotones has a relatively straightforward solution. 
Following this, we introduce the 
resource theory that arises naturally when one has uncertainty about the functional dependence. 
We describe a linear program for deciding  the question of whether one resource (i.e., state of knowledge about the functional dependence) can be converted to another.  Then, we focus on the case where the variables are binary. In this case, we  identify a triple of monotones that are complete in the sense that they capture the partial order over the set of all resources,  
and we provide an interpretation of each. 
\end{abstract}

\maketitle

\newpage

\tableofcontents

\newpage

\section{Introduction}\label{sec:intro}

Understanding causal relations between variables allows us to reason about the underlying mechanisms that generate observed data. Progress in causal inference has shown that it is possible to learn about the causal model from observations and interventions~\cite{pearl2009causality}, with the ultimate goal of understanding the functional dependencies that hold between the variables. This functional view of causality originates in how causal modeling was first introduced in fields such as genetics~\cite{wright1921correlation}, econometrics~\cite{Haavelmo1943}, and the social sciences~\cite{duncan2014introduction}. It also aligns naturally with approaches in physics and machine learning~\cite{glymour2014discovering,Schlkopf2022,allen2017quantum}. In functional causal models, each variable is described by a deterministic function of its parents, and any uncertainty in the value of a child variable is attributed to uncertainty in the value of an unobserved variable among its parents. This formulation eliminates the ambiguity between causation and inference and permits reasoning about counterfactuals. 

A causal structure can be represented as a directed acyclic graph (DAG) where nodes correspond to random variables and the set of directed edges into a node represents the possibility of a variable having functional dependence on its parents. Equivalently, it can be represented as a string diagram~\cite{jacobs2019causal,lorenz2023causal} where wires correspond to random variables and boxes represent functional dependencies between them. When the functional dependencies, the distributions over unobserved variables, and the causal structure are specified, they together define the functional causal model, which in turn determines the probability distribution over the observed variables. There exist functional causal models that give rise to the same observable distribution. However, even for such models, their underlying causal mechanisms, and the amount of causal influence they encode, can differ substantially. 

Consider the following example, where two different causal mechanisms are consistent with the same observational data. Imagine a simple diagram in which a binary observed variable $X$ influences another binary observed variable $Y$, and the exact functional dependence characterizing their relation is unknown (due to ignorance of an unobserved variable that influences $Y$ in the functional causal model). Instead, this dependence is represented by a probability distribution over functions, i.e., each time $Y$ is observed, the functional relationship is believed to be characterized by one of several possible functions, each assigned a probability weight. Consider two different probability distributions over functions. The first is an equal mixture of the identity and the flip functions---half of the time $Y$ equals $X$, and half of the time it is flipped. The second is an equal mixture of reset-to-one and reset-to-zero functions. Although both probability distributions over functions yield the same stochastic map from $X$ to $Y$ or, in other words, the same conditional distribution $\PP_{Y|X}$, their causal interpretation differs fundamentally. Under the identity-flip mixture, the output $Y$ still depends on $X$: if the mechanism determining which function is active was observed, the original input could be recovered by learning the output.  In contrast, under the reset functions, the output is independent of the input, hence $Y$ is causally disconnected from $X$, and all information about the input is lost. This example demonstrates that even when two functional causal models give rise to the same observable distribution, the probability distribution over the underlying functional dependences provides additional information, distinguishing mechanisms that involve nontrivial causal influences from mechanisms that do not. 

In this work, we develop a framework for quantifying  causal influence and  the knowledge about causal influence 
using the formalism of resource theories~\cite{coecke2016mathematical}. 
We introduce two distinct resource theories: the Resource Theory of Causal Influence, denoted by $\RT$, where the  resources are functional dependences, 
and the Resource Theory of Knowledge of Causal Influence, denoted by $\RTKnow$, where the resources are probability distributions over functional dependences.  We analyze the simplest nontrivial scenario: a diagram with two observed variables, $X$ and $Y$, with $X$ potentially influencing $Y$ and no confounding.  In both 
resource theories, 
the free operations are chosen such that they do not introduce new possibilities for causal influence from $X$ to $Y$. These resource theories provide a unified framework for comparing, manipulating, and classifying knowledge about causal mechanisms.

We first consider the case when the function relating $X$ and $Y$ is known. Within this setting,  we define the natural class of free resources and free operations thereon. 
We begin by considering the case where $X$ and $Y$ are binary variables. 
We characterize the conditions under which one function can be transformed into another, and describe the resulting pre-order induced by these transformations. We also extend this resource-theoretic framework beyond binary variables.

Next, we extend the framework to settings where there can be uncertainty about the functional dependence between variables, represented by a probability distribution over possible functions. 
Again, we identify the natural definition of the free resources and the free operations thereon. 
We show that resource convertibility in this setting can be decided using a linear program. In the binary case, we fully characterize the structure of the resource theory by introducing three resource monotones and proving that they form a complete set. We also provide interpretations of each of the three monotones. 
The first is the weight of functions that describe causal connection, the second is the degree of polarization between the two causally connected functions  (bias towards one or the other), 
and the third quantifies the ratio of the decrease in probability of causal connection to the increase in the degree of polarization between the two causally disconnected functions under a free operation that maximizes the latter.  The third monotone also has an alternative interpretation, that will be described in Sec.~\ref{sec_flag}. 

This work introduces a systematic framework for treating  causal dependence and knowledge of causal dependence  
as resources, offering a principled way to  quantify both.  

Although we focus here on the simplest case of two variables without confounding, the ideas presented herein lay the groundwork for a resource-theoretic treatment of causal dependence and knowledge of causal dependence in arbitrary causal structures.

\subsection{Motivation: communication channel that leaks a flag variable to the environment}
\label{sec_flag}

In this section, we present an example that distinguishes cases in which each of our resource theories is needed. 

Imagine that Alice and Bob share a communication channel that might change the content of Alice's message. That is, Alice inputs a message $X$ to the channel, and Bob receives a potentially different message $Y$. Apart from outputting $Y$, this communication channel also leaks to the environment information about which specific function it is applying to $X$ in each round. For the purpose of this work, we imagine that this information is encoded in what we call a \emph{flag variable}, which is a variable that is perfectly correlated with the mechanisms of the channel.\footnote{The flag variable here is the same as the \emph{response-function variable} from e.g., Ref.~\cite{BALKE199446}. It is encoded in the unobserved variable that influences $Y$ in the functional causal model.} In the case where $X$ and $Y$ are binary variables, this flag variable takes four possible values, associated to the four possible deterministic functions that the communication channel can apply:
\begin{itemize}
	\item Identity: $\Id$   \quad $(y=x)$, 
	\item Flip: $\F$   \quad $(y=x\oplus 1)$, 
	\item Reset to zero: $\Rzero$ \quad $(y=0)$,
	\item Reset to one: $\Rone$ \quad $(y=1)$.
\end{itemize}

In the first stage of the protocol, Alice and Bob use an additional, perfect communication channel termed a \emph{side-channel} to characterize the channel of interest. Their goal is to obtain the most accurate possible description of the channel of interest so that, in the second stage --- when the perfect side-channel is no longer available --- they can nonetheless use the imperfect (but characterized) channel to communicate effectively. 

We will consider a few different cases regarding how much information Alice and Bob have about the value of the flag variable in each run. 

Consider first the case where they do not have \emph{any} information about the flag variable. In this case, the best characterization of the communication channel that they can have is the conditional probability distribution $\PP_{Y|X}$, which they can obtain by sending many messages via the channel of interest and then using the perfect side-channel to compare the values of $X$ and $Y$ in each round. For example, they can learn that the imperfect channel is completely randomizing, that is, that there is a probability of $1/2$ that $Y=0$ and a probability of $1/2$ that $Y=1$, regardless of the value of $X$. Here and throughout the paper, we will represent a probability distribution that assigns a probability $p_{a_i}$ to each possible value $a_i$ of a variable $A$ via the square brackets notation $\PP_A=\sum_i p_{a_i}[a_i]$. Following this notation, the completely randomizing channel is
\begin{equation}
\label{eq_randomizing_channel}
    \PP_{Y|X=0}=\PP_{Y|X=1}=\frac{1}{2}[0]+\frac{1}{2}[1].
\end{equation}
Clearly, this channel is useless for transmitting information. As we will discuss in Sec.~\ref{sec_Shannon}, the relevant resource theory for the case where Alice and Bob do not have any information about the flag variable is Shannon theory~\cite{Shannon,coecke2016mathematical}. For instance, the fact that the conditional probability distribution of Eq.~\eqref{eq_randomizing_channel} is useless for communication is reflected in the fact that it is a free resource in Shannon theory. 

Consider now the case on the other extreme, where Bob learns \emph{everything} about the flag variable. That is, he knows exactly what function was applied by the communication channel in each round. If Bob learns that in a specific round the function applied was $\F$, he just needs to take the negation of the value of $Y$ in that round to learn the value of $X$.  On the other hand, if Bob learns that in a specific round the function applied was $\Rzero$,  then he knows that  there is nothing he can do to learn more about $X$. In this way, learning the flag variable allows Bob to decode Alice's message when he learns that the channel applies $\Id$ or $\F$,  whereas he knows that he cannot achieve such a decoding when he learns that the channel applies $\Rzero$ or $\Rone$. We say that  $\Id$ or $\F$ transmit causal influence, while  $\Rzero$ and $\Rone$ do not. As we will see, this is reflected in our Resource Theory of Causal Influence $\RT$ (Sec.~\ref{se:rtf}): there, $\Id$ and $\F$ are more resourceful than $\Rzero$ and $\Rone$.

Finally, consider an intermediary case, where Alice and Bob  extract \emph{partial} information about the flag variable from the environment. For example, they might learn that the probability that the communication channel implements $\Id$ is $1/2$ and that it implements $\F$ is $1/2$. In this example, their description of the channel given their partial information about the flag variable is the following \emph{probability distribution over functions}: 
\begin{equation}
    \PP_B^1 = \frac{1}{2} \, [\Id] + \frac{1}{2} \, [\F],
    \label{eq_Id_Flip} 
\end{equation}
where we again use the square brackets notation for a probability distribution, and we use the notation $\PP_B$ for probability distributions over binary functions (the subscript $B$ stands for binary). 

It is not hard to see that the probability distribution over functions $\PP_B^1$ of  Eq.~\eqref{eq_Id_Flip} gives rise to the conditional probability distribution $\PP_{Y|X}$ of Eq.~\eqref{eq_randomizing_channel}. There are, however, {\em other} probability distributions over functions that {\em also} give rise to the same $\PP_{Y|X}$. For example, $\PP_{Y|X}$ could also be obtained by a channel that always ignores $X$, and resets $Y$ to $0$ in half of the runs and to $1$ in the other half. The probability distribution over functions that describes this is denoted 
\begin{equation}
    \PP_B^2 = \frac{1}{2} \, [\Rzero] + \frac{1}{2} \, [\Rone].
    \label{eq_Rzero_Rone} 
\end{equation}

Therefore, the characterization of the channel as a probability distribution over functions is {\em strictly more informative} than the characterization of the channel as a conditional probability distribution $\PP_{Y|X}$. 

This more fine-grained information can sometimes serve as a useful resource for communication between Alice and Bob, in particular, when they can obtain nontrivial partial information about the flag variable. This is because having partial information about the flag variable can be sufficient to deduce whether or not there is any value in obtaining further information about the  flag variable. For example, the probability distribution over functions $\PP_B^1$ of  Eq.~\eqref{eq_Id_Flip} tells Bob that, \emph{if he were to extract complete information about the flag variable from the environment}, then he could infer $X$ perfectly in every run. In contrast, the probability distribution over functions  $\PP_B^2$ of Eq.~\eqref{eq_Rzero_Rone} tells Bob that further refining his knowledge of the flag variable is useless for inferring $X$. (Meanwhile, merely knowing the conditional probability distribution  $\PP_{Y|X}$ of Eq.~\eqref{eq_randomizing_channel} does not provide any such information.) In general, Bob values refining his information about the flag more highly when the probability distribution over functions has a higher weight on the functions that transmit causal influence ($\Id$ and $\F$). We will see this reflected in our Resource Theory of Knowledge of Causal Influence $\RTKnow$ (Sec.~\ref{se:rtfd}): one of our three   monotones quantifies how much weight one's probability distribution assigns to functions that transmit causal influence.

Alice and Bob might also want to know  what information Bob can extract about $X$ from the probability distribution over functions alone, \textit{without} extracting \textit{full}  information about the flag variable from the environment. In the case of the probability distributions over functions $\PP_B^1$ and $\PP_B^2$, there is nothing that Bob can learn about $X$ (even though in the case of $\PP_B^1$ he is sure he \textit{could} learn $X$ if he were to learn the flag variable). On the other hand, consider a probability distribution over functions such as 
\begin{equation}
    \PP_B^3 = \frac{2}{3} \, [\Id] + \frac{1}{3} \, [\F],
    \label{eq_Id_Flip_bias} 
\end{equation}
where there is some bias towards $\Id$. In this case, in any given run Bob knows that he has a $2/3$ chance of being correct if he sets $X=Y$. Therefore, for a fixed weight on the functions that transmit causal influence, the higher the bias towards $\Id$ or $\F$, the higher the probability that Bob can guess $X$ correctly from $Y$ without learning the flag variable.  As we will see, this is also reflected in $\RTKnow$ (Sec.~\ref{se:rtfd}): one of our three monotones quantifies the bias between $\Id$ and $\F$.

Consider now the following two probability distributions over functions:
\begin{align}
    &\PP_B^4 = \frac{1}{3} \, [\F] + \frac{2}{3} \, [\Rzero],     \label{eq_Flip_Rzero} \\
     &\PP_B^5 = \frac{1}{3} \, [\F] + \frac{1}{3} \, [\Rzero] + \frac{1}{3} \, [\Rone].
    \label{eq_Flip_Rzero_Rone} 
\end{align}

Those two probability distributions over functions have the same weight on functions that transmit causal influence (in this case, $1/3$), and the same bias between $\Id$ and $\F$ (in this case, complete bias towards $\F$). In other words, they will have the same value of the two monotones that were mentioned above. Moreover, 
the probability of Bob guessing the value of $X$ correctly from $Y$ is equal to $2/3$ for both $\PP_B^4$ and $\PP_B^5$ (see Appendix~\ref{app:guess}). Although the overall guessing probability is the same in each case, note that when the channel is characterized by $\PP_B^4$ and Bob postselects on the rounds where $Y=1$, he infers that the function implemented by the channel is $\F$. Therefore, in this case Bob has certainty that there was causal connectivity in at least some specific 
rounds (those where he receives $Y=1$), while  $\PP_B^5$ never affords Bob such certainty. Because of that, in terms of what they can learn about the functional dependence, Alice and Bob prefer a channel described by $\PP_B^4$ over a channel described by $\PP_B^5$. As we will see, this is reflected in $\RTKnow$ (Sec.~\ref{se:rtfd}): the last one of our three monotones is higher for  $\PP_B^4$ than for  $\PP_B^5$.

Note that in the specific case of $\PP_B^4$, when $Y=1$, Bob is not only certain that there was causal connectivity, but also certain about what is the \emph{value} of $X$. However, in general the third monotone only quantifies his certainty about the \emph{causal connectivity} after postselection on some specific value of $Y$, not about the value of $X$. This is illustrated by the fact (that will become clear in Sec.~\ref{se:rtfd}) that the probability distribution over functions
\begin{align}
\label{eq_P6}
    &\PP_B^6 = \frac{1}{3} \, \left(\frac{1}{2} [\Id]+\frac{1}{2}[\F]\right) + \frac{2}{3} \, [\Rzero].
\end{align}
 has the same value of the third monotone as $\PP_B^4$. Here, when $Y = 1$, Bob only knows that the channel implemented some causally connected function; he cannot tell whether it was $\Id$ or $\F$, and therefore he does not know the value of $X$. The fact that the third monotone indeed quantifies certainty about causal connectivity after postselection will be proven in Prop.~\ref{prop_thirdmon_interpretation}.

In summary, here we have outlined three different cases that Alice and Bob could be in:
\begin{enumerate}
    \item \textit{They do not know anything about the flag variable:} In this case, the conditional probability distribution $\PP_{Y|X}$ contains all of the relevant information. The relevant resource theory for this case is Shannon theory, as we discuss in Sec.~\ref{sec_Shannon}.
    \item  \textit{They know everything about the flag variable:} In this case, they know exactly which function ($\Id$, $\F$, $\Rzero$ or $\Rone$) was applied in each round. The relevant resource theory for this case is $\RT$, discussed in Sec.~\ref{se:rtf}.
    \item \textit{They have partial information about the flag variable:} In this case, their characterization of the channel is given by a probability distribution over functions. The relevant resource theory for this case is $\RTKnow$, discussed in Sec.~\ref{se:rtfd}.
\end{enumerate}

\section{Preliminaries}

\subsection{Causality}\label{sec:caus}

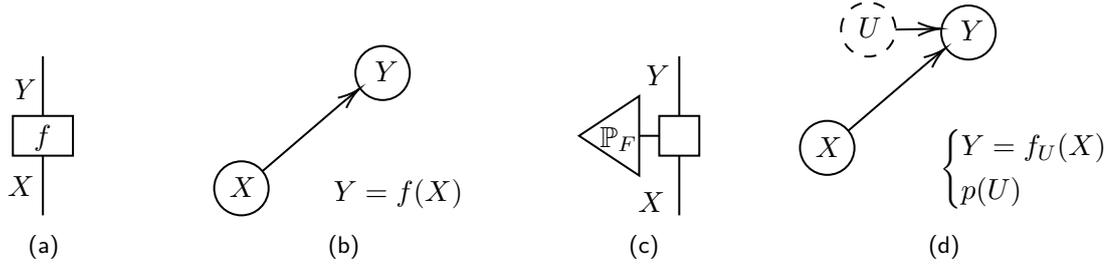
\begin{figure}[ht]
\centering
\begin{subfigure}[b]{0.24\textwidth}
  \centering
  \tikzset{every picture/.style={line width=0.75pt}}
\begin{tikzpicture}[x=0.75pt,y=0.75pt,yscale=-1,xscale=1]

\draw[line width=0.75] (55,80) -- (85,80) -- (85,100) -- (55,100) -- cycle;

\draw[line width=0.75] (70,50) -- (70,80);

\draw[line width=0.75] (70,100) -- (70,130);

\draw (64,83) node[anchor=north west][inner sep=0.75pt] {$f$};

\draw (51,109) node[anchor=north west][inner sep=0.75pt] {$X$};

\draw (52,60) node[anchor=north west][inner sep=0.75pt] {$\mathit{Y}$};

\end{tikzpicture}
  \caption{}
  \label{fig:a}
\end{subfigure}
\begin{subfigure}[b]{0.24\textwidth}
  \centering
  \tikzset{every picture/.style={line width=0.75pt}} 
\begin{tikzpicture}[x=0.75pt,y=0.75pt,yscale=-1,xscale=1]

\draw   (87.33,24.92) .. controls (87.33,17.41) and (93.41,11.33) .. (100.92,11.33) .. controls (108.42,11.33) and (114.5,17.41) .. (114.5,24.92) .. controls (114.5,32.42) and (108.42,38.5) .. (100.92,38.5) .. controls (93.41,38.5) and (87.33,32.42) .. (87.33,24.92) -- cycle ;
\draw   (16.83,82.92) .. controls (16.83,75.41) and (22.91,69.33) .. (30.42,69.33) .. controls (37.92,69.33) and (44,75.41) .. (44,82.92) .. controls (44,90.42) and (37.92,96.5) .. (30.42,96.5) .. controls (22.91,96.5) and (16.83,90.42) .. (16.83,82.92) -- cycle ;
\draw    (40.67,74.5) -- (80.61,40.02) -- (87.62,34.27) ;
\draw [shift={(89.17,33)}, rotate = 140.62] [color={rgb, 255:red, 0; green, 0; blue, 0 }  ][line width=0.75]    (10.93,-3.29) .. controls (6.95,-1.4) and (3.31,-0.3) .. (0,0) .. controls (3.31,0.3) and (6.95,1.4) .. (10.93,3.29)   ;

\draw (95.5,17.5) node [anchor=north west][inner sep=0.75pt]   [align=left] {$Y$};
\draw (23,76) node [anchor=north west][inner sep=0.75pt]   [align=left] {$X$};
\draw (75,77) node [anchor=north west][inner sep=0.75pt]   [align=left] {$Y=f(X)$};

\end{tikzpicture}
  \caption{}
  \label{fig:b}
\end{subfigure}
\begin{subfigure}[b]{0.24\textwidth}
  \centering    \tikzset{every picture/.style={line width=0.75pt}}

\begin{tikzpicture}[x=0.75pt,y=0.75pt,yscale=-1,xscale=1]

\draw [line width=0.75]    (70,45) -- (70,75) ;
\draw [line width=0.75]    (70,95) -- (70,125) ;
\draw  [line width=0.75]  (60,75) -- (80,75) -- (80,95) -- (60,95) -- cycle ;

\draw  [line width=0.75]  (50,85) -- (60,85) ;

\draw  [line width=0.75]  (20,85) -- (50,65) -- (50,105) -- cycle ;
\draw (29,79) node [anchor=north west][inner sep=0.75pt]    {$\mathbb{P}_{F}$};

\draw (48,112) node [anchor=north west][inner sep=0.75pt]    {$X$};
\draw (52,48) node [anchor=north west][inner sep=0.75pt]    {$Y$};

\end{tikzpicture}
  \caption{}
  \label{fig:c}
\end{subfigure}
\begin{subfigure}[b]{0.24\textwidth}
  \centering
  \tikzset{every picture/.style={line width=0.75pt}} 
\begin{tikzpicture}[x=0.75pt,y=0.75pt,yscale=-1,xscale=1,,baseline=(current bounding box.center)]

\draw   (87.33,24.92) .. controls (87.33,17.41) and (93.41,11.33) .. (100.92,11.33) .. controls (108.42,11.33) and (114.5,17.41) .. (114.5,24.92) .. controls (114.5,32.42) and (108.42,38.5) .. (100.92,38.5) .. controls (93.41,38.5) and (87.33,32.42) .. (87.33,24.92) -- cycle ;
\draw   (16.83,82.92) .. controls (16.83,75.41) and (22.91,69.33) .. (30.42,69.33) .. controls (37.92,69.33) and (44,75.41) .. (44,82.92) .. controls (44,90.42) and (37.92,96.5) .. (30.42,96.5) .. controls (22.91,96.5) and (16.83,90.42) .. (16.83,82.92) -- cycle ;
\draw    (40.67,74.5) -- (80.61,40.02) -- (87.62,34.27) ;
\draw [shift={(89.17,33)}, rotate = 140.62] [color={rgb, 255:red, 0; green, 0; blue, 0 }  ][line width=0.75]    (10.93,-3.29) .. controls (6.95,-1.4) and (3.31,-0.3) .. (0,0) .. controls (3.31,0.3) and (6.95,1.4) .. (10.93,3.29)   ;
\draw  [dash pattern={on 4.5pt off 4.5pt}] (37.33,23.42) .. controls (37.33,15.91) and (43.41,9.83) .. (50.92,9.83) .. controls (58.42,9.83) and (64.5,15.91) .. (64.5,23.42) .. controls (64.5,30.92) and (58.42,37) .. (50.92,37) .. controls (43.41,37) and (37.33,30.92) .. (37.33,23.42) -- cycle ;
\draw    (64.5,23.42) -- (84.67,23.04) ;
\draw [shift={(86.67,23)}, rotate = 178.92] [color={rgb, 255:red, 0; green, 0; blue, 0 }  ][line width=0.75]    (10.93,-3.29) .. controls (6.95,-1.4) and (3.31,-0.3) .. (0,0) .. controls (3.31,0.3) and (6.95,1.4) .. (10.93,3.29)   ;

\draw (95.5,17.5) node [anchor=north west][inner sep=0.75pt]   [align=left] {$Y$};
\draw (23,76) node [anchor=north west][inner sep=0.75pt]   [align=left] {$X$};
\draw (85,70) node [anchor=north west][inner sep=0.75pt]   [align=left] {$\begin{cases}Y=f_U(X)\\p(U)\end{cases}$};
\draw (44.5,15.75) node [anchor=north west][inner sep=0.75pt]   [align=left] {$U$};

\end{tikzpicture}
  \caption{}
  \label{fig:d}
\end{subfigure}
\caption{Different styles of graphical representation of causal structures. The string diagram depicted in (a) corresponds to the DAG depicted in (b), and the string diagram depicted in (c) corresponds to the DAG depicted in (d). Here, (a) (equivalently (b)) represent a causal structure with no unobserved variables, while (c) (equivalently (d)) represent a causal structure where the functional dependence of $Y$ on $X$ is sampled from a probability distribution over functions. The dashed node $U$ represents an unobserved variable. 
}
\label{fig:main}
\end{figure}

We represent causal relationships between variables using a graphical framework motivated by process theories~\cite{coecke2018picturing,Selby2021,gogioso2017fantastic} and the causal-inferential framework~\cite{omlet}. In these string diagrams~\cite{mellies2006functorial,jacobs2019causal,lorenz2023causal}, wires correspond to random variables, and boxes represent functional dependencies between variables. For example, Fig.~\ref{fig:a} illustrates two variables, where $Y$ depends on $X$ through a functional relation $Y = f(X)$. For a set of variables $V=\{X_1, \dots, X_d \}$, we say that a variable $X_1 \in V$ is a \textit{parent} of another variable $X_2 \in V$ in a diagram if $X_2$ has a direct functional dependence on $X_1$, i.e., $X_1$ and $X_2$ are respectively an input and output of the same box. The string diagram representation of causal structure can be translated to the standard DAG representation: the topological structure of the diagram corresponds to the DAG, and the labels on the diagram correspond to parameters of the causal model~\cite{lorenz2023causal}.

When a causal structure is represented by a DAG, nodes correspond to random variables and directed edges indicate direct causal influences. The string diagram illustrated in Fig.~\ref{fig:a} can be equivalently represented as a DAG as in Fig.~\ref{fig:b}. Let $V=\{X_1, \dots, X_d \}$ be the variables associated with the nodes of the DAG $G$. A node $X_1 \in V$ is called a \textit{parent} of another node $X_2 \in V$ if there is a directed edge $X_1 \to X_2$ in $G$; similarly, in this case $X_2$ is a \emph{child} of $X_1$. 

The nodes of a DAG can correspond to observed or unobserved variables (where unobserved variables are those that are not probed), which we represent by solid and dashed lines, respectively. In this work, we will only study causal structures without confounders. Therefore, the only unobserved variables we will allow for are ones that do not have any parent and have only one child in the DAG, i.e., the local error variable of a node (see Fig.~\ref{fig:d}).

To specify how variables depend on each other, one must go beyond the causal structure and also specify parameters of the causal model. A \textit{functional causal model} is defined as the specification of a causal structure (represented by a DAG or a string diagram) together with a set of functions for each observed variable and a probability distribution over each unobserved variable. Consider a causal structure that consists of observed variables $V=\{X_1, \dots, X_d \}$ and unobserved error variables $U=\{U_1, \dots, U_d \}$, where each $U_i$ is local to $X_i$ (i.e., its only child is $X_i$). Let $\mathrm{pa}(X_i)$ denote the set of all parents of $X_i$ in this causal structure. The functions $f_i$, defined by $X_i = f_i(\mathrm{pa}(X_i), U_i)$ for $i \in \{1,\dots,d\}$ represent the causal mechanism fixing the value of $X_i$ based on the values of its parents and the local unobserved variable.

Consider the causal structure represented in Fig.~\ref{fig:a}. We say that $X$ and $Y$ are \textit{causally connected} if the causal influence represented by the function $f:X \to Y$ is nontrivial, i.e., the output of the function $f$ has a nontrivial dependence on the value of $X$. If the relation is trivial, i.e., if the output of the function $f$ is \emph{independent} of the value of $X$, we say that $X$ and $Y$ are causally disconnected.

In this work, we only study the simple case of two variables $X$ and $Y$ such that $X\to Y$ and there are no confounders. We will consider two possibilities, depending on our knowledge about the functional dependence between the variables.

The first possibility is that the functional dependence of $Y$ on $X$ is described by a deterministic function whose identity does not depend on unobserved variables. This is represented by Fig.~\ref{fig:a}, where the form of the string diagram encodes causal structure and the label encodes the functional dependence, and equivalently by the DAG of Fig.~\ref{fig:b} together with the functional equation $Y=f(X)$. 

The second possibility is that the identity of the function from $X$ to $Y$ depends on an unobserved local error variable $U$, i.e., $Y=f_U(X)$. Since the value of $U$ is unknown, this equivalently describes a situation  where the functional relationship of $Y$ on $X$ is not known with certainty, but instead is described by a probability distribution over possible functions. Let $F(X\to Y)$ be the space of all possible deterministic functions from $X$ to $Y$. We represent a probability distribution over functions from $X$ to $Y$ as $\PP_{F}=\sum_u \PP_{F}(f_u) [f_u]$, with $\sum_u \PP_{F}(f_u) =1$, where $[f]$ denotes a point distribution associated with the function $f:X \to Y$. This case is depicted in Fig.~\ref{fig:c}, which explicitly encodes both the causal structure and the probability distribution over functions, and equivalently by the DAG of Fig.~\ref{fig:d} together with the functional equation $Y=f_U(X)$ and the probability distribution over the error variable, $\PP_U$. 

In this paper, we sometimes restrict the scenario to binary $X$ and $Y$. We denote the set of all bit-to-bit functions by $B := F(\text{bit} \to \text{bit})$  and the corresponding probability distributions over functions by $\PP_B$.

At the level of observed statistics, we can obtain the stochastic map that determines the value of $Y$ given $X$ from the probability distribution over functions $\PP_{F}$ as $\PP_{Y|X}= \sum_{f\in \mathcal{F}} \PP_{F}(f)\,\delta_{Y,f(X)}$. Note that two different probability distributions over functions can give rise to the same stochastic map, e.g., those of Eqs.~\eqref{eq_Id_Flip} and~\eqref{eq_Rzero_Rone}, which both give rise to the stochastic map of Eq.~\eqref{eq_randomizing_channel}. 
In this work, we are motivated by the fact that even probability distributions over functions that give rise to the same stochastic map
can differ in the amount of causal influence they represent between variables. For this reason, we study probability distributions over functions as meaningful objects in their own right, even when they give rise to the same stochastic map.

\subsection{Resource theories}
\label{sec_RT_preliminaries}

Resource theories provide a general framework for analyzing properties of processes (including states, channels, measurements, combs etc)~\cite{coecke2016mathematical}. The basic idea of this framework is to classify processes into two classes: \emph{free processes}, which are assumed to be freely available and constitute the free resources, and all others, which constitute the nonfree resources. Formally, a resource theory is defined by identifying an \emph{enveloping theory}, which defines a set of processes that is closed under parallel and sequential composition, and a \emph{free subtheory}, which consists of a closed subset of processes considered to be free. When a free process is used to achieve interconversion of one resource into another, it is referred to as a \emph{free operation}.

A central question in any resource theory is whether one resource can be converted into another using only free operations. Given two resources $R_1$ and $R_2$, we write $R_1 \to R_2$ if such a conversion is possible, and $ R_1 \not\to R_2$ otherwise. These relations define a \textit{pre-order} on the set of resources: they are reflexive (any resource can be converted into itself) and transitive (if $R_1  \to  R_2$ and $ R_2 \to  R_3 $, then $ R_1 \to  R_3 $). For any resource $R$, the set of all resources to which it can be converted is called its \textit{downward closure}. Resources that are mutually interconvertible form \textit{equivalence classes}. The set of these equivalence classes, along with the induced conversion relation between them, forms a \textit{partial order}.

A useful tool in analyzing resource conversion is the notion of a resource monotone. A monotone is a real-valued function $M$ that is non-increasing under free operations, that is, one that satisfies the constraint that if $R_1 \to R_2$, then $M(R_1) \geq M(R_2)$. In particular, if for a pair of resources $R_1$ and $R_2$ we find that $M(R_1) < M(R_2)$, this certifies that $R_1 \not\to R_2$. However, a single monotone is generally insufficient to fully characterize the pre-order; instead, a family of monotones is typically required.

One method for constructing resource monotones is called the \textit{cost-construction}, which takes any function (monotone or not) together with a set of resources, and produces a resource monotone from them. More specifically, take a set of resources $\mathbf{S}$ and a function $f$ valued on the resources of your enveloping theory. For any resource $R$, the following is the cost-construction that defines a resource monotone~\cite{gonda2019monotones}:
\begin{align}\label{eq:cost}
    M^{\text{cost}}_{f,\mathbf{S}}[R] := \mathop{\mathrm{min}}_{R^* \in \mathbf{S}}
\, \left\{f(R^*) \quad \text{s.t.} \quad  R^* \to R \right\}\,.
\end{align}

\section{The Resource Theory of Causal Influence}\label{se:rtf}

In this section, we develop the \textit{Resource Theory of Causal Influence}, which we denote by $\RT$. We study functional causal structures in which the relation between two variables is given by a single deterministic function, and we aim to quantify the strength of this causal influence using a resource-theoretic approach. We focus on the simplest case, where we have only two  observed finite cardinality variables ($X$ and $Y$) and where the causal influence goes from $X$ to $Y$.

\subsection{Resources and transformations}

In this resource theory, the objects of interest are \emph{deterministic processes}. The first type of deterministic process we consider is the set of  functions $f \, : \, X \to Y$. The second type of deterministic process we consider is the set of \emph{deterministic combs}, represented by $\tau$. A \emph{deterministic comb} is an operation taking a function $f \, : \, X \to Y$ into another function $g \, : \, X' \to Y'$ (whose domain and co-domain might have different cardinalities):
\begin{align}\label{fig:f_trans}
\tikzset{every picture/.style={line width=0.75pt}}
\begin{tikzpicture}[x=0.75pt,y=0.75pt,yscale=-1,xscale=1]

\draw[line width=0.75] (55,80) -- (85,80) -- (85,101) -- (55,101) -- cycle;

\draw[line width=0.75] (70,51) -- (70,80);

\draw[line width=0.75] (70,101) -- (70,130);

\draw[line width=0.75] (190,80) -- (220,80) -- (220,100) -- (190,100) -- cycle;

\draw[line width=0.75] (205,60) -- (205,80);

\draw[line width=0.75] (205,100) -- (205,120);

\draw[line width=0.75] (205,135) -- (205,159);

\draw[line width=0.75] (205,21) -- (205,45);

\draw[color={rgb, 255:red, 155; green, 155; blue, 155},draw opacity=1,fill={rgb, 255:red, 155; green, 155; blue, 155},fill opacity=1] (226,60) -- (226,120) -- (184,120) -- (184,60) -- (226,60) -- cycle (166,46) -- (166,135) -- (243,135) -- (243,46) -- (166,46) -- cycle;

\draw (65,86) node[anchor=north west][inner sep=0.75pt] {$g$};

\draw (48,109) node[anchor=north west][inner sep=0.75pt] {$X'$};

\draw (49,59) node[anchor=north west][inner sep=0.75pt] {$\mathit{Y}'$};

\draw (115,86) node[anchor=north west][inner sep=0.75pt] {$=$};

\draw (199,83) node[anchor=north west][inner sep=0.75pt] {$f$};

\draw (188,103) node[anchor=north west][inner sep=0.75pt] {$X$};

\draw (188,62) node[anchor=north west][inner sep=0.75pt] {$\mathit{Y}$};

\draw (184,141) node[anchor=north west][inner sep=0.75pt] {$X'$};

\draw (184,21) node[anchor=north west][inner sep=0.75pt] {$\mathit{Y}'$};

\draw (230,120) node[anchor=north west][inner sep=0.75pt] {$\tau$};

\end{tikzpicture}.
\end{align}
(This diagrammatic representation is motivated by process theories~\cite{coecke2018picturing,Selby2021,gogioso2017fantastic} and by the framework of combs~\cite{chiribella2009theoretical}, where operations on channels are depicted abstractly in this way.) 

The deterministic combs always consist of a pre-processing of $X'$ by some function $f_\text{pre}:X'\to X\times Z$ and a post-processing of $Y$ by some function $f_\text{post}:Y\times Z\to Y'$ potentially connected by a causal link $Z$. That is, they can be written as:
\begin{align}\label{fig:det_comb}
\tikzset{every picture/.style={line width=0.75pt}}
\begin{tikzpicture}[x=0.75pt,y=0.75pt,yscale=-1,xscale=1]

\draw (65,95) node[anchor=north west][inner sep=0.75pt] {$X$};

\draw (65,49) node[anchor=north west][inner sep=0.75pt] {$Y$};

\draw (75,134) node[anchor=north west][inner sep=0.75pt] {$X'$};

\draw (75,8) node[anchor=north west][inner sep=0.75pt] {$Y'$};

\draw (107,70) node[anchor=north west][inner sep=0.75pt] {$Z$};

\draw (80,112) node[anchor=north west][inner sep=0.75pt] {$f_\text{pre}$};

\draw (78,30) node[anchor=north west][inner sep=0.75pt] {$f_\text{post}$};

\draw[line width=0.75] (84,48) -- (84,64);

\draw[line width=0.75] (94,129) -- (94,153);

\draw (77,110) rectangle (112,129);

\draw (77,28) rectangle (112,48);

\draw[line width=0.75] (94,12) -- (94,28);

\draw[line width=0.75] (84,93) -- (84,110);

\draw[line width=0.75] (104,48) -- (104,110);

\end{tikzpicture}.
\end{align} 

An example of such an operation applied on a function is given by
\begin{align}\label{fig:f_id}
\tikzset{every picture/.style={line width=0.75pt}}
\begin{tikzpicture}[x=0.75pt,y=0.75pt,yscale=-1,xscale=1]

\draw[color={rgb, 255:red, 155; green, 155; blue, 155},draw opacity=1,fill={rgb, 255:red, 155; green, 155; blue, 155},fill opacity=1] (243,46) -- (243,136) -- (166,136) -- (166,46) -- (243,46) -- cycle (185,71) -- (185,111) -- (226,111) -- (226,71) -- (185,71) -- cycle;

\draw[line width=0.75] (55,80) -- (85,80) -- (85,101) -- (55,101) -- cycle;

\draw[line width=0.75] (70,51) -- (70,80);

\draw[line width=0.75] (70,101) -- (70,130);

\draw[line width=0.75] (190,80) -- (220,80) -- (220,100) -- (190,100) -- cycle;

\draw[line width=0.75] (205,66) -- (205,80);

\draw[line width=0.75] (205,100) -- (205,115);

\draw[line width=0.75] (205,135) -- (205,158);

\draw[line width=0.75] (205,21) -- (205,45);

\draw (175,57) -- (175,122);

\draw (175,57) .. controls (175,47) and (204,56) .. (205,45);

\draw (175,122) .. controls (175,132) and (204,125) .. (205,135);

\draw[line width=0.75] (196,66) -- (214,66);

\draw[line width=0.75] (199,61) -- (211,61);

\draw[line width=0.75] (202,57) -- (208,57);

\draw (205,125) -- (192,115) -- (218,115) -- cycle;

\draw (65,86) node[anchor=north west][inner sep=0.75pt] {$\Id$};

\draw (48,109) node[anchor=north west][inner sep=0.75pt] {$X'$};

\draw (49,59) node[anchor=north west][inner sep=0.75pt] {$Y'$};

\draw (115,86) node[anchor=north west][inner sep=0.75pt] {$=$};

\draw (199,83) node[anchor=north west][inner sep=0.75pt] {$f$};

\draw (200,115) node[anchor=north west][inner sep=0.75pt] {$_x$};

\draw (184,141) node[anchor=north west][inner sep=0.75pt] {$X'$};

\draw (184,21) node[anchor=north west][inner sep=0.75pt] {$Y'$};

\end{tikzpicture}.
\end{align}
which effectively maps any function $f \, : \, X \to Y$ to the identity function from $X'$ to $Y'$. This operation uses a side channel that feeds the input $X'$ to the output $Y'$, independent of what the function $f$ was.

All of the processes we consider admit of a bipartition of the involved systems between a sender and a receiver. For example, a function $f:X\to Y$ admits of the bipartition where $X$ belongs to the sender and $Y$ belongs to the receiver. A deterministic comb $\tau:F(X\to Y)\to F(X'\to Y')$ admits of the bipartition where $X$ and $X'$ belong to the sender and $Y$ and $Y'$ belong to the receiver. Our aim here is to study the causal dependence \emph{across} this bipartition. This allows us to introduce the notion of a \emph{free} process in a type-independent way, that is, a definition that can be applied to any type of deterministic process (functions, deterministic combs, etc): the free processes are those for which there is no cause-effect connection across the bipartition.

We begin by specializing this notion of freeness to the case of functions. A function $f:X\to Y$ is considered free if $Y$ has no causal dependence on $X$, that is, if
\begin{align}\label{fig:f_reset}
\tikzset{every picture/.style={line width=0.75pt}}
\begin{tikzpicture}[x=0.75pt,y=0.75pt,yscale=-1,xscale=1]

\draw[line width=0.75] (20,50) -- (20,75);

\draw[line width=0.75] (20,104) -- (20,129);

\draw[line width=0.75] (120,100) -- (120,129);

\draw[line width=0.75] (120,51) -- (120,80);

\draw[line width=0.75] (111,100) -- (129,100);

\draw[line width=0.75] (114,97) -- (126,97);

\draw[line width=0.75] (117,95) -- (123,95);

\draw (120,91) -- (110,80) -- (130,80) -- cycle;

\draw[line width=0.75] (105,75) -- (135,75) -- (135,104) -- (105,104) -- cycle;

\draw[line width=0.75] (5,75) -- (35,75) -- (35,104) -- (5,104) -- cycle;

\draw[line width=0.75] (210,50) -- (210,75);

\draw[line width=0.75] (210,104) -- (210,129);

\draw[line width=0.75] (195,75) -- (225,75) -- (225,104) -- (195,104) -- cycle;

\draw (97,40) .. controls (92,40) and (90,42) .. (90,47) -- (90,80) .. controls (90,87) and (87,90) .. (83,90) .. controls (87,90) and (90,93) .. (90,100) -- (90,133) .. controls (90,138) and (92,140) .. (97,140);

\draw (234,140) .. controls (238,140) and (241,138) .. (241,133) -- (241,100) .. controls (241,93) and (243,90) .. (247,90) .. controls (243,90) and (240,87) .. (240,80) -- (240,47) .. controls (240,42) and (238,40) .. (233,40);

\draw (6,82) node[anchor=north west][inner sep=0.75pt] {$f_{\mathrm{free}}$};

\draw (0,112) node[anchor=north west][inner sep=0.75pt] {$X$};

\draw (0,55) node[anchor=north west][inner sep=0.75pt] {$\mathit{Y}$};

\draw (45,82) node[anchor=north west][inner sep=0.75pt] {$\in$};

\draw (116,80) node[anchor=north west][inner sep=0.75pt] {$_{y}$};

\draw (101,112) node[anchor=north west][inner sep=0.75pt] {$X$};

\draw (101,54) node[anchor=north west][inner sep=0.75pt] {$\mathit{Y}$};

\draw (157,84) node[anchor=north west][inner sep=0.75pt] {$=:$};

\draw (201,83) node[anchor=north west][inner sep=0.75pt] {$\mathtt{R_y}$};

\draw (188,112) node[anchor=north west][inner sep=0.75pt] {$X$};

\draw (189,55) node[anchor=north west][inner sep=0.75pt] {$\mathit{Y}$};

\draw (250,130) node[anchor=north west][inner sep=0.75pt] {$y\in Y$};

\end{tikzpicture}.
\end{align}
We will refer to such an operation as ``reset to $y$'' and denote it by $\mathtt{R_y}$.  

Next, we specialize the notion of freeness to the case of deterministic combs. A deterministic comb $\tau$ is considered free if it does not have a wire from the pre-processing to the post-processing. Such a comb, when applied to a function, does not introduce any causal influence between the input and output variable that is not already present in the function. It is clear, then, that the operation depicted in Eq.~\eqref{fig:f_id} is \emph{not}  an example of a free deterministic comb. 
The free deterministic combs consist of the set wherein 
there is a pre-processing of the input of $f$ and a post-processing of the output of $f$; these  can be depicted as follows:
\begin{align}\label{fig:f_free}
\tikzset{every picture/.style={line width=0.75pt}}
\begin{tikzpicture}[x=0.75pt,y=0.75pt,yscale=-1,xscale=1]

\draw[color={rgb, 255:red, 155; green, 155; blue, 155},draw opacity=1,fill={rgb, 255:red, 155; green, 155; blue, 155},fill opacity=1] (395,45) -- (395,135) -- (318,135) -- (318,45) -- (395,45) -- cycle (337,70) -- (337,110) -- (378,110) -- (378,70) -- (337,70) -- cycle;

\draw[line width=0.75] (205,66) -- (205,86);

\draw[line width=0.75] (205,94) -- (205,114);

\draw[line width=0.75] (205,129) -- (205,153);

\draw[line width=0.75] (205,30) -- (205,54);

\draw[color={rgb, 255:red, 155; green, 155; blue, 155},draw opacity=1,fill={rgb, 255:red, 155; green, 155; blue, 155},fill opacity=1] (244,45) -- (244,135) -- (167,135) -- (167,45) -- (244,45) -- cycle (185,70) -- (185,110) -- (227,110) -- (227,70) -- (185,70) -- cycle;

\draw[line width=0.75] (357,65) -- (357,85);

\draw[line width=0.75] (357,93) -- (357,113);

\draw[line width=0.75] (357,128) -- (357,152);

\draw[line width=0.75] (357,28) -- (357,50);

\draw (340,50) -- (374,50) -- (374,65) -- (340,65) -- cycle;

\draw (340,113) -- (374,113) -- (374,128) -- (340,128) -- cycle;

\draw (115,85) node[anchor=north west][inner sep=0.75pt] {$:=$};

\draw (189,95) node[anchor=north west][inner sep=0.75pt] {$X$};

\draw (189,70) node[anchor=north west][inner sep=0.75pt] {$\mathit{Y}$};

\draw (185,136) node[anchor=north west][inner sep=0.75pt] {$X'$};

\draw (185,27) node[anchor=north west][inner sep=0.75pt] {$\mathit{Y}'$};

\draw (77,85) node[anchor=north west][inner sep=0.75pt] [align=left] {$\displaystyle \tau_{\mathrm{free}}$};

\draw (271,84) node[anchor=north west][inner sep=0.75pt] {$=$};

\draw (222,120) node[anchor=north west][inner sep=0.75pt] [align=left] {$\displaystyle_{\tau_\mathrm{free}}$};

\draw (342,95) node[anchor=north west][inner sep=0.75pt] {$X$};

\draw (342,70) node[anchor=north west][inner sep=0.75pt] {$\mathit{Y}$};

\draw (337,135) node[anchor=north west][inner sep=0.75pt] {$X'$};

\draw (337,26) node[anchor=north west][inner sep=0.75pt] {$\mathit{Y}'$};

\draw (345,52) node[anchor=north west][inner sep=0.75pt] [align=left] {$\displaystyle_{f_\mathrm{post}}$};

\draw (345,115) node[anchor=north west][inner sep=0.75pt] [align=left] {$\displaystyle_{f_\mathrm{{pre}}}$};

\end{tikzpicture}.
\end{align}
The action of the deterministic comb given in Eq.~\eqref{fig:f_free} on a function $f$  is the sequential composition of the pre-processing $f_{\mathrm{pre}}$, the function $f$, and the post-processing $ f_{\mathrm{post}}$: $\tau_{\mathrm{free}}(f)= f_{\mathrm{post}} \circ f \circ f_{\mathrm{pre}}$. We thus obtain a new function that can be depicted as
\begin{align}\label{fig:f_seq}
\tikzset{every picture/.style={line width=0.75pt}}
\begin{tikzpicture}[x=0.75pt,y=0.75pt,yscale=-1,xscale=1]

\draw[line width=0.75] (84,39) -- (84,64);

\draw[line width=0.75] (57,64) rectangle (110,93);

\draw[line width=0.75] (84,93) -- (84,118);

\draw[line width=0.75] (215,39) -- (215,64);

\draw[line width=0.75] (215,93) -- (215,118);

\draw[line width=0.75] (166,64) rectangle (265,93);

\draw (127,78) node[anchor=north west][inner sep=0.75pt] {$=$};

\draw (63,100) node[anchor=north west][inner sep=0.75pt] {$X'$};

\draw (63,44) node[anchor=north west][inner sep=0.75pt] {$Y'$};

\draw (60,70) node[anchor=north west][inner sep=0.75pt] {$\tau_{\mathrm{free}}(f)$};

\draw (169,71) node[anchor=north west][inner sep=0.75pt] {$f_{\mathrm{post}} \circ f \circ f_{\mathrm{pre}}$};

\draw (195,100) node[anchor=north west][inner sep=0.75pt] {$X'$};

\draw (195,44) node[anchor=north west][inner sep=0.75pt] {$Y'$};

\end{tikzpicture}.
\end{align}

Since we have defined freeness 
in a type-independent manner, the sets of free processes of different types will automatically be consistent with each other. 
For pedagogical purposes, it is useful to see this consistency in the case of deterministic combs and functions: 
the free functions are exactly the ones that can be obtained from the 
 free deterministic combs
(the RHS of Eq.~\eqref{fig:f_free}) by taking $X$ and $Y$ to be trivial variables, that is, single-valued.
There is only one possible pre-processing function from $X'$ to the singleton set $X$---the one that maps every possible value of $X'$ to the unique element of the singleton set (which is the discard map). The set of post-processing functions from the singleton set $Y$ to $Y'$ is isomorphic to $Y'$; that is, there are $|Y|$ possible choices of post-processing functions: those that map the single element into $y \in Y$ for each possible value that the variable $Y$ can take.   These are what we described earlier as the ``reset to $y'$'' functions, denoted $\mathtt{R}_{y'}$.

A central problem in any resource theory is to determine the partial order of resources, where one resource is above another if the first can be converted to the second by free operations.  In this paper, the objects that we will conceptualize as resources are the \emph{functions}, so our goal is to determine the partial order over functions in $\RT$.   
Because the most general operation from a function to a function is a deterministic comb, it is the set of  free deterministic combs that defines the set of free operations on functions. 

 The partial order of deterministic combs or of higher-order resources in $\RT$ can be studied, but it is outside of the scope of this paper.

\subsection{Example: bit-to-bit functions
}\label{se:fb2b}

In order to get some intuition for the Resource Theory of Causal Influence, let us discuss in detail the case where both $X$ and $Y$ are binary variables, with values in the set $\{0,1\}$. In this case, the \textit{set of all resources} has four elements, corresponding to the functions $\Id$, $\F$, $\Rzero$ and $\Rone$ defined in Sec.~\ref{sec_flag}. The functions $\Id$ and $\F$ can transmit causal influence, as their output depends on their input, while for the functions $\Rzero$ and $\Rone$ there is no causal influence from input to output. Thus the equivalence class of free resources is the set $\{\Rzero,\Rone\}$, while the nonfree resources are $\Id$ and $\F$, which also form an equivalence class because one can freely transform one into the other simply by post-processing with $\F$ or by pre-processing with $\F$.

Therefore, the resources in the bit-to-bit case constitute a total order: the equivalence class $\{\Id,\F\}$ strictly dominates the equivalence class $\{\Rzero,\Rone\}$. A useful monotone that completely characterizes the total order in this example can be defined using the size of the \textit{image} of the function (the cardinality of the set of all possible output values the function can produce from its domain). While $\Rzero$ and $\Rone$ have $|\Imm(\Rzero)|=|\Imm(\Rone)|=1$, $\Id$ and $\F$ have $|\Imm(\Id)|=|\Imm(\F)|=2$. We define the monotone as $M(f)=\log_2|\Imm(f)|$. In the resource theory of bit-to-bit functions, it can be interpreted as the number of bits of causal influence: zero bits for $\Rzero$ and $\Rone$ and one bit for $\Id$ and $\F$.

\subsection{Resource conversion assessment}
\label{sec_functionRT_resource_conversion}

We now turn to the generic case, where $X$ and $Y$ can have any finite cardinalities, and ask: given two functions  $f \, : \, X \to Y$ and  $g \, : \, X' \to Y'$, can one find a free operation $\tau_{\textrm{free}}$ that takes $f$ to $g$, i.e., such that $g = \tau_{\textrm{free}}(f)$? In the case of binary $X$ and $Y$, we saw that whether resource conversion is possible depends on the size of the image of $f$ and $g$. As we show next, this holds in general.

\begin{prop}
	\label{prop_increase_imagesize}
	Let $f \, : \, X \to Y$ and $g \, : \, X' \to Y'$ be two functions whose domains $X$, $X'$  and co-domains $Y$, $Y'$ are of finite cardinality. Then, $f \to g$ in $\RT$ if and only if $|\Imm(f)|\geq|\Imm(g)|$. 

    Therefore, the partial order over equivalence classes in $\RT$ is a total order.
\end{prop}
\begin{proof}

\textbf{(If side):} For this part of the proof, assume that $|\Imm(f)|\geq|\Imm(g)|$.
    Any function $f \, : \, X \to Y$ can be written as $f=f_i \circ f_s$, where $f_s: X \rightarrow \Imm(f)$ is a surjective function and $f_i: \Imm(f)\rightarrow Y$ is an injective function, which can be depicted diagrammatically as
    \begin{align}\label{fig:f_composition}
    \tikzset{every picture/.style={line width=0.75pt}}
\begin{tikzpicture}[x=0.75pt,y=0.75pt,yscale=-1,xscale=1]

\draw[line width=0.75] (55,81) rectangle (85,101);

\draw[line width=0.75] (70,52) -- (70,81);

\draw[line width=0.75] (70,101) -- (70,131);

\draw[line width=0.75] (218,107) rectangle (248,127);

\draw[line width=0.75] (233,78) -- (233,107);

\draw[line width=0.75] (233,127) -- (233,157);

\draw[line width=0.75] (218,57) rectangle (248,77);

\draw[line width=0.75] (233,28) -- (233,57);

\draw (64,83) node[anchor=north west][inner sep=0.75pt] {$f$};

\draw (48,110) node[anchor=north west][inner sep=0.75pt] {$X$};

\draw (49,60) node[anchor=north west][inner sep=0.75pt] {$Y$};

\draw (135,87) node[anchor=north west][inner sep=0.75pt] {$=$};

\draw (226,110) node[anchor=north west][inner sep=0.75pt] {$f_s$};

\draw (210,136) node[anchor=north west][inner sep=0.75pt] {$X$};

\draw (188,84) node[anchor=north west][inner sep=0.75pt] {$ \Imm(f)$};

\draw (227,60) node[anchor=north west][inner sep=0.75pt] {$f_i$};

\draw (212,37) node[anchor=north west][inner sep=0.75pt] {$Y$};

\end{tikzpicture}.
    \end{align} 
    It is easy to see that this is the case: we can choose the function $f_s$ to have the same action as $f$, mapping $x$ to $f(x)$; the function $f_i$ then embeds this result inside the co-domain of $f$. We will similarly decompose $g \, : \, X' \to Y'$ into $g=g_i \circ g_s$, with $g_s: X'\rightarrow \Imm(g)$ surjective and $g_i: \Imm(g)\rightarrow Y'$ injective.

	 Since $|\Imm(f)|\geq|\Imm(g)|$, it is possible to find an injective function 
     $i:\Imm(g)\rightarrow \Imm(f)$ and its left inverse $\pi:\Imm(f)\rightarrow \Imm(g)$ such that the composition $\pi\circ i$ is the identity function on $\Imm(g)$, that is, $\pi(i(z))=z \,\,\forall z\in \Imm(g)$. That is, one can reversibly map a smaller set into a bigger set.
    We will now show that it is possible to transform $f$ into $g$ by pre- and post-processings with the functions defined above. All surjective functions have a right inverse, that is, there exists $f^{-1}_s$ such that $f_s \circ f^{-1}_s = \Id$ for $f_s$ surjective. All injective functions have a left inverse, that is, there exists $f^{-1}_i$ such that $ f^{-1}_i \circ f_i = \Id$ for $f_i$ injective. A protocol that transforms $f$ into $g$ is 
    the following:
\begin{align}\label{fig:protocol}
\tikzset{every picture/.style={line width=0.75pt}}
\begin{tikzpicture}[x=0.75pt,y=0.75pt,yscale=-1,xscale=1]

\draw[line width=0.75] (56,133) rectangle (85,153);

\draw[line width=0.75] (70,71) -- (70,91);

\draw[line width=0.75] (57,92) rectangle (87,112);

\draw[line width=0.75] (71,112) -- (71,133);

\draw[line width=0.75] (58,51) rectangle (88,71);

\draw[line width=0.75] (71,30) -- (71,51);

\draw[color={rgb, 255:red, 144; green, 19; blue, 254},draw opacity=1,fill={rgb, 255:red, 255; green, 255; blue, 255},fill opacity=1,line width=0.75] (56,174) rectangle (85,194);

\draw[color={rgb, 255:red, 144; green, 19; blue, 254},draw opacity=1,fill={rgb, 255:red, 255; green, 255; blue, 255},fill opacity=1,line width=0.75] (71,153) -- (71,174);

\draw[line width=0.75] (56,215) rectangle (85,235);

\draw[color={rgb, 255:red, 144; green, 19; blue, 254},draw opacity=1,fill={rgb, 255:red, 255; green, 255; blue, 255},fill opacity=1,line width=0.75] (71,194) -- (71,215);

\draw[line width=0.75] (56,256) rectangle (85,277);

\draw[line width=0.75] (70,235) -- (70,256);

\draw[line width=0.75] (56,297) rectangle (85,318);

\draw[line width=0.75] (70,277) -- (70,297);

\draw[line width=0.75] (70,318) -- (70,338);

\draw[line width=0.75] (165,113) rectangle (194,133);

\draw[line width=0.75] (166,72) rectangle (196,92);

\draw[line width=0.75] (167,31) rectangle (197,52);

\draw[line width=0.75] (180,10) -- (180,31);

\draw[line width=0.75] (165,154) rectangle (194,174);

\draw[line width=0.75] (165,195) rectangle (194,216);

\draw[line width=0.75] (180,174) -- (180,195);

\draw[line width=0.75] (165,236) rectangle (194,257);

\draw[line width=0.75] (165,278) rectangle (194,298);

\draw[line width=0.75] (179,298) -- (179,319);

\draw[line width=0.75] (181,52) -- (181,72);

\draw[line width=0.75] (181,92) -- (181,113);

\draw[line width=0.75] (180,133) -- (180,154);

\draw[line width=0.75] (180,216) -- (180,236);

\draw[line width=0.75] (180,257) -- (180,278);

\draw[line width=0.75] (165,319) rectangle (195,340);

\draw[line width=0.75] (180,340) -- (180,360);

\draw (206,261) .. controls (210,261) and (213,258) .. (213,254) -- (213,196) .. controls (213,189) and (215,185) .. (220,185) .. controls (215,185) and (213,182) .. (213,176) -- (213,118) .. controls (213,113) and (211,111) .. (206,111);

\draw[line width=0.75] (330,159) rectangle (360,180);

\draw[line width=0.75] (331,119) rectangle (361,139);

\draw[line width=0.75] (344,98) -- (344,119);

\draw[line width=0.75] (331,200) rectangle (360,221);

\draw[line width=0.75] (331,242) rectangle (361,262);

\draw[line width=0.75] (346,139) -- (346,159);

\draw[line width=0.75] (345,180) -- (345,200);

\draw[line width=0.75] (346,221) -- (346,242);

\draw[line width=0.75] (346,262) -- (346,283);

\draw (370,221) .. controls (375,221) and (377,218) .. (377,213) -- (377,200) .. controls (377,194) and (380,190) .. (384,190) .. controls (380,190) and (377,187) .. (377,180) -- (377,167) .. controls (377,162) and (375,160) .. (370,160);

\draw
[color={rgb, 255:red, 144; green, 19; blue, 254},opacity=1] 
(492,175) rectangle (522,195);

\draw
[color={rgb, 255:red, 144; green, 19; blue, 254},opacity=1] 
(507,195) -- (507,216);

\draw
[color={rgb, 255:red, 144; green, 19; blue, 254},opacity=1] 
(507,154) -- (507,175);

\draw (60,132) node[anchor=north west][inner sep=0.75pt] {$f_i^{-1}$};

\draw (49,199) node[anchor=north west][inner sep=0.75pt] 
{$X$};

\draw (50,33) node[anchor=north west][inner sep=0.75pt] {$Y'$};

\draw (112,179) node[anchor=north west][inner sep=0.75pt] {$=$};

\draw (67,99) node[anchor=north west][inner sep=0.75pt] {$\pi$};

\draw (67,57) node[anchor=north west][inner sep=0.75pt] {$g_i$};

\draw (64,177) node[anchor=north west][inner sep=0.75pt] [color={rgb, 255:red, 144; green, 19; blue, 254},opacity=1] 
{$f$};

\draw (59,216) node[anchor=north west][inner sep=0.75pt] {$f_s^{-1}$};

\draw (66,260) node[anchor=north west][inner sep=0.75pt] {$i$};

\draw (64,303) node[anchor=north west][inner sep=0.75pt] {$g_s$};

\draw (30,74) node[anchor=north west][inner sep=0.75pt] {$\Imm(g)$};

\draw (30,115) node[anchor=north west][inner sep=0.75pt] {$\Imm(f)$};

\draw (49,157) node[anchor=north west][inner sep=0.75pt] 
{$Y$};

\draw (30,239) node[anchor=north west][inner sep=0.75pt] {$\Imm(f)$};

\draw (30,280) node[anchor=north west][inner sep=0.75pt] {$\Imm(g)$};

\draw (49,324) node[anchor=north west][inner sep=0.75pt] {$X'$};

\draw (159,13) node[anchor=north west][inner sep=0.75pt] {$Y'$};

\draw (175,78) node[anchor=north west][inner sep=0.75pt] {$\pi$};

\draw (173,158) node[anchor=north west][inner sep=0.75pt] {$f_i$};

\draw (173,199) node[anchor=north west][inner sep=0.75pt] {$f_s$};

\draw (174,281) node[anchor=north west][inner sep=0.75pt] {$i$};

\draw (140,54) node[anchor=north west][inner sep=0.75pt] {$\Imm(g)$};

\draw (140,95) node[anchor=north west][inner sep=0.75pt] {$\Imm(f)$};

\draw (158,138) node[anchor=north west][inner sep=0.75pt] {$Y$};

\draw (158,220) node[anchor=north west][inner sep=0.75pt] {$X$};

\draw (140,301) node[anchor=north west][inner sep=0.75pt] {$\Imm(g)$};

\draw (159,345) node[anchor=north west][inner sep=0.75pt] {$X'$};

\draw (175,326) node[anchor=north west][inner sep=0.75pt] {$g_s$};

\draw (170,237) node[anchor=north west][inner sep=0.75pt] {$f_s^{-1}$};

\draw (169,113) node[anchor=north west][inner sep=0.75pt] {$f_i^{-1}$};

\draw (176,38) node[anchor=north west][inner sep=0.75pt] {$g_i$};

\draw (140,178) node[anchor=north west][inner sep=0.75pt] {$\Imm(f)$};

\draw (140,260) node[anchor=north west][inner sep=0.75pt] {$\Imm(f)$};

\draw (229,178) node[anchor=north west][inner sep=0.75pt] {$\Id$};

\draw (272,179) node[anchor=north west][inner sep=0.75pt] {$=$};

\draw (323,101) node[anchor=north west][inner sep=0.75pt] {$Y'$};

\draw (339,165) node[anchor=north west][inner sep=0.75pt] {$\pi$};

\draw (340,204) node[anchor=north west][inner sep=0.75pt] {$i$};

\draw (305,141) node[anchor=north west][inner sep=0.75pt] {$\Imm(g)$};

\draw (305,183) node[anchor=north west][inner sep=0.75pt] {$\Imm(f)$};

\draw (305,224) node[anchor=north west][inner sep=0.75pt] {$\Imm(g)$};

\draw (325,267) node[anchor=north west][inner sep=0.75pt] {$X'$};

\draw (341,248) node[anchor=north west][inner sep=0.75pt] {$g_s$};

\draw (340,125) node[anchor=north west][inner sep=0.75pt] {$g_i$};

\draw (393,181) node[anchor=north west][inner sep=0.75pt] {$\Id$};

\draw (434,179) node[anchor=north west][inner sep=0.75pt] {$=$};

\draw (487,154) node[anchor=north west][inner sep=0.75pt] 
{$Y'$};

\draw (485,201) node[anchor=north west][inner sep=0.75pt] 
{$X'$};

\draw (502,182) node[anchor=north west][inner sep=0.75pt][color={rgb, 255:red, 144; green, 19; blue, 254},opacity=1] 
{$g$};

\end{tikzpicture}.
\end{align}

\noindent
    \textbf{(Only if side)}   For this part of the proof, we assume that  $|\Imm(f)|<|\Imm(g)|$ and show that $f\not\rightarrow g$.
    
    The size of the image of a function is at most the size of its domain. 
In the decomposition $f = f_i \circ f_s$, the domain of $f_i$ is necessarily $\Imm(f)$.  
Therefore, any pre- and post-processing $\tau=(\tau_{\text{pre}},\tau_{\text{post}})$ yields a composition whose effective domain cannot exceed $\Imm(f)$ and can never increase the size of the image of the original function. Hence,

	\begin{equation}
		|\Imm(\tau_{\text{post}}\circ f \circ \tau_{\text{pre}})| \leq |\Imm(f)|\, \quad\forall \tau_{\text{pre}},\tau_{\text{post}}.
	\end{equation}
    If there existed $\tau_{\text{pre}}$ and $\tau_{\text{post}}$ such that $g = \tau_{\text{post}}\circ f \circ \tau_{\text{pre}}$, we could conclude that $|\Imm(f)| \geq |\Imm(g)|$, which is a contradiction with the assumption that $ |\Imm(f)|< |\Imm(g)|$. Hence, $f$ cannot be transformed into $g$ via free operations.

\end{proof}
It follows that this resource theory is perfectly captured by a single monotone 
\begin{equation}
    M(f)=\log_2|\Imm(f)|,
\end{equation}
and that the partial order of functions is a total order, no matter the cardinalities of their inputs and outputs.

\section{The Resource Theory of Knowledge of Causal Influence}\label{se:rtfd}

Motivated by the considerations in the introduction, we now transition from considering functions as the resources of interest, to considering {\em probability distributions over functions} as the resources of interest. We refer to the resulting resource theory as the \textit{Resource Theory of Knowledge of Causal Influence} and denote it by $\RTKnow$. 
(We discuss the relationship between this resource theory and $\RT$ in Section~\ref{sec_relation_RTs}.)
Similarly to Sec. \ref{se:rtf}, we will focus on the simplest case of two observed variables, $X$ and $Y$, where the causal influence goes from $X$ to $Y$.  However, we now imagine that there is an unobserved error variable that influences $Y$, as in Figs.~\ref{fig:c} and~\ref{fig:d}.

\subsection{Resources}

In contrast to $\RT$,  the objects  of interest now do not consist of deterministic processes (functions, deterministic combs, etc), but instead consist of \emph{probability distributions over deterministic processes}. 

The first type of object, which is the one that our work will focus on in the partial order characterization, consists of the probability distributions over the space of possible functions from a finite cardinality random variable $X$ to a finite cardinality random variable $Y$. Mathematically, it is a probability distribution over a random variable $F$ with values in the finite set $F(X\to Y)$ of all possible functions from $X$ to $Y$.  We denote such a distribution by $\PP_F$, and the set of all probability distributions over functions from $X$ to $Y$ is denoted by $\mathcal{D}(F(X\to Y))$. Recalling the notation $[f]=\delta_{F,f}$, we represent a probability distributions over functions as
\begin{align}
   \PP_F\,=\sum_{f\in F(X\to Y)} \PP_F(f)[f], 
\end{align}
with $\PP_{F}(f)$ being a weight coefficient associated to each function and $\sum_{f\in F(X\to Y)} \PP_{F}(f) =1$.

The sequential composition of probability distributions over functions  can be 
well defined
such that it 
is also a probability distribution over functions:
\begin{align}
\PP_{F}\circ\PP'_{F'}=\left(\sum_{f\in F(Y\to Z)}\PP_{F}(f)[f]\right)\circ \left(\sum_{f'\in F(X\to Y)}\PP'_{F'}(f')[f']\right):= \sum_{f}\sum_{f'}\PP_{F}(f) \PP'_{F'}(f')[f\circ f'] .
\end{align}

Diagrammatically, a probability distribution over functions will be represented as: 
\begin{align}\label{fig:Pfdef}
\tikzset{every picture/.style={line width=0.75pt}}

\begin{tikzpicture}[x=0.75pt,y=0.75pt,yscale=-1,xscale=1]

\draw [line width=0.75]    (70,45) -- (70,75) ;
\draw [line width=0.75]    (70,95) -- (70,125) ;
\draw  [line width=0.75]  (60,75) -- (80,75) -- (80,95) -- (60,95) -- cycle ;

\draw  [line width=0.75]  (50,85) -- (60,85) ;

\draw  [line width=0.75]  (20,85) -- (50,65) -- (50,105) -- cycle ;
\draw (29,79) node [anchor=north west][inner sep=0.75pt]    {$\mathbb{P}_{F}$};

\draw (48,112) node [anchor=north west][inner sep=0.75pt]    {$X$};
\draw (52,48) node [anchor=north west][inner sep=0.75pt]    {$Y$};

\end{tikzpicture}
\end{align}

Here, we are using the notation from the framework of causal-inferential theories~\cite{omlet}, where vertical wires are part of the causal theory and horizontal wires are part of the inferential theory. (That is, a process with vertical inputs and outputs is a function representing some actual dynamics happening in the world; the horizontal wire in the diagram is a variable ranging over all such possible functions, and the triangle represents an agent's knowledge of what the function is.)

Note that this can also be represented as in Eq.~\eqref{fig:prob_over_functions}, where each value of $\Lambda$ is associated with one deterministic function from $X$ to $Y$.  The process from $X \times \Lambda$ to $Y$ may be termed a \emph{universal function applier}: the $\Lambda$ input acts as a program that ranges over all possible functions from $X$ to $Y$ (the decoration on the box signifies which input is the program). In other words, the output $Y$ for a given input $X$ and $\Lambda$ of the universal function applier is the result of applying the deterministic function associated with $\Lambda$ to $X$.

\begin{align}\label{fig:prob_over_functions}
\tikzset{every picture/.style={line width=0.75pt}} 

\begin{tikzpicture}[x=0.75pt,y=0.75pt,yscale=-1,xscale=1]

\draw   (140.73,73.6) -- (200,73.6) -- (200,98) -- (140.73,98) -- cycle ;
\draw    (170.4,55.8) -- (170.4,73.47) ;
\draw    (190.4,97.8) -- (190.4,125) ;
\draw    (159.7,97.5) -- (159.7,152.5) ;
\draw   (159.5,164.67) -- (148,152.33) -- (171,152.33) -- cycle ;
\draw    (141,159.33) -- (154.33,159.33) ;
\draw   (107.33,159.33) -- (140.33,139.33) -- (140.33,179.33) -- cycle ;
\draw   (154,91) -- (165.33,91) -- (165.33,97.67) -- (154,97.67) -- cycle ;

\draw (154.73,55.67) node [anchor=north west][inner sep=0.75pt]    {$Y$};
\draw (174.07,105.33) node [anchor=north west][inner sep=0.75pt]    {$X$};
\draw (145.33,105.33) node [anchor=north west][inner sep=0.75pt]    {$\Lambda $};
\draw (119.33,154.33) node [anchor=north west][inner sep=0.75pt]    {$\mathbb{P}_{\Lambda }$};

\end{tikzpicture}
\end{align}

The stochastic map from $X$ to $Y$ associated to a given probability distribution over functions is given by $\PP_{Y|X}= \sum_{f\in F(X\to Y)} \PP_{F}(f)\,\delta_{Y,f(X)}$.

The second type of object of interest consists of \emph{probability distributions over deterministic combs}, represented by $\PP_\tau$.  These are represented as

\begin{align}    \label{fig:prob_over_combs}
    \tikzset{every picture/.style={line width=0.75pt}} 

\begin{tikzpicture}[x=0.75pt,y=0.75pt,yscale=-1,xscale=1]

\draw   (180.73,100.27) -- (240,100.27) -- (240,124.67) -- (180.73,124.67) -- cycle ;
\draw    (210.4,69.47) -- (210.4,100.13) ;
\draw    (230.4,124.47) -- (230.4,146.67) ;
\draw    (190.7,124.17) -- (190.7,179.17) ;
\draw    (150,154.8) -- (163.33,154.8) ;
\draw   (116.33,154.8) -- (149.33,134.8) -- (149.33,174.8) -- cycle ;
\draw   (180.73,179.27) -- (240,179.27) -- (240,203.67) -- (180.73,203.67) -- cycle ;
\draw    (210.4,208.47) -- (210.4,228.4) -- (210.4,230.67) ;
\draw    (230.4,157.47) -- (230.4,179.67) ;
\draw  [dash pattern={on 4.5pt off 4.5pt}][line width=0.75]  (259.4,214.6) -- (163.8,214.6) -- (163.8,91.4) ;
\draw  [dash pattern={on 4.5pt off 4.5pt}][line width=0.75]  (200.2,172.6) -- (259.4,172.6) -- (259.4,214.6) ;
\draw  [dash pattern={on 4.5pt off 4.5pt}][line width=0.75]  (163.8,91) -- (260.2,91) -- (260.2,137.4) ;
\draw  [dash pattern={on 4.5pt off 4.5pt}][line width=0.75]  (200.2,172.6) -- (200.2,135.57) -- (260.6,135.57) ;

\draw (190.73,69.33) node [anchor=north west][inner sep=0.75pt]    {$Y'$};
\draw (213.07,136) node [anchor=north west][inner sep=0.75pt]    {$Y$};
\draw (128.33,149.8) node [anchor=north west][inner sep=0.75pt]    {$\mathbb{P}_{\tau }$};
\draw (189.73,215.2) node [anchor=north west][inner sep=0.75pt]    {$X'$};
\draw (213.07,157) node [anchor=north west][inner sep=0.75pt]    {$X$};

\end{tikzpicture}
\end{align}

They can be equivalently represented using the universal function applier as in Eq.~\eqref{fig:direct_cause_comb}. Here, the upper  universal function applier is such that $\Lambda'$ ranges over functions from $Y$ to $Y'$, and the lower universal function applier is such that $\Lambda$ ranges over functions from $X'$ to $X\times\Lambda'$.

\begin{align}
\label{fig:direct_cause_comb} \tikzset{every picture/.style={line width=0.75pt}} 

\begin{tikzpicture}[x=0.75pt,y=0.75pt,yscale=-1,xscale=1]

\draw   (160.73,93.6) -- (220,93.6) -- (220,118) -- (160.73,118) -- cycle ;
\draw    (190.4,75.8) -- (190.4,93.47) ;
\draw    (210.4,117.8) -- (210.4,140) ;
\draw    (179.7,117.5) -- (179.7,172.5) ;
\draw   (179.5,256.67) -- (168,244.33) -- (191,244.33) -- cycle ;
\draw    (161,251.33) -- (174.33,251.33) ;
\draw   (127.33,251.33) -- (160.33,231.33) -- (160.33,271.33) -- cycle ;
\draw   (174,111) -- (185.33,111) -- (185.33,117.67) -- (174,117.67) -- cycle ;
\draw   (160.73,172.6) -- (220,172.6) -- (220,197) -- (160.73,197) -- cycle ;
\draw    (210.4,196.8) -- (210.4,219) ;
\draw    (179.7,196.5) -- (179.7,244.5) ;
\draw   (174,190) -- (185.33,190) -- (185.33,196.67) -- (174,196.67) -- cycle ;
\draw    (210.4,150.8) -- (210.4,173) ;

\draw (171.73,75.67) node [anchor=north west][inner sep=0.75pt]    {$Y'$};
\draw (193.07,125.33) node [anchor=north west][inner sep=0.75pt]    {$Y$};
\draw (163.33,125.33) node [anchor=north west][inner sep=0.75pt]    {$\Lambda '$};
\draw (139.33,246.33) node [anchor=north west][inner sep=0.75pt]    {$\mathbb{P}_{\Lambda }$};
\draw (191.07,204.33) node [anchor=north west][inner sep=0.75pt]    {$X'$};
\draw (164.33,204.33) node [anchor=north west][inner sep=0.75pt]    {$\Lambda $};
\draw (193.07,158.33) node [anchor=north west][inner sep=0.75pt]    {$X$};

\end{tikzpicture}
\end{align}

In the next subsection, we will define freeness of a deterministic process in our resource theory in a type-independent manner, just like we did for $\RT$.  In this paper, the objects that we will conceptualize as resources are the \emph{probability distributions over functions}, so our goal is to determine the partial order over these in $\RTKnow$.  The free set of {\em probability distributions over deterministic combs} (which is determined by specializing our type-independent definition of freeness to deterministic combs)  then defines the set of free operations that governs interconversion among the resources.

The partial order of probability distributions over deterministic combs or of higher-order resources in $\RTKnow$ can also be studied, but it is outside of the scope of this paper.

\subsection{ Free resources and free operations}\label{se:skft}

Regardless of what type of process one considers (knowledge of function, knowledge of deterministic comb, etcetera), the definition of the free subset of processes is the same: it is those for which there is no cause-effect connection across the bipartition. This definition can be applied in a type-independent manner, as we will exemplify now.

First, consider the type consisting of probability distributions over functions. The free resources are the ones where:  
\begin{align}\label{fig:Pf_reset}
\tikzset{every picture/.style={line width=0.75pt}}

\begin{tikzpicture}[x=0.75pt,y=0.75pt,yscale=-1,xscale=1]

\draw [line width=0.75]    (70,45) -- (70,75) ;
\draw [line width=0.75]    (70,95) -- (70,125) ;
\draw  [line width=0.75]  (60,75) -- (80,75) -- (80,95) -- (60,95) -- cycle ;

\draw  [line width=0.75]  (50,85) -- (60,85) ;
\draw  [line width=0.75]  (20,85) -- (50,65) -- (50,105) -- cycle ;
\draw (29,79) node [anchor=north west][inner sep=0.75pt]    {$\mathbb{P}_{F}$};

\draw [line width=0.75]    (161,100) -- (180,100) ;
\draw [line width=0.75]    (164,98) -- (177,98) ;
\draw [line width=0.75]    (167,96) -- (173,96) ;

\draw [line width=0.75]    (170,45) -- (170,75) ;
\draw [line width=0.75]    (170,95) -- (170,125) ;
\draw   (170,92) -- (155,75) -- (185,75) -- cycle ;

\draw  [line width=0.75]  (140,85) -- (150,85) ;
\draw  [line width=0.75]  (110,85) -- (140,65) -- (140,105) -- cycle ;
\draw (119,77) node [anchor=north west][inner sep=0.75pt]    {$\mathbb{P}_{F}$};

\draw  [line width=0.75, dashed]  (150,70) -- (190,70) -- (190,104) -- (150,104) -- cycle ;

\draw (48,112) node [anchor=north west][inner sep=0.75pt]    {$X$};
\draw (52,50) node [anchor=north west][inner sep=0.75pt]    {$Y$};
\draw (151,113) node [anchor=north west][inner sep=0.75pt]    {$X$};
\draw (175,48) node [anchor=north west][inner sep=0.75pt]    {$Y$};
\draw (91,80) node [anchor=north west][inner sep=0.75pt]    {$=$};

\end{tikzpicture}.
\end{align}

That is, a probability distribution over functions $\PP_F$ is considered free if it is of the form $\PP_{F}  = \sum_{y} \PP_{F}(\mathtt{R_y}) [\mathtt{R_y}]$. In other words, a free probability distribution over functions is given by a distribution that only has weight on the causally disconnected `reset to $y$' functions.

Now, consider the type consisting of probability distributions over deterministic combs. The free operations of this type are again those that do not have cause-effect connection across the bipartition, that is, those of the form
\begin{align}
    \label{fig:prob_over_combs_free}
    \tikzset{every picture/.style={line width=0.75pt}} 

\begin{tikzpicture}[x=0.75pt,y=0.75pt,yscale=-1,xscale=1]

\draw   (180.73,100.27) -- (240,100.27) -- (240,124.67) -- (180.73,124.67) -- cycle ;
\draw    (210.4,69.47) -- (210.4,100.13) ;
\draw    (230.4,124.47) -- (230.4,146.67) ;

\draw    (150,154.8) -- (163.33,154.8) ;
\draw   (116.33,154.8) -- (149.33,134.8) -- (149.33,174.8) -- cycle ;
\draw   (180.73,179.27) -- (240,179.27) -- (240,203.67) -- (180.73,203.67) -- cycle ;
\draw    (210.4,208.47) -- (210.4,228.4) -- (210.4,230.67) ;
\draw    (230.4,157.47) -- (230.4,179.67) ;
\draw  [dash pattern={on 4.5pt off 4.5pt}][line width=0.75]  (259.4,214.6) -- (163.8,214.6) -- (163.8,91.4) ;
\draw  [dash pattern={on 4.5pt off 4.5pt}][line width=0.75]  (200.2,172.6) -- (259.4,172.6) -- (259.4,214.6) ;
\draw  [dash pattern={on 4.5pt off 4.5pt}][line width=0.75]  (163.8,91) -- (260.2,91) -- (260.2,137.4) ;
\draw  [dash pattern={on 4.5pt off 4.5pt}][line width=0.75]  (200.2,172.6) -- (200.2,135.57) -- (260.6,135.57) ;

\draw (190.73,69.33) node [anchor=north west][inner sep=0.75pt]    {$Y'$};
\draw (213.07,136) node [anchor=north west][inner sep=0.75pt]    {$Y$};
\draw (128.33,149.8) node [anchor=north west][inner sep=0.75pt]    {$\mathbb{P}_{\tau }$};
\draw (189.73,215.2) node [anchor=north west][inner sep=0.75pt]    {$X'$};
\draw (213.07,157) node [anchor=north west][inner sep=0.75pt]    {$X$};

\end{tikzpicture}
\end{align}

Those can be equivalently represented by Eq.~\eqref{fig:common_cause_comb}, where it is more evident that the causal structure of the comb consists only of a common-cause. Here, the upper  universal function applier is such that $\Lambda_2$ ranges over functions from $Y$ to $Y'$, and the lower universal function applier is such that $\Lambda_1$ ranges over functions from $X'$ to $X$.

\begin{align} \label{fig:common_cause_comb}
\tikzset{every picture/.style={line width=0.75pt}} 

\begin{tikzpicture}[x=0.75pt,y=0.75pt,yscale=-1,xscale=1]

\draw   (180.73,107.27) -- (240,107.27) -- (240,131.67) -- (180.73,131.67) -- cycle ;
\draw    (210.4,89.47) -- (210.4,107.13) ;
\draw    (230.4,131.47) -- (230.4,153.67) ;
\draw   (174.5,271.33) -- (163,259) -- (186,259) -- cycle ;
\draw    (156,266) -- (169.33,266) ;

\draw   (95,266) -- (155.33,240) -- (155.33,292) -- cycle ;

\draw   (194,124.67) -- (205.33,124.67) -- (205.33,131.33) -- (194,131.33) -- cycle ;
\draw   (180.73,186.27) -- (240,186.27) -- (240,210.67) -- (180.73,210.67) -- cycle ;
\draw    (230.4,210.47) -- (230.4,232.67) ;
\draw   (194,203.67) -- (205.33,203.67) -- (205.33,210.33) -- (194,210.33) -- cycle ;
\draw    (230.4,164.47) -- (230.4,186.67) ;

\draw (170.0,259)
  .. controls (165.0,150.0) and (190.0,170.0)
  .. (200.2,131.0);

\draw (179.0,259)
  .. controls (178.0,235.0) and (205.0,245.0)
  .. (201.0,210.2);

\draw (191.73,89.33) node [anchor=north west][inner sep=0.75pt]    {$Y'$};
\draw (213.07,139) node [anchor=north west][inner sep=0.75pt]    {$Y$};

\draw (166.73,139) node [anchor=north west][inner sep=0.75pt]    {$\Lambda_2$};
\draw (176.73,218) node [anchor=north west][inner sep=0.75pt]    {$\Lambda_1$};

\draw (211.07,218) node [anchor=north west][inner sep=0.75pt]    {$X'$};
\draw (213.07,172) node [anchor=north west][inner sep=0.75pt]    {$X$};

\draw (114,261) node [anchor=north west][inner sep=0.75pt]    {$\PP_{\Lambda_1,\Lambda_2}$};

\end{tikzpicture}
\end{align}

Note that the free operations of this type enable the possibility of uncertainty about the pre- and post-processings, in particular, the possibility of \emph{correlated} pre- and post-processings. This stands in contrast to $\RT$ introduced in Sec.~\ref{se:rtf}, in which uncertainty is not introduced into the resource. 

Algebraically, we can write this free operation as
\begin{align}
\Gamma_{\textrm{free}}(\PP_F)
= \sum_{\substack{\lambda_1 \in \Lambda_1 \\  \lambda_2 \in \Lambda_2}} \PP_{\Lambda_1,\Lambda_2}(\lambda_1, \lambda_2)\,
[f^{\lambda_2}_{\mathrm{post}}]
\circ \PP_F
\circ [f^{\lambda_1}_{\mathrm{pre}}],
\end{align}
where each $f^{\lambda_1}_{\mathrm{pre}}$ is a function from $X'$ to $X$, each $f^{\lambda_2}_{\mathrm{post}}$ is a function from $Y$ to $Y'$,  and, as always, $[\cdot]$ indicates a point distribution. The values of $\lambda_1$ and $\lambda_2$ might be correlated, as indicated by the joint probability distribution $\PP_{\Lambda_1,\Lambda_2}$. 

Since we have defined free operations in a type independent manner, the sets of free processes of different types will automatically be consistent with each other.  That is, when we apply free operations on the trivial probability distribution over functions, we obtain the free probability distributions over functions. For pedagogical purposes, in Eq.~\eqref{fig:sanity_check} we explicitly show this consistency.

\begin{align}\label{fig:sanity_check}
    \tikzset{every picture/.style={line width=0.75pt}} 

\begin{tikzpicture}[x=0.75pt,y=0.75pt,yscale=-1,xscale=1]

\draw  [dash pattern={on 0.84pt off 2.51pt}]  (267.2,131.07) -- (267.2,142.27) ;
\draw    (256.47,128.27) -- (277.47,128.27) ;
\draw    (259.47,126.27) -- (274.13,126.27) ;
\draw    (261.8,124.6) -- (271.47,124.6) ;


\draw   (79.73,63.27) -- (139,63.27) -- (139,87.67) -- (79.73,87.67) -- cycle ;
\draw    (109.4,45.47) -- (109.4,63.13) ;
\draw  [dash pattern={on 0.84pt off 2.51pt}]  (129.4,87.47) -- (129.4,109.67) ;

\draw   (73.5,227.33) -- (62,215) -- (85,215) -- cycle ;
\draw    (55,222) -- (68.33,222) ;

\draw   (-6,222) -- (54.33,196) -- (54.33,248) -- cycle ;

\draw   (93,80.67) -- (104.33,80.67) -- (104.33,87.33) -- (93,87.33) -- cycle ;

\draw   (79.73,142.27) -- (139,142.27) -- (139,166.67) -- (79.73,166.67) -- cycle ;
\draw    (129.4,166.47) -- (129.4,188.67) ;
\draw   (93,159.67) -- (104.33,159.67) -- (104.33,166.33) -- (93,166.33) -- cycle ;
\draw  [dash pattern={on 0.84pt off 2.51pt}]  (129.4,120.47) -- (129.4,142.67) ;

\draw (69.0,215)
  .. controls (65.0,110.0) and (90.0,130.0)
  .. (99.2,87.0);

\draw (78.0,215)
  .. controls (77.0,195.0) and (104.0,205.0)
  .. (100.0,166.2);

\draw (150,130.67) node [anchor=north west][inner sep=0.75pt] {$=$};
\draw (90.73,45.33) node [anchor=north west][inner sep=0.75pt] {$Y'$};
\draw (112.07,95) node [anchor=north west][inner sep=0.75pt] {$Y$};
\draw (112.07,128) node [anchor=north west][inner sep=0.75pt] {$X$};
\draw (109.07,174) node [anchor=north west][inner sep=0.75pt] {$X'$};

\draw (70.73,95)  node [anchor=north west][inner sep=0.75pt] {$\Lambda_2$};
\draw (77.73,174) node [anchor=north west][inner sep=0.75pt] {$\Lambda_1$};

\draw (11,217) node [anchor=north west][inner sep=0.75pt] {$\PP_{\Lambda_1,\Lambda_2}$};


\draw   (216.93,63.67) -- (276.2,63.67) -- (276.2,88.07) -- (216.93,88.07) -- cycle ;
\draw    (246.6,45.87) -- (246.6,63.53) ;
\draw  [dash pattern={on 0.84pt off 2.51pt}]  (266.6,87.87) -- (266.6,110.07) ;

\draw   (210.7,227.73) -- (199.2,215.4) -- (222.2,215.4) -- cycle ;
\draw    (192.2,222.4) -- (205.53,222.4) ;

\draw   (131.2,222.4) -- (191.53,196.4) -- (191.53,248.4) -- cycle ;

\draw   (230.2,81.07) -- (241.53,81.07) -- (241.53,87.73) -- (230.2,87.73) -- cycle ;

\draw   (216.93,142.67) -- (276.2,142.67) -- (276.2,167.07) -- (216.93,167.07) -- cycle ;
\draw    (266.6,166.87) -- (266.6,189.07) ;
\draw   (230.2,160.07) -- (241.53,160.07) -- (241.53,166.73) -- (230.2,166.73) -- cycle ;

\draw (206.0,215.4)
  .. controls (202.0,110.0) and (227.0,130.0)
  .. (236.4,87.4);

\draw (215.0,215.4)
  .. controls (214.0,195.0) and (241.0,205.0)
  .. (237.2,166.6);

\draw (227.93,45.73) node [anchor=north west][inner sep=0.75pt] {$Y'$};
\draw (249.27,95.4) node [anchor=north west][inner sep=0.75pt] {$Y$};
\draw (249.27,128.4) node [anchor=north west][inner sep=0.75pt] {$X$};
\draw (246.27,174.4) node [anchor=north west][inner sep=0.75pt] {$X'$};

\draw (207.93,95.4)  node [anchor=north west][inner sep=0.75pt] {$\Lambda_2$};
\draw (214.93,174.4) node [anchor=north west][inner sep=0.75pt] {$\Lambda_1$};

\draw (150.0,217.4) node [anchor=north west][inner sep=0.75pt] {$\PP_{\Lambda_1,\Lambda_2}$};

\draw (288.6,130.67) node [anchor=north west][inner sep=0.75pt] {$=$};


\draw    (398.07,165.87) -- (419.07,165.87) ;
\draw    (401.07,163.87) -- (415.73,163.87) ;
\draw    (403.4,162.2) -- (413.07,162.2) ;

\draw   (358.53,63.27) -- (417.8,63.27) -- (417.8,87.67) -- (358.53,87.67) -- cycle ;
\draw    (388.2,45.47) -- (388.2,63.13) ;
\draw  [dash pattern={on 0.84pt off 2.51pt}]  (408.2,87.47) -- (408.2,109.67) ;

\draw   (352.3,227.33) -- (340.8,215) -- (363.8,215) -- cycle ;
\draw    (333.8,222) -- (347.13,222) ;

\draw   (272.8,222) -- (333.13,196) -- (333.13,248) -- cycle ;

\draw   (371.8,80.67) -- (383.13,80.67) -- (383.13,87.33) -- (371.8,87.33) -- cycle ;

\draw (348.0,215.0)
  .. controls (345.0,110.0) and (370.0,130.0)
  .. (378.0,87.0);


\draw (369.53,45.33) node [anchor=north west][inner sep=0.75pt] {$Y'$};
\draw (390.87,95) node [anchor=north west][inner sep=0.75pt] {$Y$};
\draw (387.87,174) node [anchor=north west][inner sep=0.75pt] {$X'$};
\draw    (408.2,166.47) -- (408.2,188.67) ;

\draw (349.53,95)  node [anchor=north west][inner sep=0.75pt] {$\Lambda_2$};

\draw (299.0,217.0) node [anchor=north west][inner sep=0.75pt] {$\PP_{\Lambda_2}$};

\draw (434.2,130.67) node [anchor=north west][inner sep=0.75pt] {$=$};


\draw   (522.3,163.73) -- (510.8,151.4) -- (533.8,151.4) -- cycle ;
\draw    (503.8,158.4) -- (517.13,158.4) ;

\draw   (442.8,158.4) -- (503.13,132.4) -- (503.13,184.4) -- cycle ;

\draw    (522.5,119.57) -- (522.5,151.57) ;

\draw (503.73,133.73) node [anchor=north west][inner sep=0.75pt] {$Y'$};
\draw (469.0,153.4) node [anchor=north west][inner sep=0.75pt] {$\PP_{\Lambda_2}$};

\end{tikzpicture}
\end{align}

\subsection{Resource conversion assessment: linear program} 

We now describe a linear program that can be used to determine the ordering relation that holds between any two resources. (This linear program is closely analogous to the one in Ref.~\cite{Wolfe2020quantifyingbell}.) In Sec.~\ref{se:b2brc}, we study the ordering relation for the bit-to-bit scenario in detail. 

We are given two resources, $\PP_{F}$ and $\PP'_{F'}$, and we ask whether $\PP_{F}$ can be transformed into $\PP'_{F'}$ by means of free operations. Formally, this is equivalent to asking whether $\PP'_{F'}$ is in the downward closure of $\PP_{F}$. Since the set of free operations forms a convex set with a finite number of extremal points, i.e., a polytope, and since free operations act linearly on resources, the downward closure of any given resource $\PP_{F}$ also forms a polytope. Hence, one only needs to check whether $\PP'_{F'}$ belongs to this polytope. 
Note that every vertex of the downward-closure polytope is necessarily an image of $\PP_{F}$ under an extremal free operation. Each extremal free operation is specified by a point distribution on a pair of functions---a pre-processing function from $X'$ to $X$ and a post-processing function from $Y$ to $Y'$. Consequently, we can check whether $\PP'_{F'}$ belongs to this polytope via three steps:
\begin{compactenum}
   \item List all the extremal free operations by listing all possible such pre- and post-processings.
   \item Find the image of $\PP_{F}$ under all the elements in the list of extremal free operations. 
   \item Check whether $\PP'_{F'}$ is in the convex hull of this image. 
\end{compactenum}
This last step is a well-known linear program, sometimes called the convex hull problem~\cite{Kalantari2015}.

We have hence specified a linear program that assesses resource conversion: 

\noindent \texttt{Program:} \textbf{Resource conversion in the resource theory of knowledge of causal influence}.\\
\texttt{Input}: two resources $\PP_{F}$ and $\PP'_{F'}$. \\
\texttt{Output}: yes/no answer to question of whether $\PP_{F} \, {\rightarrow} \, \PP'_{F'}$ in $\RTKnow$.

Finally, using two instances of this linear program one can check both whether $\PP_{F} \, \overset{\textrm{free}}{\rightarrow} \, \PP'_{F'}$ and whether $\PP'_{F'} \, \overset{\textrm{free}}{\rightarrow} \, \PP_{F}$, which fully specifies the exact ordering relation between the two resources.

\subsection{Bit-to-bit scenario: complete characterization}\label{sec:b2bc}

In this section we present the details of the Resource Theory of Knowledge of Causal Influence for the case of binary $X$ and $Y$, where $X$ is the parent of $Y$.

As mentioned before, we denote the set of all possible bit-to-bit functions as $B:=F(\text{bit}\to~\text{bit})$; in this case there are four possible functions from $X$ to $Y$: $\Id$, $\F$, $\Rzero$ and $\Rone$. Unlike Sec. \ref{se:fb2b}, here we do not restrict ourselves to cases where we have complete knowledge on the functional dependence: that is, our resource (probability distribution over functions) is in general of the form 
\begin{align}
   & \PP_{B} = \PP_B(\mathtt{\Id}) \, [\Id] + \PP_B(\mathtt{\F}) \, [\F] + \PP_B(\mathtt{\Rzero}) \, [\Rzero] + \PP_B(\mathtt{\Rone}) \, [\Rone]\, , \\
   & \text{where } \PP_B(\mathtt{\Id}) +\PP_B(\mathtt{\F}) +\PP_B(\mathtt{\Rzero}) +\PP_B(\mathtt{\Rone}) =1\, . \nonumber
\end{align} 
The space of probability distributions over any four-element set can be geometrically viewed as a simplex in three dimensions. It immediately follows that we can geometrically consider the possibilities for $\PP_{B}$ as a point in such a three-dimensional simplex. Such a geometrical view of the space of resources is discussed in more detail in Sec.~\ref{se:geom}. A useful parametrization of $\PP_{B}$ in terms of three parameters is given by 
\begin{align} \label{generalparametrization}
\PP_{B} = \beta\left(\frac{1-\alpha}{2}[\Id]+\frac{1+\alpha}{2}[\F]\right) +(1-\beta)\left(\frac{1-\gamma}{2}[\Rzero]+\frac{1+\gamma}{2}[\Rone]\right),
\end{align}
where $\beta \in [0,1]$, $\alpha \in [-1,1]$, and $\gamma \in [-1,1]$, and where $\gamma$ is not defined when $\beta=1$, and $\alpha$ is not defined when $\beta=0$. In this parametrization, $\beta$ represents the weight on the sector of functions that carry causal influence, while $\alpha$ represents bias towards the $\F$ function within this sector. Similarly, $\gamma$ represents bias towards $\Rone$ in the sector that does not carry causal influence was implemented. 
In this parametrization, one can readily see the set of free resources is the one corresponding to the line segment picked out by $\beta=0$, containing convex combinations of $[\Rzero]$ and $[\Rone]$.

\subsubsection{Resource conversion}\label{se:b2brc}

Our free operations are all combinations of a  probabilistic  pre- and a post-processing, where each is drawn from the set of four functions $\Id$, $\F$, $\Rzero$, and $\Rone$ (since it is sufficient to consider the case where $X'$ and $Y'$ are also binary variables) and where the choice of function may be correlated by a common cause. The set of free operations, denoted $\mathcal{T}$,  is then a polytope with 16 extremal points. 

As an example, imagine that the original resource is $\PP_B^4$ from Eq.~\eqref{eq_Flip_Rzero}, given by
\begin{equation}
    \PP_B^4=\frac{1}{3}[\F]+\frac{2}{3}[\Rzero],
\end{equation}
and we apply to it a free operation of the form of Eq.~\eqref{fig:common_cause_comb} where $\PP_{\Lambda_1,\Lambda_2}=\frac{1}{2}[00]+\frac{1}{2}[11]$, that is, $\Lambda_1=\Lambda_2=\Lambda$, and both the pre- and post- processing are equal to $\Id$ when $\Lambda=0$ and to $\F$ when $\Lambda=1$. In this case, one can check that $\PP_B^4$ gets transformed to 
\begin{equation}
    \PP_B^5=\frac{1}{3}[\F]+\frac{1}{3}[\Rzero]+\frac{1}{3}[\Rone].
\end{equation}

As discussed in Sec.~\ref{sec_RT_preliminaries}, a resource $\PP_{F}$ can be freely transformed to the set of resources that corresponds to its downward closure, i.e., the polytope defined by $\{\Gamma_{\textrm{free}}[\PP_{F}] \, | \, \Gamma_{\textrm{free}} \in \mathcal{T}\}$, here denoted by $\mathcal{P}_{\downarrow}$. We will now flesh out this polytope for the resources in the bit-to-bit case by studying the action of the extremal free operations on a generic resource $\PP_B$. In the next section we provide a geometrical depiction of the space of resources and their free operations. 

Let us consider a resource in the form of Eq.~\eqref{generalparametrization}, and see explicitly what happens when we act on it with different choices of the extremal free operations, i.e., with $\Gamma_{\textrm{pre}},\Gamma_{\textrm{post}} \in \{\Id,\F,\Rzero,\Rone\}$. 
\begin{itemize}
    \item[(i)] \textbf{$\Gamma_{\textrm{pre}}=\Gamma_{\textrm{post}}=\Id$}: \hfill \textit{(1  extremal free operation)}\\
    The image of a resource when one pre-processes and post-processes by the identity is, of course, itself.\\
    \phantom{A}  \hfill $\PP_{B} \rightarrow \PP_{B}$.
    
    \item[(ii)] \textbf{$\Gamma_{\textrm{pre}}=$anything, $\Gamma_{\textrm{post}}=\mathtt{R_y}$}: \hfill \textit{(8   extremal free operations)}\\
    If one post-processes by $\mathtt{R_y}$, then the image of $\PP_{B}$ is always the function $\mathtt{R_y}$, regardless of what the pre-processing is or what $\PP_{B}$ is.   
    Notice that even though there are 8  extremal free operations in this category, they can be classified into two groups of ``operationally equivalent operations'', namely, the one where $\Gamma_{\textrm{post}}=\mathtt{R_0}$ and the one where  $\Gamma_{\textrm{post}}=\mathtt{R_1}$,  in the sense that the output resource coincides for all the operations in the same group.   \\
    \phantom{A} \hfill $\PP_{B} \rightarrow [\mathtt{R_y}]$.
    
    \item[(iii)] \textbf{$\Gamma_{\textrm{pre}}=\mathtt{R_y}$, $\Gamma_{\textrm{post}}=\Id$}: \hfill \textit{(2 extremal free operations)}\\
    Here, the image of $\PP_{B}$ will be a probability distribution over the functions $\Rzero$ and $\Rone$. 
    Since the image of this operation on $\PP_{B}$ can be obtained as a mixture of the images of the operations from item (ii), it is not an extreme point of the downward closure $\mathcal{P}_{\downarrow}$ for any value of $y$.  \\
    \phantom{A} \hfill $\PP_{B} \rightarrow \frac{1}{2}(1-\gamma + \beta(\gamma-\alpha))[\Rzero]+\frac{1}{2}(1+\gamma - \beta(\gamma-\alpha))[\Rone]$.
    
    \item[(iv)] \textbf{$\Gamma_{\textrm{pre}}=\mathtt{R_y}$, $\Gamma_{\textrm{post}}=\F$}: \hfill \textit{(2 extremal free operations)}\\
    Here, the image of $\PP_{B}$ will be a probability distribution over the functions $\Rzero$ and $\Rone$ like in item (iii), but with the weights of  $\Rzero$ and $\Rone$ interchanged. Again, it is not an extreme point of $\mathcal{P}_{\downarrow}$ for any value of $y$. \\
    \phantom{A} \hfill $\PP_{B} \rightarrow \frac{1}{2}(1-\gamma + \beta(\gamma-\alpha))[\Rone]+\frac{1}{2}(1+\gamma - \beta(\gamma-\alpha))[\Rzero]$.
    
    \item[(v)] \textbf{$\Gamma_{\textrm{pre}}=\F$, $\Gamma_{\textrm{post}}=\Id$}: \hfill \textit{(1  extremal free operation )}\\
    In this case, the resource $\PP_{B}$ is transformed in a way that the sign of $\alpha$ is flipped, while $\beta$ and $\gamma$ remain invariant. 
    \\
    \phantom{A} \hfill $\PP_{B} \rightarrow \beta(\frac{1+\alpha}{2}[\Id]+\frac{1-\alpha}{2}[\F]) +(1-\beta)(\frac{1-\gamma}{2}[\Rzero]+\frac{1+\gamma}{2}[\Rone])$.
    
    \item[(vi)] \textbf{$\Gamma_{\textrm{pre}}=\Id$, $\Gamma_{\textrm{post}}=\F$}: \hfill \textit{(1 extremal free operation )}\\
    In this case, the sign of $\alpha$ is flipped {\em and} the sign of $\gamma$ is flipped, while $\beta$ is left invariant. 
    \\
    \phantom{A} \hfill $\PP_{B} \rightarrow \beta(\frac{1+\alpha}{2}[\Id]+\frac{1-\alpha}{2}[\F]) +(1-\beta)(\frac{1+\gamma}{2}[\Rzero]+\frac{1-\gamma}{2}[\Rone])$.
    
    \item[(vii)] \textbf{$\Gamma_{\textrm{pre}}=\F$, $\Gamma_{\textrm{post}}=\F$}: \hfill \textit{(1  extremal free operation )}\\
    In this case, the sign of $\gamma$ is flipped, while $\alpha$ and $\beta$ remain invariant. 
    \\
    \phantom{A} \hfill $\PP_{B} \rightarrow \beta(\frac{1-\alpha}{2}[\Id]+\frac{1+\alpha}{2}[\F]) +(1-\beta)(\frac{1+\gamma}{2}[\Rzero]+\frac{1-\gamma}{2}[\Rone])$.
\end{itemize}

We see then that the polytope $\mathcal{P}_{\downarrow}$, describing the downward closure of $\PP_B$, is given by the convex hull of (at most) the 6 resources  given in items (i), (ii), (v), (vi), and (vii). These are summarized in Table \ref{table1}. Note that even though the free operations in (iii) and (iv) are extremal, their images are not vertices of the polytope $\mathcal{P}_{\downarrow}$. 

\begin{table}
\centering
	\renewcommand{\arraystretch}{1.5}
    \small{
    \begin{adjustbox}{width=0.7\textwidth}
	\begin{tabular}{cccc}
		Item  & Pre & Post & Final Resource \\ \hline
        (i)   & $\Id$     & $\Id$      &  $\beta\left(\frac{1-\alpha}{2}[\Id]+\frac{1+\alpha}{2}[\F]\right)+(1-\beta)\left(\frac{1-\gamma}{2}[\Rzero]+\frac{1+\gamma}{2}[\Rone]\right)$ \\
        (ii)  & Any     & $\mathtt{R_y}$      & $[\mathtt{R_y}]$      \\       
        (v)  & $\F$     & $\Id$      & $\beta\left(\frac{1+\alpha}{2}[\Id]+\frac{1-\alpha}{2}[\F]\right)+(1-\beta)\left(\frac{1-\gamma}{2}[\Rzero]+\frac{1+\gamma}{2}[\Rone]\right)$     \\      
        (vi)   & $\Id$     & $\F$      & $\beta\left(\frac{1+\alpha}{2}[\Id]+\frac{1-\alpha}{2}[\F]\right)+(1-\beta)\left(\frac{1+\gamma}{2}[\Rzero]+\frac{1-\gamma}{2}[\Rone]\right)$ \\  
         (vii)    & $\F$     & $\F$      & $\beta(\frac{1-\alpha}{2}[\Id]+\frac{1+\alpha}{2}[\F]) +(1-\beta)(\frac{1+\gamma}{2}[\Rzero]+\frac{1-\gamma}{2}[\Rone])$   
	\end{tabular}
    \end{adjustbox}}
    \caption{List of 6 extremal points of the downward  closure polytope $\mathcal{P}_{\downarrow}$ of a resource $\PP_B=\beta\left(\frac{1-\alpha}{2}[\Id]+\frac{1+\alpha}{2}[\F]\right)+(1-\beta)\left(\frac{1-\gamma}{2}[\Rzero]+\frac{1+\gamma}{2}[\Rone]\right)$. Here Any $\in \{\Id, \F, \Rzero, \Rone \}$, and $\mathtt{R_y} \in \{\Rzero,\Rone\}$. Notice that although in principle Item (ii) has 8 transformations, they effectively reduce to two distinct transformations of the resource.  }
	\label{table1}
\end{table}

Since one can freely change the sign of $\alpha$ while leaving the other parameters unchanged (by the operation of item (v)), the sign of $\alpha$ is not relevant for determining the equivalence class of a resource. In addition, since the sign of $\gamma$ can be similarly changed, the sign of $\gamma$ is also irrelevant for determining the equivalence class. Therefore, 
only the values of $(|\alpha|, \beta, |\gamma|)$ are relevant to define the equivalence class of any given resource. 
However, recall that some resources do not have all such values defined. Specifically, when $\beta=1$ the resources do not have any specified value for $\gamma$, and when $\beta=0$ the resources do not have any specified value for $\alpha$. The case of $\beta=0$ singles out the set of free resources, so here the equivalence class can be identified without having to rely on any values for $\alpha$ or $\gamma$. To deal with the case $\beta=1$, note that free operations cannot increase the value of $\beta$. Therefore, if we start with a resource $\PP_{B}$ with $\beta=1$, all the resources that can potentially be in the same equivalence class as $\PP_{B}$ will also have $\beta=1$, and so also have $\gamma$ unspecified. Hence, when $\beta=1$, to decide whether a resource $\PP_{B}'$ is in the equivalence class of $\PP_{B}$ one would rely on the value of $|\alpha|$. Hence the fact that $\gamma$ is unspecified in this case creates no problems. We conclude then that the values of $(|\alpha|, \beta, |\gamma|)$ are sufficient to specify the equivalence classes of resources, even in the cases of $(|\alpha|, 1, \text{unspecified})$ and of $(\text{unspecified}, 0, |\gamma|)$.

In Sec.~\ref{se:b2bcsm} we will discuss in which sense  the set of parameters $(|\alpha|, \beta, |\gamma|)$ is also \emph{necessary} to specify the equivalence classe, so we can take $(|\alpha|, \beta, |\gamma|)$ to be the {\em canonical form} of a resource (see Lemma \ref{lem:canform}).

Let us finish this subsection by showing that the partial order over equivalence classes in this resource theory is \emph{not} a total order. 

\begin{prop}\label{prop:nottot}
    The partial order over equivalence classes in $\RTKnow$ is not a total order.
\end{prop}
\begin{proof}
    We will prove the claim by presenting an example of two resources in the bit-to-bit scenario that cannot be transformed into each other. Consider the following two resources:
    \begin{align}
        \PP_{B}^{1}=\frac{1}{2}[\Id]+\frac{1}{2}[\F] &  \qquad \leftrightarrow \qquad  \beta=1, \alpha = 0, \gamma \text{ undefined} \label{eq_F1}\\
        \PP_{B}^{2}=\frac{1}{2}[\Id]+\frac{1}{2}[\Rzero] & \qquad \leftrightarrow \qquad    \beta=1/2, \alpha = -1, \gamma=-1
    \label{eq_F2}
    \end{align}

    None of the extremal free operations can possibly increase the value of $\beta$, as can be seen in Table~\ref{table1}. Therefore, we have that $\PP_{B}^{2}\not\to \PP_{B}^{1}$ in the pre-order of this resource theory.   Furthermore, from Table~\ref{table1} we can also see that none of the extremal free operations can increase $|\alpha|$ (only make it undefined, in the case of post-processing with $\mathtt{R_y}$). Therefore, we also have $\PP_{B}^{1}\not\to \PP_{B}^{2}$ in the pre-order of this resource theory.
\end{proof}

\subsubsection{Depiction of resources and their downward closure}\label{se:geom}

The space of probability distributions over a variable that takes four possible values can be depicted as a simplex with four vertices (a tetrahedron). In particular, each one of our probability distributions over bit-to-bit functions, $\PP_{B}$, corresponds to a point inside the tetrahedron whose vertices are $[\Id],[\F],[\Rzero]$, and $[\Rone]$. Recall the parametrization of $\PP_{B}$ in terms of the parameters $(\alpha,\beta,\gamma)$ given in Eq.~\eqref{generalparametrization}. 

\begin{figure}[htbp]
\centering
\begin{subfigure}[b]{0.3\textwidth}
  \centering
  \includegraphics[width=\textwidth]{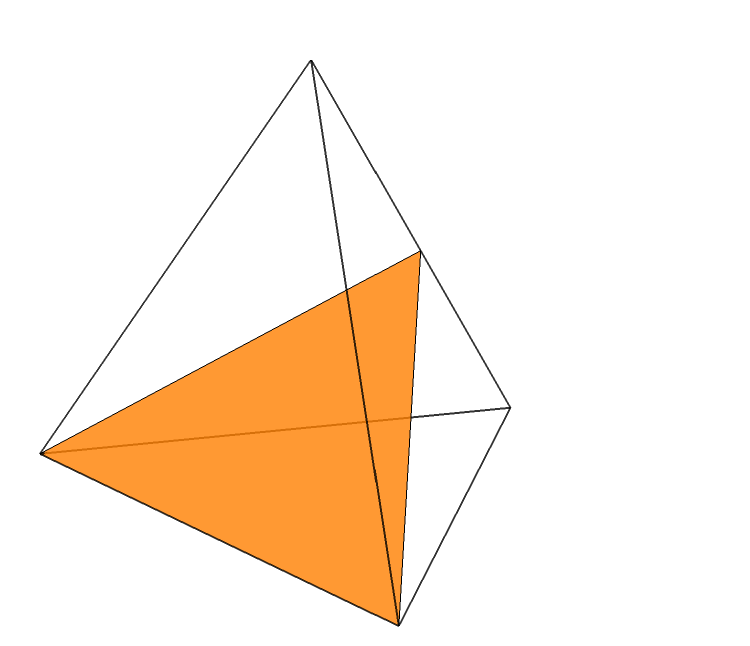}
       \put(-145,30){$[\Rzero]$}
        \put(-65,-7){$[\Rone]$}
        \put(-85,115){$[\Id]$}
        \put(-40,40){$[\F]$} 
  \caption{}
  \label{fig:abca}
\end{subfigure}
\hfill
\begin{subfigure}[b]{0.3\textwidth}
  \centering
  \includegraphics[width=\textwidth]{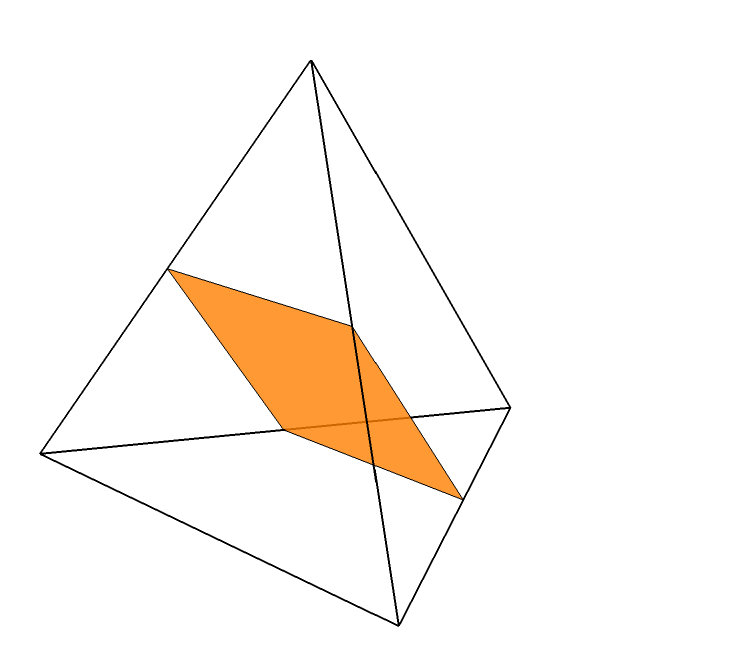}
       \put(-145,30){$[\Rzero]$}
        \put(-67,-7){$[\Rone]$}
        \put(-85,115){$[\Id]$}
        \put(-40,40){$[\F]$}
  \caption{}
  \label{fig:abcb}
\end{subfigure}
\hfill
\begin{subfigure}[b]{0.3\textwidth}
  \centering
  \includegraphics[width=\textwidth]{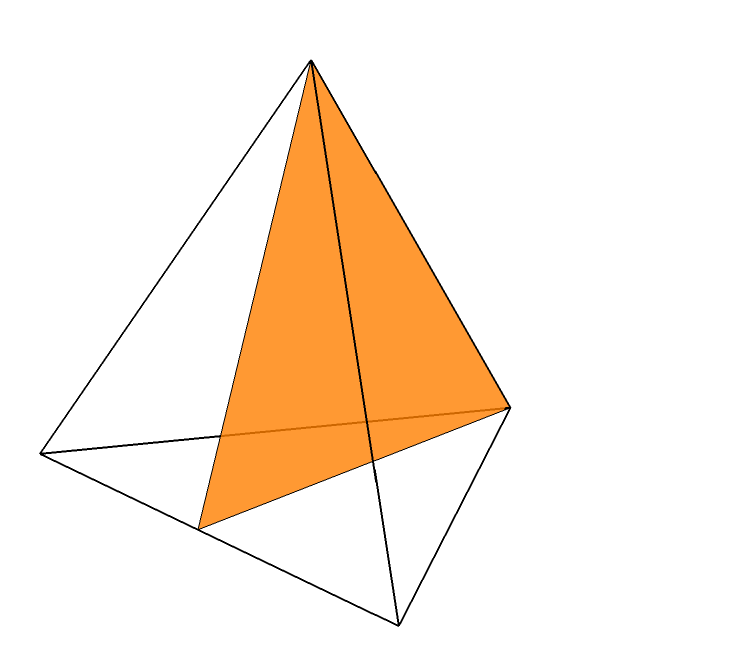}
       \put(-145,30){$[\Rzero]$}
        \put(-68,-7){$[\Rone]$}
        \put(-85,115){$[\Id]$}
        \put(-40,40){$[\F]$}
  \caption{}
  \label{fig:abcc}
\end{subfigure}
     \caption{The set of all possible resources can be illustrated as a tetrahedron formed by the convex hull of the functions ${\Id,\F,\Rzero,\Rone}$. (a) The plane representing all resources with $\alpha=0$. (b) The plane that corresponds to all resources with $\beta=\frac{1}{2}$. (c) The plane representing all resources with $\gamma=0$.}
\end{figure}
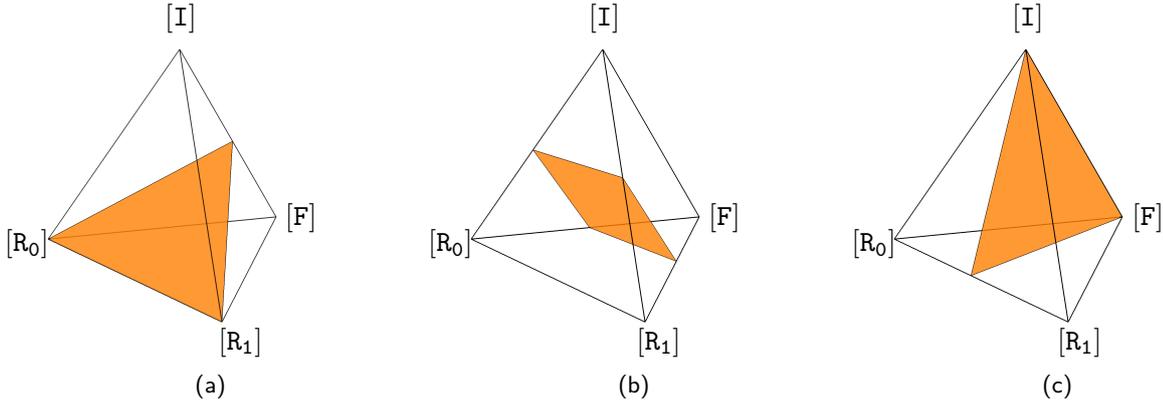

The parameter $\alpha$ corresponds to the bias towards the function $[\Id]$ or $[\F]$, relative to the plane defined by the points  $\frac{1}{2}[\Id] + \frac{1}{2}[\F]$, $[\Rzero]$, and $[\Rone]$. Here, $\alpha=0$ corresponds to this plane, $\alpha = -1$ corresponds to the $[\Id]$ vertex and $\alpha = 1$ corresponds to the $[\F]$ vertex. The intersection of the tetrahedron with the $\alpha=0$ plane is depicted in Fig.~\ref{fig:abca}.

The parameter $\beta$ represents the weight of the functions that carry causal influence. Resources with $\beta = 1$ lie on the edge connecting $[\Id]$ and $[\F]$, while those with $\beta = 0$ lie on the edge connecting $[\Rzero]$ and $[\Rone]$. Increasing $\beta$ from 0 to 1 can be visualized as moving a plane from the edge $[\Rzero]$--$[\Rone]$ to the edge $[\Id]$--$[\F$]. The intersection between the tetrahedron and the $\beta = \frac{1}{2}$ plane is depicted in Fig.~\ref{fig:abcb}.

Finally, the parameter $\gamma$ corresponds to the bias towards the function $[\Rzero]$ and $[\Rone]$, relative to the plane defined by the points $\frac{1}{2}[\Rzero] + \frac{1}{2}[\Rone]$, $[\Id]$, and $[\F]$. Here, $\gamma=0$ corresponds to this plane, $\gamma = -1$ corresponds to the $[\Rzero]$ vertex and $\gamma = 1$ corresponds to the $[\Rone]$ vertex. The intersection of the tetrahedron with the $\gamma=0$ plane is depicted in Fig.~\ref{fig:abcc}.

Let us now revisit the extremal free operations described in Sec.~\ref{se:b2brc} and discuss how they can be visualized geometrically. 

Transformations (i) and (ii) are easy to visualize. Transformation (i), given by ($\Gamma_{\textrm{pre}}=\Id$, $\Gamma_{\textrm{post}}=\Id$), maps the resource to itself, and transformation (ii), given by ($\Gamma_{\textrm{pre}}=\Id$, $\Gamma_{\textrm{post}}=\mathtt{R_y}$), maps any resource to $[\mathtt{R_y}]$, given by one of the vertices $[\Rzero]$ or $[\Rone]$.

\begin{figure}[htbp]
    \centering
    \begin{subfigure}[b]{0.22\textwidth}
        \centering
        \includegraphics[width=\textwidth]{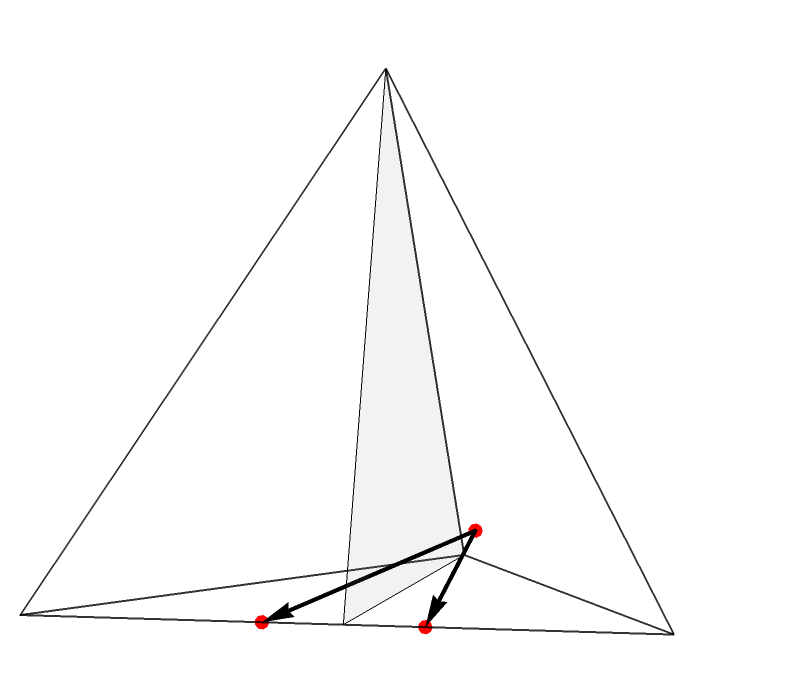}
        \put(-102,0){$[\Rzero]$}
        \put(-20,-3){$[\Rone]$}
        \put(-57,85){$[\Id]$}
        \put(-40,15){$_{(iii)}$}
        \put(-70,20){$_{(iv)}$}
        \caption{}
        \label{fig:freetransa}
    \end{subfigure}
    \hfill
    \begin{subfigure}[b]{0.23\textwidth}
        \centering
        \includegraphics[width=\textwidth, trim=0 55 0 0, clip]{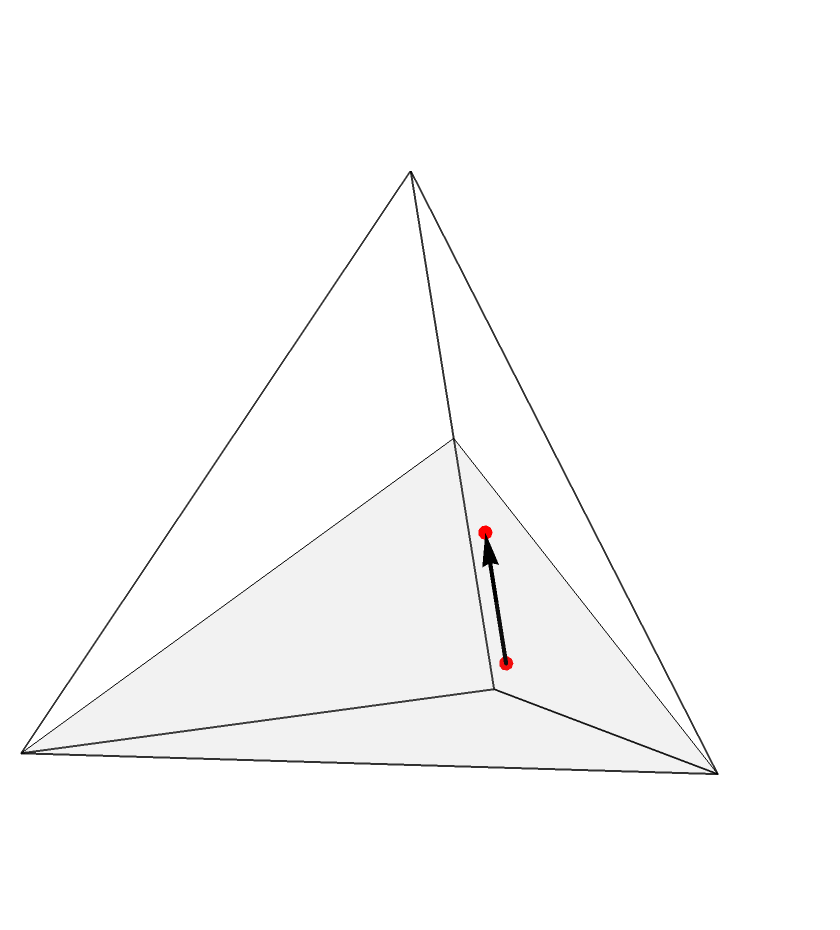}
         \put(-108,-2){$[\Rzero]$}
        \put(-18,-4){$[\Rone]$}
        \put(-59,85){$[\Id]$}
        \put(-60,25){$_{(v)}$}
        \caption{}
        \label{fig:freetransb}
    \end{subfigure}
    \hfill
    \begin{subfigure}[b]{0.24\textwidth}
        \centering
        \includegraphics[width=\textwidth]{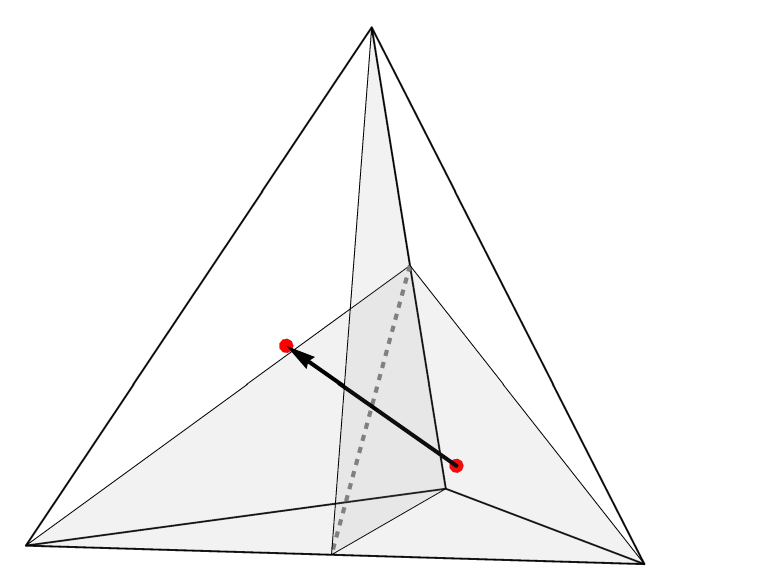}
        \put(-112,-5){$[\Rzero]$}
        \put(-18,-8){$[\Rone]$}
        \put(-62,85){$[\Id]$}
        \put(-73,20){$_{(vi)}$}
        \caption{}
        \label{fig:freetransc}
    \end{subfigure}
    \hfill
    \begin{subfigure}[b]{0.24\textwidth}
        \centering
        \includegraphics[width=\textwidth]{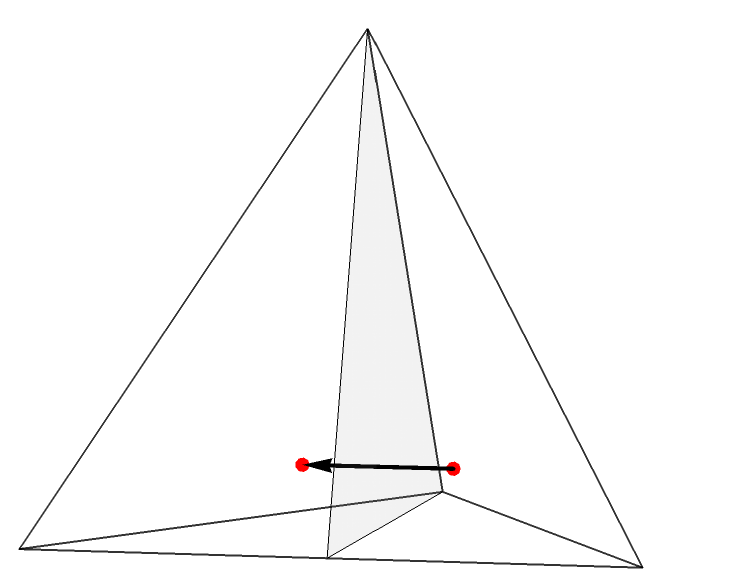}
        \put(-112,-4){$[\Rzero]$}
        \put(-18,-8){$[\Rone]$}
        \put(-61,85){$[\Id]$}
        \put(-65,25){$_{(vii)}$}
        \caption{}
        \label{fig:freetransd}
    \end{subfigure}
    
    \caption{Geometrical visualization of the extremal free operations applied to a resource $\PP_{B} = \frac{1}{8} [\Id] + \frac{3}{8}[\F] + \frac{1}{8}[\Rzero] + \frac{3}{8}[\Rone]$, which corresponds to $(\alpha,\beta,\gamma)=(\frac{1}{2},\frac{1}{2},\frac{1}{2})$. (a) Transformation (iii), given by ($\Gamma_{\textrm{pre}}=\mathtt{R_y}$, $\Gamma_{\textrm{post}}=\Id$) and transformation (iv), given by $(\Gamma_{\textrm{pre}}=\mathtt{R_y}$, $\Gamma_{\textrm{post}}=\F)$. The plane $\gamma=0$ is depicted. (b) Transformation (v), given by ($\Gamma_{\textrm{pre}}=\F$, $\Gamma_{\textrm{post}}=\Id$). The plane $\alpha=0$ is depicted. (c) Transformation (vi), given by ($\Gamma_{\textrm{pre}}=\Id$, $\Gamma_{\textrm{post}}=\F$). The plane $\gamma=0$ and the plane $\alpha=0$, as well as the line $\alpha=\gamma=0$, are depicted. (d) Transformation (vii), given by ($\Gamma_{\textrm{pre}}=\F$, $\Gamma_{\textrm{post}}=\F$). The plane $\gamma=0$ is depicted.}
\end{figure}

Both transformations (iii) and (iv) map any resource to a free resource. These transformations are represented in Fig.~\ref{fig:freetransa}, where we start from the initial resource $\PP_{B} = \frac{1}{8} [\Id] + \frac{3}{8}[\F] + \frac{1}{8}[\Rzero] + \frac{3}{8}[\Rone]$, which corresponds to $(\alpha,\beta,\gamma)=\left(\frac{1}{2},\frac{1}{2},\frac{1}{2}\right)$. Transformation (iii), given by ($\Gamma_{\textrm{pre}}=\mathtt{R_y}$, $\Gamma_{\textrm{post}}=\Id$), maps $\PP_{B}$ to $\frac{1}{4}[\Rzero] + \frac{3}{4}[\Rone]$, which corresponds to $(\alpha',\beta',\gamma')=\left(\textrm{unspecified},0,\frac{1}{2}\right)$. Transformation (iv), given by $(\Gamma_{\textrm{pre}}=\mathtt{R_y}$, $\Gamma_{\textrm{post}}=\F)$,  gives the  output resource  $(\alpha'',\beta'',\gamma'')=\left(\textrm{unspecified},0,- \frac{1}{2}\right)$, where the sign of $\gamma$ is flipped. In Fig.~\ref{fig:freetransa} the $\gamma=0$ plane is shown for reference. 

Transformations (v), (vi) and (vii) can be understood as reflections with respect to an appropriate plane or line in the polytope. Transformation (v), given by ($\Gamma_{\textrm{pre}}=\F$, $\Gamma_{\textrm{post}}=\Id$), is shown in Fig.~\ref{fig:freetransb}: the resource is simply flipped around the $\alpha=0$ plane. In the figure, we depict the $\alpha=0$ plane and illustrate how a resource defined by $(\alpha,\beta,\gamma)=\left(\frac{1}{2},\frac{1}{2},\frac{1}{2}\right)$ is transformed to $(\alpha',\beta',\gamma')=\left(-\frac{1}{2},\frac{1}{2},\frac{1}{2}\right)$. Transformation (vi), given by ($\Gamma_{\textrm{pre}}=\Id$, $\Gamma_{\textrm{post}}=\F$), is illustrated in Fig.~\ref{fig:freetransc}: the resource is flipped around the $\alpha=\gamma=0$ line (the $\alpha=0$ plane and the $\gamma=0$ plane are also illustrated) and transformed into $(\alpha',\beta',\gamma')=\left(-\frac{1}{2},\frac{1}{2},-\frac{1}{2}\right)$. Finally, consider transformation (vii), ($\Gamma_{\textrm{pre}}=\F$, $\Gamma_{\textrm{post}}=\F$), showed in Fig.~\ref{fig:freetransd}. Here,  the resource is first flipped around the $\alpha=0$ plane, and then flipped around the $\alpha=\gamma=0$ line. This transformation leaves $\alpha$ and $\beta$ invariant, and changes the sign of $\gamma$. This effectively corresponds to flipping the resource around the $\gamma=0$ plane.

\begin{figure}[htbp]
    \centering
    \begin{subfigure}[b]{0.4\textwidth}
        \centering
        \includegraphics[width=\textwidth]{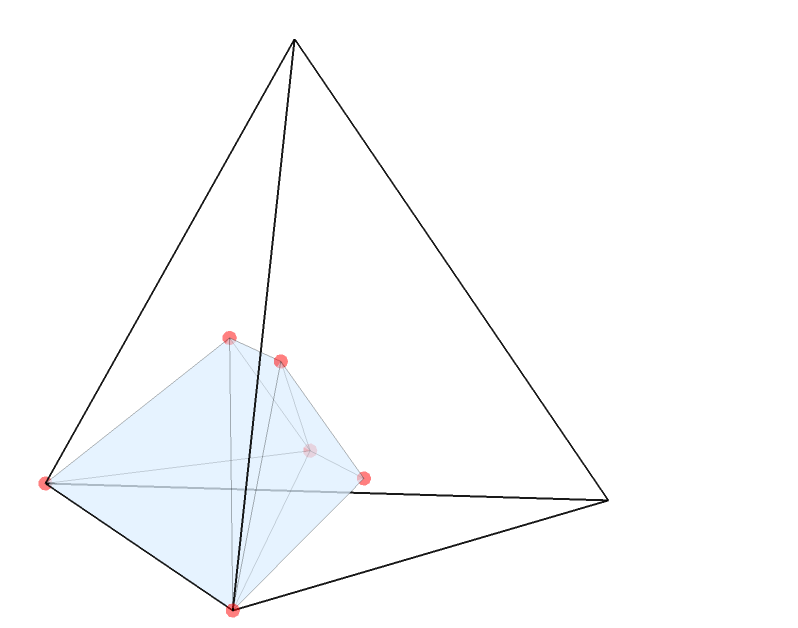}
        \put(-185,20){$[\Rzero]$}
        \put(-135,-10){$[\Rone]$}
        \put(-120,140){$[\Id]$}
        \put(-40,25){$[\F]$}
        \put(-95,35){$\PP_B$}
        \caption{}
        \label{fig_downward_closure_bits}
    \end{subfigure}
    \hfill
    \begin{subfigure}[b]{0.4\textwidth}
        \centering
        \includegraphics[width=\textwidth]{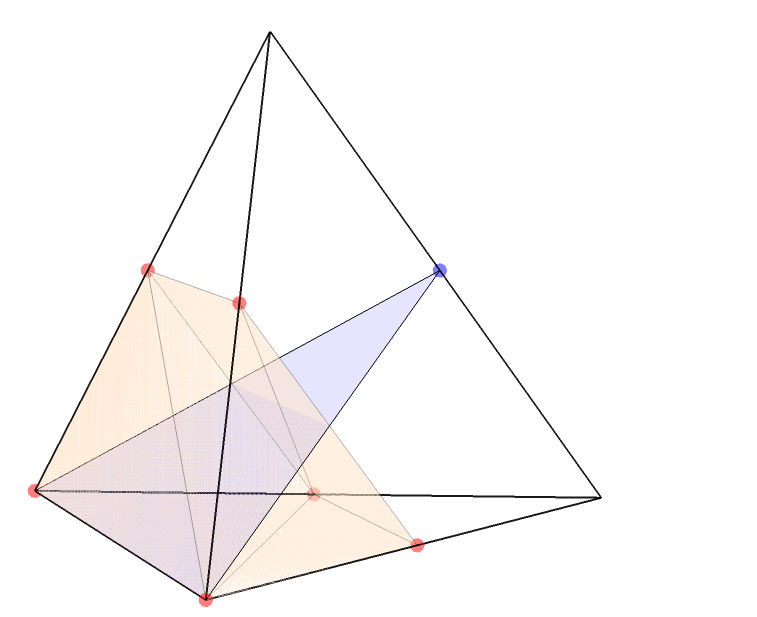}
        \put(-190,20){$[\Rzero]$}
        \put(-137,-12){$[\Rone]$}
        \put(-123,145){$[\Id]$}
        \put(-35,25){$[\F]$}
        \put(-70,85){$\PP_{B}^{1}$}
        \put(-165,85){$\PP_{B}^{2}$}
        \caption{}
        \label{fig_downward_closure_prop3}
    \end{subfigure}

    \caption{(a) The downward closure of the resource $\PP_{B} = \frac{1}{8} [\Id] + \frac{3}{8}[\F] + \frac{1}{8}[\Rzero] + \frac{3}{8}[\Rone]$ is the convex hull of the 6 red points of this figure. (b) The downward closure of $\PP_{B}^{1}=\frac{1}{2}[\Id]+\frac{1}{2}[\F]$ is indicated by the blue triangle, and the downward closure of $\PP_{B}^{2}=\frac{1}{2}[\Id]+\frac{1}{2}[\Rzero]$ is the convex hull of the six red points, indicated in orange.  Note that for the pair $\PP_{B}^{1}$ and $\PP_{B}^{2}$, neither resource is included in the downward closure of the other.}
\end{figure}

By analyzing the extremal free operations, we can see that the downward closure $\mathcal{P}_{\downarrow}$ is given by the convex hull of (at most) the 6 resources  given by transformations (i), (ii), (v), (vi), and (vii). The six extremal points that define the polytope $\mathcal{P}_{\downarrow}$ for the resource $(\alpha,\beta,\gamma)=\left(\frac{1}{2},\frac{1}{2},\frac{1}{2}\right)$ are represented by the six red points in Fig.~\ref{fig_downward_closure_bits}. 

This geometrical visualization provides an alternative proof of Prop.~\ref{prop:nottot}: it suffices to note that $\PP_{B}^{1}$ and $\PP_{B}^{2}$ are not contained in the downward closure of the other, as made explicit in Fig.~\ref{fig_downward_closure_prop3}.

\subsubsection{Complete set of monotones}\label{se:b2bcsm}

In this section, we prove that three resource monotones characterize the pre-order of resources in  $\RTKnow$ when $X$ and $Y$ are binary.

\begin{restatable}{prop}{Prop}\label{prop:complete}
    Consider resources in a bit-to-bit scenario, parametrised as in Eq.~\eqref{generalparametrization}. The following three functions of the parameters  $(|\alpha|,\beta,|\gamma|)$ form a complete set of monotones for this scenario:
\begin{align}
&M_\beta(\PP_B) :=\beta\,, \label{eq_beta_monotone}\\
&M_{|\alpha|}(\PP_{B}) := |\alpha|\,,\label{eq_alpha_monotone}
\\
&M_{|\gamma|,\beta}(\PP_{B}) := \begin{cases} \frac{\beta}{1-|\gamma|(1-\beta)} \quad &\text{if} \quad \beta \neq 1\,,\\
\qquad 1 \quad &\text{if} \quad \beta = 1\,.
\end{cases}\label{eq_third_monotone}
\end{align} 
\end{restatable}

The proof is given in Appendix~\ref{app:proof-monotones}.

Next, we turn to the question: are all three monotones {\em necessary} to specify the pre-order of this resource theory? As we show now, the answer is yes;  it is not possible to specify the pre-order with less information than that provided by all  three monotones, or equivalently, by all three parameters $(|\alpha|,\beta,|\gamma|)$.

\begin{restatable}{lem}{Lemma}\label{lem:canform}
The three monotones $M_\beta$, $M_{|\alpha|}$ and $M_{|\gamma|,\beta}$, or equivalently the three parameters $(|\alpha|,\beta,|\gamma|)$, are necessary and sufficient for specifying the equivalence class of a resource, and hence may be taken as a canonical form of a resource. 
\end{restatable}
The proof is given in Appendix~\ref{app:proof-lemma}. We first note that one can compute the values of $(|\alpha|,\beta,|\gamma|)$ from the values of $(M_\beta,M_{|\alpha|},M_{|\gamma|,\beta})$, so that these two choices are really interconvertible.  We prove that all three parameters  are necessary by giving explicit examples of pairs of resources for which each of the three is necessary to determine the relative ordering of the pair. These monotones in turn completely specify the values of $(|\alpha|,\beta,|\gamma|)$.

Note that, unlike the other two monotones, $M_{|\alpha|}$ is a  ``partial monotone'': it is not defined on all the resources from the enveloping theory; it is only defined for nonfree resources.

\begin{figure}[htbp]
\centering
\begin{subfigure}[b]{0.3\textwidth}
  \centering
  \includegraphics[width=\textwidth]{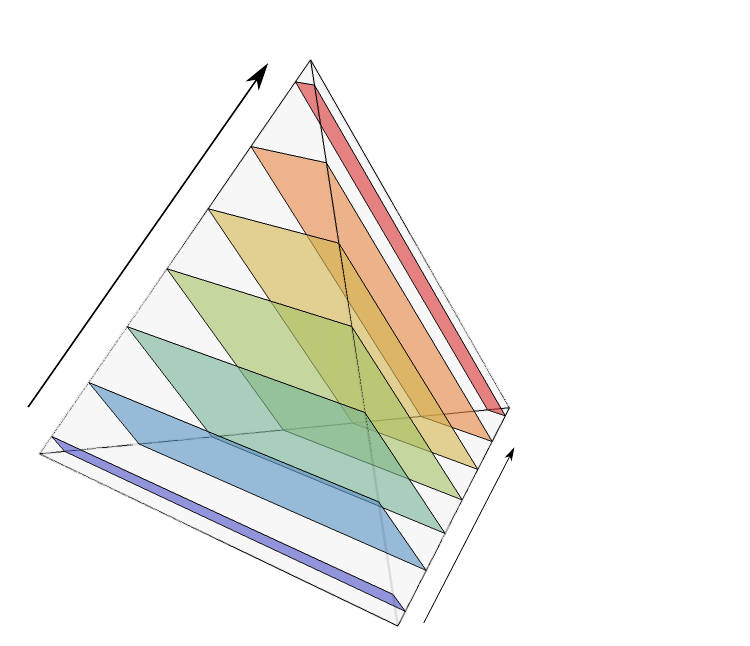}
       \put(-145,30){$[\Rzero]$}
        \put(-65,-7){$[\Rone]$}
        \put(-85,115){$[\Id]$}
        \put(-40,40){$[\F]$}
  \caption{}
  \label{fig:betalevels}
\end{subfigure}
\hfill
\begin{subfigure}[b]{0.3\textwidth}
  \centering
  \includegraphics[width=\textwidth]{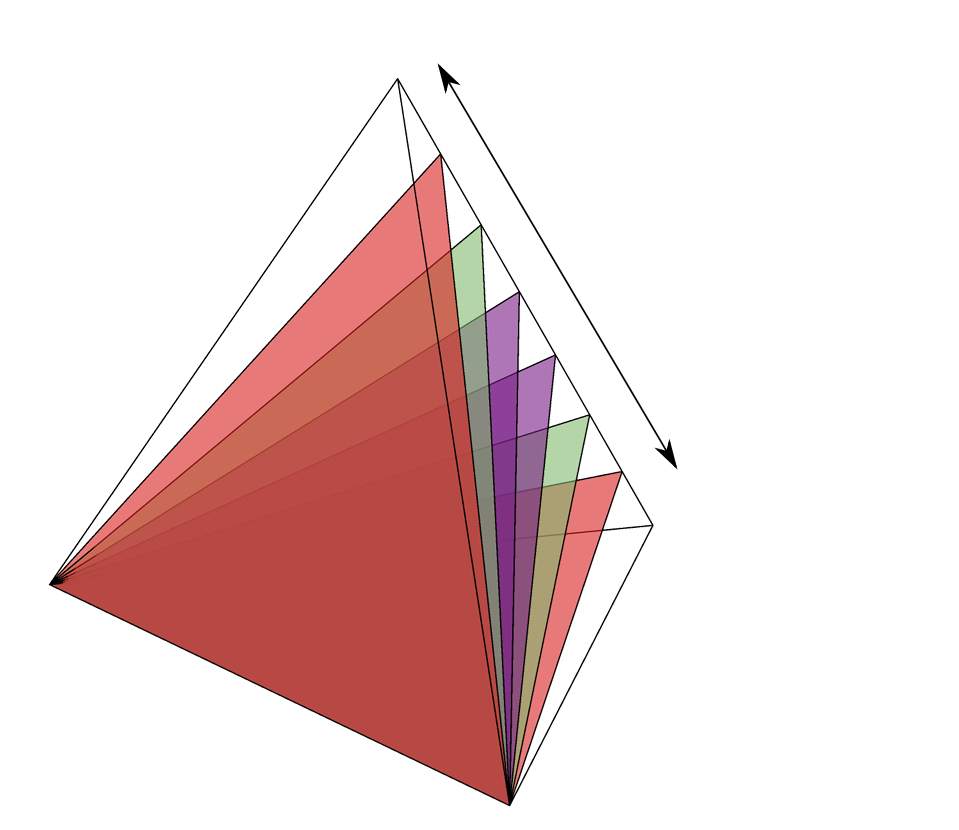}
       \put(-145,30){$[\Rzero]$}
        \put(-67,-7){$[\Rone]$}
        \put(-85,115){$[\Id]$}
        \put(-40,40){$[\F]$}
  \caption{}
  \label{fig:alphalevels}
\end{subfigure}
\hfill
\begin{subfigure}[b]{0.3\textwidth}
  \centering
  \includegraphics[width=\textwidth, trim=0 25 0 0, clip]{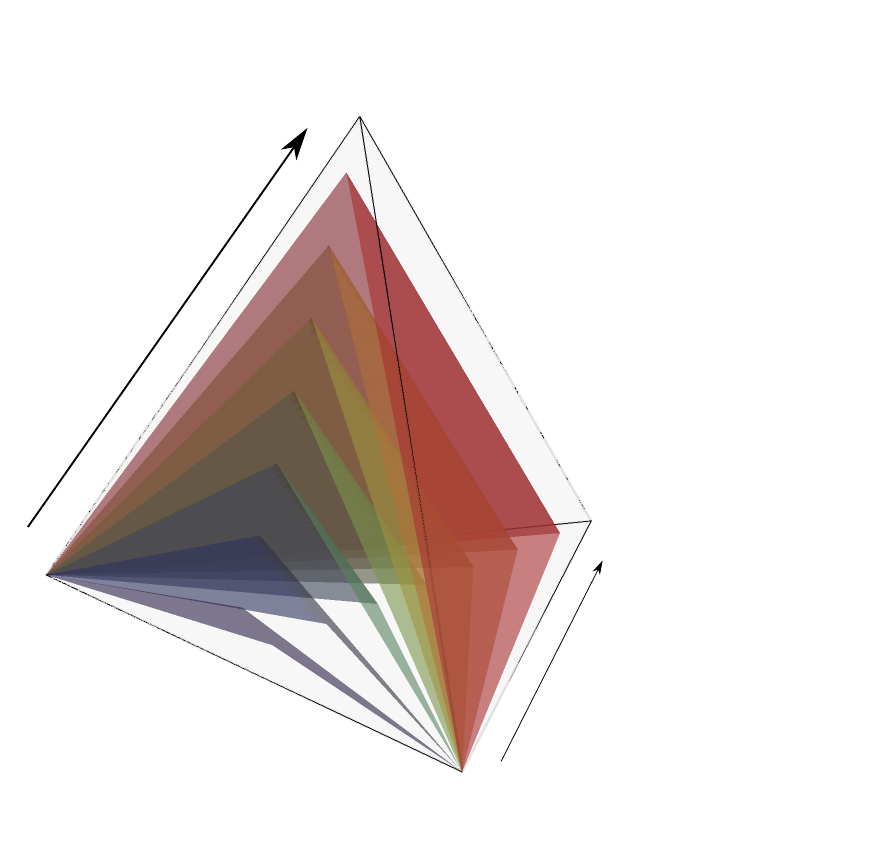}
       \put(-145,30){$[\Rzero]$}
        \put(-68,-7){$[\Rone]$}
        \put(-85,115){$[\Id]$}
        \put(-40,40){$[\F]$}
  \caption{}
  \label{fig:gammalevels}
\end{subfigure}
     \caption{Level curves of the monotones (a) $M_\beta$, (b) $M_{|\alpha|}$ and (c) $M_{|\gamma|,\beta}$.  The arrows indicate the direction towards which the level curves have higher values of the monotone. }
\end{figure}
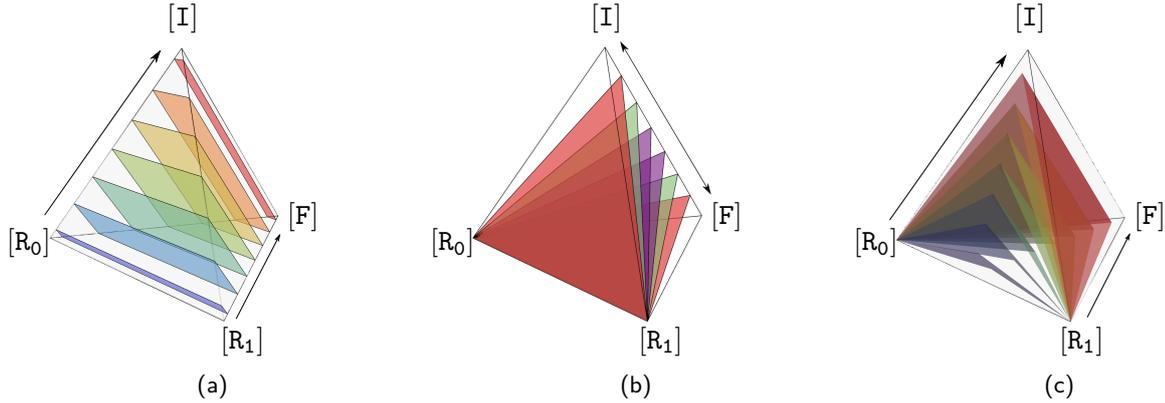

We finish this subsection by discussing the level curves of each of the monotones. For the case of $M_\beta(\PP_B)$, the level curves are depicted in Fig.~\ref{fig:betalevels}. They are simply given by the horizontal planes defined by different values of $\beta$. 
The free operations cannot take a resource that lies on a lower plane (closer to the $[\Rzero]$-$[\Rone]$ line) to a higher plane (closer to the  $[\Id]$-$[\F]$ line).

The level curves of the monotone $M_{|\alpha|}(\PP_{B})$ are shown in Fig.~\ref{fig:alphalevels}. They correspond to pairs of planes defined by different values of $|\alpha|$. The free operations cannot take a resource that lies on an inner pair of planes (closer to $|\alpha|=0$) to one on an outer pair of planes (further from $|\alpha|=0$).  

Finally, for the case of $M_{|\gamma|,\beta}(\PP_{B})$, the level curves are the nested polytopes of Fig.~\ref{fig:gammalevels}. The green lines of Fig.~\ref{fig:third_theta} are a cross section of the level curves of Fig.~\ref{fig:gammalevels}. The free operations cannot take a resource that lies on an inner polytope to one on an outer polytope.

Each face of the downward closure of a resource lies on a level curve of one of the monotones. For example, the downward closure of Fig.~\ref{fig_downward_closure_bits} is such that the face that does not include  $[\Rzero]$ nor $[\Rone]$ lies on a level curve of $M_\beta$, the two faces that include both $[\Rzero]$ and $[\Rone]$ lie on a level curve of $M_{|\alpha|}$ and the two faces that respectively include only $[\Rzero]$ and  only $[\Rone]$ lie on a level curve of $M_{|\gamma|,\beta}$. This fact provides an alternative proof, one that is more geometrical, of the fact that $M_\beta$, $M_{|\alpha|}$ and $M_{|\gamma|,\beta}$ are monotones.

\subsubsection{Interpretation of the monotones}

 In this section, we present interpretations of each of the monotones $(M_\beta,M_{|\alpha|},M_{|\gamma|,\beta})$. 

Let us begin with $M_\beta$.   By definition (Eq.~\eqref{eq_beta_monotone}), $M_\beta(\PP_B)$ is simply the value of $\beta$ in the parametrization of $\PP_B$. Recall that $\beta$ determines the probability that the functional dependence of the output bit on the input bit is given by a causally connected function ($\Id$ or $\F$). 
 $M_{\beta}$ can therefore be described as the {\em weight of causally connected functions} in the probability distribution over functions, or simply 
{\em the probability of causal connection}.\footnote{This corresponds to what is termed {\em the probability of necessity and sufficiency } in Ref.~\cite{pearl2009causality}}  
$M_\beta$ cannot increase under the free operations because the free operations do not provide any means of increasing 
the weight of  causally connected functions, i.e., causal connectivity cannot be generated freely. 

Next, we consider $M_{|\alpha|}$.  
By definition (Eq.~\eqref{eq_alpha_monotone}), $M_{|\alpha|}(\PP_{B})$ is simply  $|\alpha|$ in the parameterization of $\PP_{B}$. 
It is straightforward to verify that $$|\alpha | = \frac{|\PP_{B}(\Id)-\PP_{B}(\F)|}{\PP_{B}(\Id)+\PP_{B}(\F)},$$
so that $|\alpha|$ is aptly described as the {\em polarization} between the causally connected functions in the probability distribution over functions, or simply the 
{\em $\Id{-}\F$ polarization}.  It is  the bias towards either $\Id$ or $\F$ and away from the equal mixture of the two. 
The fact that the monotone $M_{|\alpha|}$ cannot increase captures the fact that one cannot freely increase  this polarization. 
Thus, as we already discussed in Sec.~\ref{sec:intro}, there is a kind of resourcefulness that corresponds to knowledge about {\em the precise nature} of the causal connectivity (and it is quantified by $|\alpha|$). 
Finally, we consider
$M_{|\gamma|, \beta}$.
 $M_{|\gamma|, \beta}(\PP_B)$ is defined in Eq.~\eqref{eq_third_monotone} and is a function of $\beta$ and $|\gamma|$. 
It is straightforward to verify that 
$$|\gamma | = \frac{|\PP_{B}(\Rzero)-\PP_{B}(\Rone)|}{\PP_{B}(\Rzero)+\PP_{B}(\Rone)},$$
so that $|\gamma|$ is aptly described as the {\em polarization} between the causally {\em disconnected} functions in the probability distribution over functions, or simply the {\em $\Rzero{-}\Rone$ polarization}. It describes the bias towards either $\Rzero$ or $\Rone$ and away from the equal mixture of the two.

Unlike $\beta$ and $|\alpha|$, the parameter $|\gamma|$ {\em can} be increased under the free operations.
 To increase the bias towards $\Rone$, for instance, it suffices to implement a mixture of $\Id$ and $\Rone$ via post-processing.
However, 
 such an operation 
necessarily decreases the value of $\beta$ if the latter is nonzero.  

 For example, consider a case where the starting resource has some bias towards $[\Rone]$, such as
\begin{equation}  \PP_B^7=\frac{1}{3}\left(\frac{1}{2}[\Id]+\frac{1}{2}[\F]\right)+\frac{2}{9}[\Rzero]+\frac{4}{9}[\Rone].
\end{equation}
It is possible to sacrifice  some probability of causal connection 
to increase the bias towards $[\Rone]$ using the free operations. For example, imagine that with a probability of $3/4$ we post-process 
with $\Id$ and with a probability of $1/4$ we post-process 
with $\Rone$. The resulting probability distribution over functions is:
\begin{equation}
\PP_B^8=\frac{1}{4}\left(\frac{1}{2}[\Id]+\frac{1}{2}[\F]\right)+\frac{1}{6}[\Rzero]+\frac{7}{12}[\Rone].
\end{equation}
The bias towards $[\Rone]$ has increased, but $\beta$  has decreased from $1/3$ to $1/4$. 

 The monotone $M_{|\gamma|, \beta}$ quantifies the amount by which $\beta$ must decrease in such an operation, in a manner that we will now make precise.

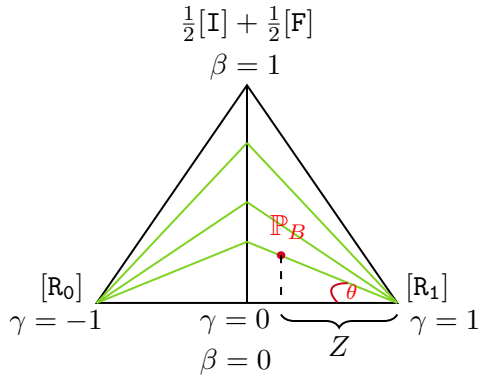
\begin{figure}[h]
    \centering
    \tikzset{every picture/.style={line width=0.75pt}}      
\usetikzlibrary{decorations.pathreplacing}

\begin{tikzpicture}[x=0.75pt,y=0.75pt,yscale=-1,xscale=1]

\draw (120,71.25) -- (195,180.75) -- (45,180.75) -- cycle ;
\draw (120,71.25) -- (120,180.75) ;

\draw [color={rgb, 255:red, 126; green, 211; blue, 33}, draw opacity=1] (120,100.25) -- (45,180.75) ;
\draw [color={rgb, 255:red, 126; green, 211; blue, 33}, draw opacity=1] (120,130) -- (45,180.75) ;
\draw [color={rgb, 255:red, 126; green, 211; blue, 33}, draw opacity=1] (120,150) -- (45,180.75) ;

\draw [color={rgb, 255:red, 126; green, 211; blue, 33}, draw opacity=1] (120,100.25) -- (195,180.75) ;
\draw [color={rgb, 255:red, 126; green, 211; blue, 33}, draw opacity=1] (120,130) -- (195,180.75) ;
\draw [color={rgb, 255:red, 126; green, 211; blue, 33}, draw opacity=1] (120,150) -- (195,180.75) ;

\draw [color={rgb, 255:red, 208; green, 2; blue, 27}, draw opacity=1]
      [fill={rgb, 255:red, 208; green, 2; blue, 27}, fill opacity=1]
      (135.38,156.75) .. controls (135.38,155.92) and (136.1,155.25) .. (137,155.25)
      .. controls (137.9,155.25) and (138.63,155.92) .. (138.63,156.75)
      .. controls (138.63,157.58) and (137.9,158.25) .. (137,158.25)
      .. controls (136.1,158.25) and (135.38,157.58) .. (135.38,156.75) -- cycle ;

\draw [color={rgb, 255:red, 208; green, 2; blue, 27}, draw opacity=1]
      (165.13,180.5) .. controls (155.81,170.17) and (170.56,170.67) .. (170.31,170.92) ;


\draw[dashed] (137,156.75) -- (137,180.75);

\draw[decorate,decoration={brace,amplitude=5pt,mirror}]
      (137,186.5) -- (195,186.5)
      node[midway, yshift=-12pt] {$Z$};


\draw (85,30) node [anchor=north west][inner sep=0.75pt] [align=left] {$\frac{1}{2}[\Id]+\frac{1}{2}[\F]$};
\draw (100,54) node [anchor=north west][inner sep=0.75pt] [align=left] {$\beta=1$};
\draw (95,182.75) node [anchor=north west][inner sep=0.75pt] [align=left] {$\gamma=0$};
\draw (95,202.25) node [anchor=north west][inner sep=0.75pt] [align=left] {$\beta=0$};

\draw (197.5,164.25) node [anchor=north west][inner sep=0.75pt] [align=left] {$[\Rone]$};
\draw (200,183.25) node [anchor=north west][inner sep=0.75pt] [align=left] {$\gamma=1$};

\draw (15,164.25) node [anchor=north west][inner sep=0.75pt] [align=left] {$[\Rzero]$};
\draw (0,183.25) node [anchor=north west][inner sep=0.75pt] [align=left] {$\gamma=-1$};

\draw (168,170) node [anchor=north west][inner sep=0.75pt] [align=left] {\red$_\theta$\blk};
\draw (130,135) node [anchor=north west][inner sep=0.75pt] [align=left] {\red$\PP_B$\blk};

\end{tikzpicture}
    \caption{For a resource $\PP_B$ living on the plane $\alpha=0$, the monotone $M_{|\gamma|,\beta}$ is the tangent of the angle $\theta$ indicated here.}
    \label{fig:third_theta}
\end{figure}

The kind of free operation that we used here, consisting of post-processing with $\Id$ with a probability $p$ and with $\Rone$ with a probability $(1-p)$, moves the resource down the line that connects it to the vertex $[\Rone]$. This line lies on a level curve of the monotone $M_{|\gamma|,\beta}$. An example is depicted in Fig.~\ref{fig:third_theta}: there, we show a resource $\PP_B$ that lives in the plane $\alpha=0$ (like $\PP_B^7$ and $\PP_B^8$), and the cross section of some of the level curves of  $M_{|\gamma|,\beta}$ on this plane, depicted as green lines. The distance between the point $[\Rone]$ and the projection of the point $\PP_B$ down to the line $[\Rzero]{-}[\Rone]$ is $Z\equiv 1-|\gamma|(1-\beta)$. The monotone  $M_{|\gamma|,\beta}$ is exactly the tangent of the angle $\theta$, that is, it measures the ratio between how much the point moves down (decrease in $\beta$) and how much it moves to the right (increase in $Z$). Intuitively, when we start with a resource that has some bias towards $[\Rzero]$ (respectively $[\Rone]$), the monotone $M_{|\gamma|,\beta}$ quantifies how much  probability of causal connection 
needs to be sacrificed to increase the bias towards $[\Rzero]$ (respectively $[\Rone]$) as measured by $Z$ by using free operations that take the form of post-processing with $\Id$ with a probability $p$ and with  $\Rzero$ (respectively $\Rone$) with a probability $(1-p)$.

Note that the monotone $M_{|\gamma|,\beta}$ does not measure the tradeoff between  the probability of causal connection 
and the $\Rzero-\Rone$ polarization
(which is measured by $|\gamma|$), but instead the tradeoff between 
 the probability of causal connection 
and  the value of $Z$.  
One can check that 
$Z=1-|\PP_B(\Rzero)-\PP_B(\Rone)|$,  so that, up to a linear transformation, it describes an {\em unnormalized} version of the $\Rzero{-}\Rone$ polarization.

Note that although it is possible, by decreasing  
the probability of causal connection, 
to obtain 
more certainty about the nature of the causal dependence (specifically, along the $[\Rzero]-[\Rone]$ direction), the reverse process is not possible, as the probability of causal connection 
cannot be freely increased by any means.

An alternative interpretation of  $M_{|\gamma|,\beta}$  is provided in Sec.~\ref{sec:examples}, which is related to the motivating example presented in Sec.~\ref{sec_flag}. 

To summarise: 
\begin{compactitem}
\item The $M_\beta$ monotone quantifies  the probability of causal connection in $\PP_B$, i.e., the weight assigned to functions that describe causal connection 
\item The $M_{|\alpha|}$ monotone quantifies the  $\Id{-}\F$ polarization in $\PP_B$. 
\item  The $M_{|\gamma|,\beta}$ monotone captures the ratio of the decrease of the probability of causal connection to the increase in an unnormalized version of the $\Rzero{-}\Rone$ polarization (quantified by  $Z= 1-|\gamma|(1-\beta)=1-|\PP_B(\Rzero)-\PP_B(\Rone)|$) in a free operation that takes $\PP_B$ to the closest resource with maximal $\Rzero{-}\Rone$ polarization (i.e., to either $[\Rzero]$ or $\Rone$).
\end{compactitem}

Note that having knowledge about whether the true causal influence is $\Rzero$ or $\Rone$, conditioned on knowing that it is one of the two, is {\em not} valuable in this resource theory.  That is, $|\gamma|$ is not a monotone.

\subsection{Back to the communication channel that leaks a flag variable to the environment}\label{sec:examples}

Let us now revisit the examples that appear in Sec.~\ref{sec_flag}, where we discussed the communication channel from a binary $X$ to a binary $Y$ that leaks a flag variable to the environment. Specifically, let us look at the case where Alice and Bob learn partial information about the flag, and thus their description of the channel is given by a probability distribution over functions. 

The probability distributions over functions that appeared in Sec.~\ref{sec_flag} were: 
\begin{alignat*}{2}
&\PP_{B}^{1} = \frac{1}{2} \, [\Id] + \frac{1}{2} \, [\F]\,, \quad\quad\; &&
\PP_{B}^{2} = \frac{1}{2} \, [\Rzero] + \frac{1}{2} \, [\Rone]\,, \\
& \PP_{B}^{3} = \frac{2}{3} \, [\Id] + \frac{1}{3} \, [\F]\, \quad\quad\; &&\PP_B^4 = \frac{1}{3} \, [\F] + \frac{2}{3} \, [\Rzero],     \\
     &\PP_B^5 = \frac{1}{3} \, [\F] + \frac{1}{3} \, [\Rzero] + \frac{1}{3} \, [\Rone],\quad\;\quad  &&\PP_B^6 = \frac{1}{3} \, \left(\frac{1}{2} [\Id]+\frac{1}{2}[\F]\right) + \frac{2}{3} \, [\Rzero].
\end{alignat*}

The values of our three monotones for each of these resources are
\begin{align}
    & \PP_{B}^{1}\leftrightarrow\left(M_\beta(\PP_{B}^{1} ),M_{|\alpha|}(\PP_{B}^{1}),M_{|\gamma|,\beta}(\PP_{B}^{1})\right)=(1,0,1), \label{eq_mon_1} \\
    &\PP_{B}^{2}\leftrightarrow\left(M_\beta(\PP_{B}^{2} ),M_{|\alpha|}(\PP_{B}^{2}),M_{|\gamma|,\beta}(\PP_{B}^{2})\right)=(0,\text{undefined},0), \label{eq_mon_2} \\
    &\PP_{B}^{3}\leftrightarrow\left(M_\beta(\PP_{B}^{3} ),M_{|\alpha|}(\PP_{B}^{3}),M_{|\gamma|,\beta}(\PP_{B}^{3})\right)=\left(1,\frac{1}{3},1\right), \label{eq_mon_3}  \\
    &\PP_{B}^{4}\leftrightarrow\left(M_\beta(\PP_{B}^{4} ),M_{|\alpha|}(\PP_{B}^{4}),M_{|\gamma|,\beta}(\PP_{B}^{4})\right)=\left(\frac{1}{3},1,1\right), \label{eq_mon_4} \\
    &\PP_{B}^{5}\leftrightarrow\left(M_\beta(\PP_{B}^{5} ),M_{|\alpha|}(\PP_{B}^{5}),M_{|\gamma|,\beta}(\PP_{B}^{5})\right)=\left(\frac{1}{3},1,\frac{1}{3}\right), \label{eq_mon_5}  \\
    &\PP_{B}^{6}\leftrightarrow\left(M_\beta(\PP_{B}^{6} ),M_{|\alpha|}(\PP_{B}^{6}),M_{|\gamma|,\beta}(\PP_{B}^{6})\right)=\left(\frac{1}{3},0,1\right). \label{eq_mon_6} 
\end{align}

\begin{figure}[h!]
    \centering
    \begin{tikzpicture}[
  node distance=1.5cm and 2cm,
  every node/.style={inner sep=1pt},
  ->, >=stealth, thick
]

\node (P3) at (0,3) {$\PP_B^3$};
\node (P1) at (0,1.5) {$\PP_B^1$};
\node (P4) at (2,1.5) {$\PP_B^4$};
\node (P5) at (2,0) {$\PP_B^5$};
\node (P6) at (0,0) {$\PP_B^6$};
\node (P2) at (1,-1) {$\PP_B^2$};

\draw (P3) -> (P1);
\draw (P4) -> (P5);
\draw (P4) -> (P6);
\draw (P1) -> (P6);
\draw (P5) -> (P2);
\draw (P6) -> (P2);

\end{tikzpicture}
    \caption{Partial order of the six examples of probability distributions over functions that appear in Sec.~\ref{sec_flag} in $\RTKnow$.}
    \label{fig_partial_order}
\end{figure}
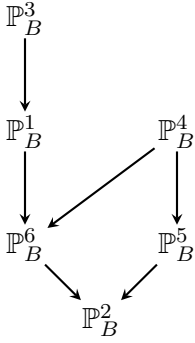

From these equations, we can infer the ordering of these six resources within the partial order of resources in $\RTKnow$, as depicted in Fig.~\ref{fig_partial_order}.

As we have seen, $\PP_B^1$ and $\PP_B^2$ give rise to the same conditional probability distribution $\PP_{Y|X}$, the completely randomizing channel. However, they encode different  probabilities of causal connection. 
As anticipated in Sec.~\ref{sec_flag}, $\PP_B^1$ has a higher value of the monotone $M_\beta$ than $\PP_B^2$. 

The resource $\PP_B^3$, on the other hand, has the same probability of causal connection 
as $\PP_B^1$. As anticipated in Sec.~\ref{sec_flag}, the difference between $\PP_B^1$ and $\PP_B^3$ is given by the monotone $M_{|\alpha|}$: 
$\PP_B^3$ has a higher $\Id{-}\F$ polarization, 
hence more resourcefulness. The same happens when comparing resources $\PP_B^4$ and $\PP_B^6$: their only difference is in the monotone $M_{|\alpha|}$, and  $\PP_B^4$ is more resourceful than $\PP_B^6$ because it has  a larger $\Id{-}\F$ polarization. 

The resources $\PP_B^4$ and $\PP_B^5$ have the same value for the monotones $M_\beta$ and $M_{|\alpha|}$. As anticipated in Sec.~\ref{sec_flag}, their difference is in the monotone $M_{\beta,|\gamma|}$. In Sec.~\ref{sec_flag} it was mentioned that, in the scenario of the communication channel that leaks a flag variable to the environment, $\PP_B^4$ is a more valuable resource than $\PP_B^5$ because it allows for Bob to know that the channel implemented a function that has causal connectivity in the rounds where $Y=1$. In the proposition below we show that, in general, the success in the task of inferring causal connectivity under postselection on a value of $Y$ is quantified by the monotone $M_{|\gamma|,\beta}$. Therefore, Prop.~\ref{prop_thirdmon_interpretation} provides an alternative interpretation of this monotone.

\begin{prop}
\label{prop_thirdmon_interpretation}
     Let Bob describe the 
    communication channel by 
    a probability distribution over functions $\PP_B$, as 
    in Sec.~\ref{sec_flag},  and assume a uniform prior over values of the channel input $X$.  Let $\PP_{B|Y}(\Id \cup \F|0)$ be the likelihood that Bob assigns to the channel having implemented a causally connected function (i.e., $\Id$ or $\F$) when he postselects on the rounds with $Y=0$.
    Analogously, let $\PP_{B|Y}(\Id \cup \F|1)$ be the likelihood that Bob assigns to the channel having implemented a causally connected function (i.e., $\Id$ or $\F$) when he postselects on the rounds with $Y=1$.
    The larger one of these two quantities is equal to the monotone $M_{|\gamma|, \beta}$.

    In words, the monotone $M_{|\gamma|, \beta}$ is the maximum probability (when varying the value of $Y$ on which one post-selects) of inferring causal connectivity.
\end{prop}
\begin{proof}
    Consider the parametrization of Eq.~\eqref{generalparametrization} for a resource:
    \begin{equation*} 
\PP_{B} = \beta\left(\frac{1-\alpha}{2}[\Id]+\frac{1+\alpha}{2}[\F]\right) +(1-\beta)\left(\frac{1-\gamma}{2}[\Rzero]+\frac{1+\gamma}{2}[\Rone]\right).
\end{equation*}

First, suppose that $\beta\neq 1$ and $\gamma>0$. In this case, the weight on $[\Rzero]$ is smaller than the weight on $[\Rone]$. Therefore, Bob will assign a higher likelihood for the function to be in the sector that carries causal influence (i.e., to $\Id$ or $\F$) if he obtains $Y=0$ than if he obtains $Y=1$. Thus, we want to calculate the likelihood that the function carries causal influence given that Bob obtained $Y=0$. That is, we want to know $\PP_{B|Y}(\Id \cup \F|0)=1-\PP_{B|Y}(\Rzero \cup \Rone|0)$. To do so, we use the following Bayesian inversion:
\begin{equation}
    \PP_{B|Y}(\Rzero \cup \Rone|0)=\frac{\PP_{Y|B}(0|\Rzero \cup \Rone)\PP_B(\Rzero \cup \Rone)}{\PP_Y(0)}. \label{eq_bayespropproof}
\end{equation}

From the parametrization, we know that $\PP_{Y|B}(0|\Rzero \cup \Rone)=\frac{1-\gamma}{2}$, and $\PP_B(\Rzero \cup \Rone)=1-\beta$. Now, we proceed to calculate $\PP_Y(0)$. To do so, we first obtain the conditional probability distribution $\PP_{Y|X}$ that is generated from our probability distribution over functions:
\begin{equation}
    \begin{cases}
        \PP_{Y|X}(0|0)=\beta\frac{(1-\alpha)}{2}+(1-\beta)\frac{(1-\gamma)}{2}\, ,\\
        \PP_{Y|X}(0|1)=\beta\frac{(1+\alpha)}{2}+(1-\beta)\frac{(1-\gamma)}{2} \, .
    \end{cases}
\end{equation}
We assume a uniform prior over $X$,  
that is, we assume that $\PP_X(0)=\PP_X(1)=1/2$. With this, we get:
\begin{equation}
    \PP_Y(0)=\PP_{Y|X}(0|0)\PP_X(0)+\PP_{Y|X}(0|1)\PP_X(1)=\frac{1-\gamma+\beta\gamma}{2}.
\end{equation}

Using this in Eq.~\eqref{eq_bayespropproof}, we get:
\begin{equation}
    \PP_{B|Y}(\Id \cup \F|0)=1-\PP_{B|Y}(\Rzero \cup \Rone|0)=1-\frac{(1-\gamma)(1-\beta)}{1-\gamma+\beta\gamma}=\frac{\beta}{1-\gamma(1-\beta)},
\end{equation}
which is precisely equal to $M_{|\gamma|,\beta}=\frac{\beta}{1-|\gamma|(1-\beta)}$ when $\gamma>0$ and thus when $|\gamma|=\gamma$. 

Now, suppose that $\beta\neq 1$ and $\gamma<0$. In this case, Bob assigns a higher likelihood for the function to be in the sector that carries causal influence (i.e., to $\Id$ or $\F$) if he obtains $Y=1$ than if he obtains $Y=0$. Thus, in this case we want to know $\PP_{B|Y}(\Id \cup \F|1)$. With analogous calculations to the ones described above, we learn that
\begin{equation}
    \PP_{B|Y}(\Id \cup \F|1)=\frac{\beta}{1+\gamma(1-\beta)},
\end{equation}
which is precisely equal to $M_{|\gamma|,\beta}=\frac{\beta}{1-|\gamma|(1-\beta)}$ when $\gamma<0$ and thus when $|\gamma|=-\gamma$. 

Finally, in the case where $\beta=1$, Bob assigns 1 to both $\PP_{B|Y}(\Id \cup \F|1)$ and $\PP_{B|Y}(\Id \cup \F|1)$. Indeed, $M_{|\gamma|,\beta}=1$ in this case. 

\end{proof}

\subsection{Relation to the Resource Theory of Causal Influence}
\label{sec_relation_RTs}

Now that we introduced $\RT$ and $\RTKnow$, one could wonder if the former is just a special case of the latter. It is certainly true that the set of resources considered in $\RT$, i.e., functions $f:X \to Y$, can be viewed as a subset of all probability distributions over functions, more precisely, they are in one-to-one correspondence with point distributions, $[f]$. Here we show that if we restrict the resources in $\RTKnow$ to be point distributions on functions, and we do not restrict the set of free operations in any way, this leads to the same pre-order of resources as in $\RT$. Note that in this ``restricted'' $\RTKnow$ the set of free operations still consists of local pre- and post-processing correlated by a common cause, but the free operations need to map a point distribution on a function to a point distribution on a function. An example of such processing that transforms $[f]$ to $[\Rzero]$ is having $\PP_{\tau}=\frac{1}{2}[\tau_1] + \frac{1}{2}[\tau_2]$ where 
$\tau_1$ is pre-processing with $\Id$ and post-processing with $\Rzero$ while $\tau_2$ is pre-processing with $\Rone$ and post-processing with $\Rzero$.
Note that while this tells us that $[f]\to[\Rzero]$ in $\RTKnow$, it doesn't immediately tell us that $f\to\Rzero$ within $\RT$, as, in this case, the pre-processing cannot be thought of as simply being a (point distribution on a) function. It turns out, however, that $f\to\Rzero$ in fact holds in $\RT$. This simple example turns out to be generic:

\begin{prop}
Consider the restriction of the Resource Theory of Knowledge of Causal Influence $\RTKnow$ to the point distribution resources. The resulting resource theory has the same pre-order as the Resource Theory of Causal Influence $\RT$.
\end{prop}
\begin{proof}
Consider two functions $f_1:X \to Y$ and $f_2:X’ \to Y’$.  If $f_1$ can be transformed to $f_2$ with the free operations in $\RT$ (local pre- and post-processing, see Eq.~\eqref{fig:f_free}), we denote it by $f_1 \to f_2$. If a point distribution on $f_1$ can be transformed to a point distribution on $f_2$ with the free operations in $\RTKnow$ (local pre- and post-processing correlated by a common-cause, illustrated in Fig.~\ref{fig:common_cause_comb}), we denote it by $[f_1]\to[f_2]$.

It is easy to see that $f_1 \to f_2 \Rightarrow [f_1]\to[f_2]$. It suffices to note that

\begin{align}\label{fig:proof3}
\tikzset{every picture/.style={line width=0.75pt}}
\begin{tikzpicture}[x=0.75pt,y=0.75pt,yscale=-1,xscale=1]


\draw[line width=0.75] (94,48) -- (94,64);
\draw[line width=0.75] (94,129) -- (94,153);

\draw (77,110) rectangle (112,129);
\draw[line width=0.75] (79,64) rectangle (110,93);
\draw (77,28) rectangle (112,48);

\draw[line width=0.75] (94,12) -- (94,28);
\draw[line width=0.75] (94,93) -- (94,110);

\draw[line width=0.75] (181,39) -- (181,64);
\draw[line width=0.75] (181,93) -- (181,118);
\draw[line width=0.75] (166,64) rectangle (197,93);


\draw[line width=0.75] (371,93) -- (371,120);      
\draw (354,120) rectangle (389,140);               
\draw[line width=0.75] (360,63) rectangle (383,92);
\draw (354,20) rectangle (388,40);                 

\draw[line width=0.75] (372,140) -- (372,165);     
\draw[line width=0.75] (371,0) -- (371,20);       
\draw[line width=0.75] (371,40) -- (371,63);       

\draw[line width=0.75] (350,80) -- (360,80);
\draw[line width=0.75] (310,80) -- (350,55) -- (350,105) -- cycle;
\draw (322,73) node[anchor=north west,inner sep=0.75pt] {$[f_1]$};

\draw[line width=0.75] (344,30) -- (354,30);              
\draw[line width=0.75] (320,30) -- (344,10) -- (344,50) -- cycle;
\draw (326,23) node[anchor=north west,inner sep=0.75pt] {$[h]$};

\draw[line width=0.75] (344,130) -- (354,130);            
\draw[line width=0.75] (320,130) -- (344,110) -- (344,150) -- cycle;
\draw (326,123) node[anchor=north west,inner sep=0.75pt] {$[g]$};

\draw [line width=0.75] (447,80) -- (457,80);
\draw [line width=0.75] (407,80) -- (447,55) -- (447,105) -- cycle;
\draw (420,73) node[anchor=north west,inner sep=0.75pt] {$[f_2]$};

\draw[line width=0.75] (468,39) -- (468,64);
\draw[line width=0.75] (468,93) -- (468,118);
\draw[line width=0.75] (457,63) rectangle (479,92);

\draw (390,78) node[anchor=north west,inner sep=0.75pt] {$=$};
\draw (446,100) node[anchor=north west,inner sep=0.75pt] {$X'$};
\draw (446,43)  node[anchor=north west,inner sep=0.75pt] {$Y'$};


\draw (127,78) node[anchor=north west,inner sep=0.75pt] {$=$};

\draw (79,95) node[anchor=north west,inner sep=0.75pt] {$X$};
\draw (79,49) node[anchor=north west,inner sep=0.75pt] {$Y$};
\draw (75,134) node[anchor=north west,inner sep=0.75pt] {$X'$};
\draw (75,8)  node[anchor=north west,inner sep=0.75pt] {$Y'$};

\draw (85,70)  node[anchor=north west,inner sep=0.75pt] {$f_1$};
\draw (87,114) node[anchor=north west,inner sep=0.75pt] {$g$};
\draw (87,30)  node[anchor=north west,inner sep=0.75pt] {$h$};

\draw (171,71) node[anchor=north west,inner sep=0.75pt] {$f_2$};
\draw (159,100) node[anchor=north west,inner sep=0.75pt] {$X'$};
\draw (160,44)  node[anchor=north west,inner sep=0.75pt] {$Y'$};

\draw (240,78) node[anchor=north west,inner sep=0.75pt] {$\iff$};

\draw (356,95)  node[anchor=north west,inner sep=0.75pt] {$X$};
\draw (356,48)  node[anchor=north west,inner sep=0.75pt] {$Y$};
\draw (351,145) node[anchor=north west,inner sep=0.75pt] {$X'$};
\draw (352,4)   node[anchor=north west,inner sep=0.75pt] {$Y'$};

\end{tikzpicture}.
\end{align}
which given the free operation that transforms $f_1$ to $f_2$ explicitly shows how to construct a free operation that transforms $[f_1]$ to $[f_2]$. 

We will now show that $[f_1]\to[f_2] \Rightarrow f_1 \to f_2$. Consider a general free operation in $\RTKnow$, which consists of pre- and post-processings correlated by a common cause, 
which can be equivalently written as a convex combination of pre- and post-processings.
Hence, a free operation that transforms $[f_1]$ to $[f_2]$ can be expressed as 
\begin{align}\label{fig:proof1}
[f_2] = \sum_i p(i)  [h_i] \circ [f_1] \circ [g_i].
\end{align}
Now note that a convex combination over possible local processings gives a point distribution on a function, $[f_2]$, if and only if each element of this convex combination individually equals to $[f_2]$. This follows from the fact that point distributions on functions are the extreme points in the simplex representing probability distributions over functions. That is, Eq.~\eqref{fig:proof1} implies that
\begin{align}
    \forall i, \quad [f_2] = [h_i] \circ [f_1] \circ [g_i].
\end{align}
To complete the proof, notice that we have the equivalence given by Eq.~\eqref{fig:proof3}.

This shows that the existence of a free operation in $\RTKnow$ that transforms $[f_1]$ to $[f_2]$ implies the existence of a free operation in the $\RT$ that transforms $f_1$ to $f_2$.

In summary, we showed that $[f_1]\to[f_2]$ if and only if $f_1 \to f_2$. Therefore, the two resource theories have the same pre-order.

\end{proof}

\subsection{Monotones beyond the bit-to-bit scenario}\label{se:n2n}

In this section, we continue studying the causal structure where we have two observed variables $X$ and $Y$, and causal influence goes from $X$ to $Y$. However, here we take the cardinality of each of these variables to be an arbitrary finite integer.

Given a resource $\PP_{F}$, let us define the total probability assigned to functions with image size $k$ by $\beta_k$, with $k = 1,\ldots,n$, where $n \coloneq |Y|$. Notice that, by normalization, one has $\sum _{k=1}^n\beta_k=1$. 

Theorem~\ref{thm:betaM} proves that a set of quantities involving those $\{\beta_k\}_k$ are all monotones of $\RTKnow$. Note, however, that these are \emph{not} a complete set of monotones: for example, in the case of bit-to-bit functions, the construction of Theorem~\ref{thm:betaM}  only gives us $M_\beta$. Before presenting Theorem~\ref{thm:betaM}, we present an useful lemma:

\begin{lem}
It is possible to increase the weight on lower image-size sectors (other than just the free sector) with free operations.
\label{lemma_example}
\end{lem}
\begin{proof}
We prove this by presenting an example.

Let $|X|=|Y|=3$. Define three functions $f_{1}:X\to Y$, $f_{2}:X\to Y$ and $f_{3}:X\to Y$ as
\begin{align}
f_1(x) =
\begin{cases}
0, & \text{if } x = 0, 1, \\
1, & \text{if } x = 2.
\end{cases}
\end{align}
\begin{align}
f_2(x) =
\begin{cases}
0, & \text{if } x = 0, 1, \\
2, & \text{if } x = 2.
\end{cases}
\end{align}
\begin{align}
f_3(x) =
\begin{cases}
0, & \text{if } x = 0, \\
2, & \text{if } x = 1,2.
\end{cases}
\end{align}
Note that all functions defined above have the same image size, $|\Imm(f_1)|=|\Imm(f_2)|=|\Imm(f_3)|=2$. Consider the probability distribution over functions 
\begin{align}
    \PP_F=\frac{1}{3} [\Id] + \frac{2}{3} \left(\,\frac{1}{2}[f_1]+\frac{1}{2}[f_2]\,\right).
\end{align}
One possible processing of $\PP_F$ is pre-processing with $f_2$ and post-processing with $f_3$:
\begin{align}
f_3 \circ \PP_F \circ f_2 = f_2
\end{align}
which serves as an example that it is possible to obtain certainty about a nonfree function by decreasing the image size of the initial resource. In this example, the weight on the sector characterized by $|\Imm(f)|=2$ is increased, while the weight on the sector characterized by $|\Imm(f)|=3$ is decreased.
    
\end{proof}

\begin{thm}\label{thm:betaM}
Consider resources in a $X\to Y$ scenario, and the figures of merit $\{\beta_k\}_k$ defined above. Then, the following functions are resource monotones: 
\begin{align}
M_n[\PP_{F}] &:= \beta_n\,, \\
M_{n-1}[\PP_{F}] &:= \beta_n+ \beta_{n-1}\,,\\
&\vdots\\
M_2[\PP_{F}] &:= \beta_n+ \beta_{n-1} + ... +\beta_2\,,\\
M_1[\PP_{F}] &:= \beta_n+ \beta_{n-1} + ... +\beta_2 + \beta_1 = 1\,.
\end{align}

\end{thm}
 \begin{proof}
First notice that pre- and post- processing with any function cannot increase the image size of a function. Hence, no extremal free operations can increase any $\beta_i$ by decreasing a $\beta_j$ with $j<i$. Rather, the only way to increase a $\beta_i$ by an extremal free operation is by decreasing some collection of $\{\beta_j\}_j$ with $j>i$ (by Lemma~\ref{lemma_example}).  Because of normalization, an increase of $\beta_i$ can only happen at the cost of an \emph{equal} decrease of $\sum_{j>i}\beta_j$. This makes the value of $M_k$ either invariant (for $k\leq i$) or lower (for $k>i$). 

To see that free operations that are not extremal also cannot increase the $\{M_k\}_{k=1,\ldots,n}$, it suffices to note two things: i) the $\beta_i$ obtained by applying a convex combination of  extremal free operations to $\PP_F$ is the corresponding convex combination of the $\beta_i$'s obtained by applying each extremal free operation to $\PP_F$, and ii) each $M_k$ is a linear function of the $\{\beta_k\}_k$. 
 \end{proof}

Note that the monotones $\{M_k\}_{k=1,\ldots,n}$  are \emph{not} Ky Fan $k$-norms, despite the superficial similarity, since a Ky Fan $k$-norm of a vector is the sum of the $k$ {\em largest} components of a vector. Hence, there is no immediate connection between these monotones and majorization theory.

\section{Discussion}

\subsection{Connection with the Average Causal Effect}

In this paper, we have developed resource-theoretic methods for quantifying causation, where the resources under consideration are taken to be probability distributions over functional dependencies. Here, we denoted these by $\PP_{F} \in \mathcal{D}(F(X\to Y))$. Because most prior works in the literature do not study probability distributions over functional dependencies, but rather conditional probability distributions, one might have expected that measures of causal influence would be a function of the latter rather than the former. For example, the quantity known as the Average Causal Effect (ACE)~\cite{Holland1988,pearl2009causality}  is a function that assigns a real number to every conditional probability distribution, and which is sometimes interpreted as a measure of causal influence.\footnote{It is also sometimes described rather as a measure of how different an outcome would be under two different treatments, which we believe to be a more appropriate description.} We dispute this interpretation, and we further show the precise relation between ACE and measures of causal inference in the special case where $X$ and $Y$ are binary: the ACE is a lower bound on the probability of causal connection for a given conditional probability distribution $\PP_{Y|X}$. Here, we will also refer to such conditional probability distributions as \emph{stochastic maps}. 

\begin{figure}[h]
    \centering
    \tikzset{every picture/.style={line width=0.75pt}}
\begin{tikzpicture}[x=0.75pt,y=0.75pt,yscale=-1,xscale=1]

\draw[color={rgb, 255:red, 65; green, 117; blue, 5},draw opacity=1,line width=1.5] (207,86) -- (207,105);
\draw[color={rgb, 255:red, 74; green, 144; blue, 226},draw opacity=1,line width=1.5] (155,63) -- (155,144);

\draw[color={rgb, 255:red, 74; green, 144; blue, 226},draw opacity=1,line width=1.5,dash pattern={on 5.63pt off 4.5pt}] (155,144) -- (155,293);

\draw (120,250) -- (250,250);

\draw (50,340) -- (120,250);

\draw (50,340) -- (180,340);

\draw (180,340) -- (250,250);

\draw[color={rgb, 255:red, 155; green, 155; blue, 155},draw opacity=1,dash pattern={on 4.5pt off 4.5pt}] (50,190) -- (50,340);

\draw[color={rgb, 255:red, 155; green, 155; blue, 155},draw opacity=1,dash pattern={on 4.5pt off 4.5pt}] (120,28) -- (120,250);

\draw[color={rgb, 255:red, 155; green, 155; blue, 155},draw opacity=1,dash pattern={on 4.5pt off 4.5pt}] (180,90) -- (180,340);

\draw[color={rgb, 255:red, 155; green, 155; blue, 155},draw opacity=1,dash pattern={on 4.5pt off 4.5pt}] (250,100) -- (250,250);

\draw (120,28) -- (250,100);

\draw (50,190) -- (120,28);

\draw (50,190) -- (180,90);

\draw (180,90) -- (250,100);

\draw (120,28) -- (180,90);

\draw (50,190) -- (250,100);


\draw[shift={(155,293)}, rotate = 90] [color={rgb, 255:red, 74; green, 144; blue, 226},draw opacity=1,fill={rgb, 255:red, 74; green, 144; blue, 226},fill opacity=1,line width=1.5] (0,0) circle [x radius=3.05, y radius=3.05];
\draw[shift={(155,63)}, rotate = 90] [color={rgb, 255:red, 74; green, 144; blue, 226},draw opacity=1,fill={rgb, 255:red, 74; green, 144; blue, 226},fill opacity=1,line width=1.5] (0,0) circle [x radius=3.05, y radius=3.05];

\draw[shift={(155,144)}, rotate = 90] [color={rgb, 255:red, 74; green, 144; blue, 226},draw opacity=1,fill={rgb, 255:red, 74; green, 144; blue, 226},fill opacity=1,line width=1.5] (0,0) circle [x radius=3.05, y radius=3.05];

\draw[color={rgb, 255:red, 65; green, 117; blue, 5},draw opacity=1,line width=1.5,dash pattern={on 5.63pt off 4.5pt}] (207,105) -- (207,290);


\draw[shift={(207,290)}, rotate = 90] [color={rgb, 255:red, 65; green, 117; blue, 5},draw opacity=1,fill={rgb, 255:red, 65; green, 117; blue, 5},fill opacity=1,line width=1.5] (0,0) circle [x radius=3.05, y radius=3.05];
\draw[shift={(207,86)}, rotate = 90] [color={rgb, 255:red, 65; green, 117; blue, 5},draw opacity=1,fill={rgb, 255:red, 65; green, 117; blue, 5},fill opacity=1,line width=1.5] (0,0) circle [x radius=3.05, y radius=3.05];

\draw[shift={(207,105)}, rotate = 90] [color={rgb, 255:red, 65; green, 117; blue, 5},draw opacity=1,fill={rgb, 255:red, 65; green, 117; blue, 5},fill opacity=1,line width=1.5] (0,0) circle [x radius=3.05, y radius=3.05];

\draw (20,185) node[anchor=north west][inner sep=0.75pt] {$[\Rone]$};

\draw (105,8) node[anchor=north west][inner sep=0.75pt] {$[\Id]$};

\draw (175,70) node[anchor=north west][inner sep=0.75pt] {$[\F]$};

\draw (258,90) node[anchor=north west][inner sep=0.75pt] {$[\Rzero]$};

\draw (-120,80) node[anchor=north west][inner sep=0.75pt] {$\boldsymbol{{ \mathcal{D}(F(X\to Y))}}$};
\draw (-100,300) node[anchor=north west][inner sep=0.75pt] {$\boldsymbol{{ \Sigma(X\to Y)}}$};
\draw (-65,200) node[anchor=north west][inner sep=0.75pt] {$\boldsymbol{{\Gamma}}$};
\draw[->, thick, black] (-70,100) -- (-70,290);

\draw (13,330) node[anchor=north west][inner sep=0.75pt] {$\scalebox{0.65}{$\begin{pmatrix} 0 & 0 \\ 1 & 1 \end{pmatrix}$}$};

\draw (185,330) node[anchor=north west][inner sep=0.75pt] {$\scalebox{0.65}{$\begin{pmatrix} 0 & 1 \\ 1 & 0 \end{pmatrix}$}$};

\draw (250,220) node[anchor=north west][inner sep=0.75pt] {$\scalebox{0.65}{$\begin{pmatrix} 1 & 1 \\ 0 & 0 \end{pmatrix}$}$};

\draw (120,220) node[anchor=north west][inner sep=0.75pt] {$\scalebox{0.65}{$\begin{pmatrix} 1 & 0 \\ 0 & 1 \end{pmatrix}$}$};

\end{tikzpicture}
    \caption{Geometrical description of the map $\Gamma[\PP_{F}]$ as a projection of the tetrahedron of probability distributions over functions $\mathcal{D}(F(X\to Y))$ down to the square of stochastic maps $\Sigma(X \to Y)$ for binary $X$ and $Y$.  
 }
    \label{fig:simplex_projection}
\end{figure}
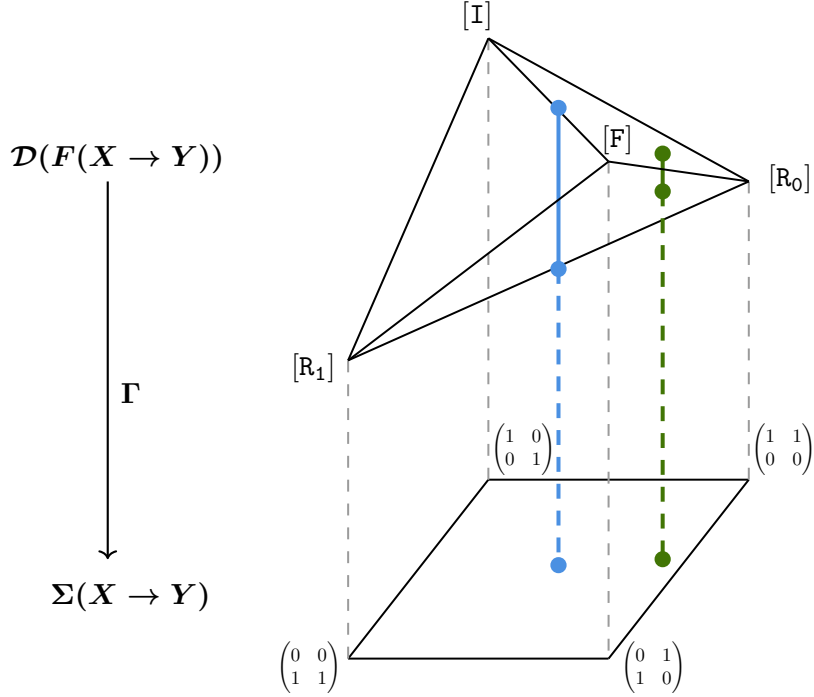

Let us first relate stochastic maps to probability distributions over functions. Denote the space of stochastic maps $\PP_{Y|X}:X\to Y$ by $\PP_{Y|X}\in\Sigma(X \to Y).$ 
Recall that there is a straightforward way to map a probability distribution over functions (an element of $\mathcal{D}(F(X\to Y))$) to a  conditional probability distribution (an element of $\Sigma(X\to Y)$), which is the following
\begin{align}
\Gamma[\PP_{F}]:=\sum_{f\in F(X\to Y)} \PP_{F}(f)\,\delta_{Y,f(X)}
\end{align}
for any $\PP_{F}=\sum_{f\in F(X\to Y)}\PP_{F}(f) [f]$. Note that $\Gamma:\mathcal{D}(F(X\to Y))\to \Sigma(X\to Y)$ is always surjective, and is injective if and only if either $X$ or $Y$ is the singleton set. 
If $X$ and $Y$ are binary variables, we can geometrically picture this map as a projection of the tetrahedron of probability distributions over functions down to the square of conditional probability distributions (or equivalently, stochastic maps), as shown in Fig.~\ref{fig:simplex_projection}.   One can also view the action of $\Gamma$ as
the linear quotient map relative  to the equivalence relation generated by the relation $\frac{1}{2} [\Id] + \frac{1}{2} [\F]\simeq \frac{1}{2}[\Rzero]+\frac{1}{2}[\Rone]$.

It is also useful to look at the preimage map associated to $\Gamma$, namely a map that takes a conditional probability distribution to a set of probability distributions over functions. We denote this map by $\bar{\Gamma}:\Sigma(X\to Y) \to \mathcal{P}(\mathcal{D}(F(X\to Y)))$ where $\mathcal{P}(Z)$ denotes the powerset of $Z$. It is defined by
\begin{align}
\bar{\Gamma}(\PP_{Y|X}):=\{\PP_{F}\ | \ \Gamma[\PP_{F}]=\PP_{Y|X}\}
\end{align}
for all $\PP_{Y|X}\in\Sigma(X\to Y)$. 
In the case where both $X$ and $Y$ are binary variables, we can geometrically picture this as mapping any given $\PP_{Y|X}$ to a set containing the points in the intersection of the vertical line through $\PP_{Y|X}$ with the tetrahedron. For example, in Fig.~\ref{fig:simplex_projection} the green (blue) dot is mapped to the set of probability distributions over functions corresponding to the points in the green (blue) line segment. 

It is clear that the ACE, as well as any other quantity that depends only on the stochastic map, is not a measure of causal influence. For example, consider the conditional probability distribution $\PP_{Y|X}$ from Eq.~\eqref{eq_randomizing_channel}, which is represented by the blue dot in the square in Fig.~\ref{fig:simplex_projection}. Then, any probability distributions over functions in the blue line segment above will give the same value of ACE. However, the top dot of the line segment, which corresponds to the equal mixture of $\Id$ and $\F$ (i.e., $\PP_B^1$ as defined in Eq.~\eqref{eq_Id_Flip}), has maximal probability of causal connection since its $M_{\beta}$ is 1, while the bottom dot of the line segment, which correspond to the equal mixture of $\Rzero$ and $\Rone$ (i.e., $\PP_B^2$ as defined in Eq.~\eqref{eq_Rzero_Rone}), has no causal connection whatsoever.
Despite this fact, we will now show that the ACE is closely connected to $M_\beta$.

For the case we are considering here, i.e., a direct causal influence from a binary variable $X$ to a binary variable $Y$ with no confounders, the ACE is defined as\footnote{The general definition of the ACE employs the do-conditional $\PP_{Y|do(X=x)}$, which denotes the distribution of $Y$ given that $X$ was {\em set} to a value of $X=x$~\cite{Holland1988,pearl2009causality}. For the case of direct causal influence of $X$ on $Y$ with no confounding, $\PP_{Y|do(X=x)}=\PP_{Y|X}$. When $Y$ is not a binary variable, the ACE is defined as $\textrm{ACE} = \mathbb{E}[\PP_{Y|X=1}-\PP_{Y|X=0}]$.}
\begin{equation}\label{eq:ACE}
\textrm{ACE}(\PP_{Y|X}) = \PP_{Y|X}(1|1)-\PP_{Y|X}(1|0).
\end{equation}
Any functional on stochastic maps (that is, a function $C:\Sigma(X\to Y)\to \mathds{R}$) uniquely defines a functional on probability distributions over functions via $C_{\text{dist}}:=C\circ \Gamma$. Indeed, when $C$ is ACE, it induces the following definition on probability distributions:
\begin{align}
    {\rm ACE}_{\text{dist}}(\PP_B) \coloneq {\rm ACE}\big(\Gamma (\PP_B) \big).
\end{align}
Since we have that
\begin{align}
    \PP_{Y|X}(Y=1|X=1) = \PP_{B}(\Id)+ \PP_{B}(\Rone), \quad \PP_{Y|X}(Y=1|X=0) = \PP_{B}(\F)+ \PP_{B}(\Rone), 
\end{align}
we conclude 
\begin{align}\label{eq:ACEIF}
{\rm ACE}_{\text{dist}}(\PP_B) &= \PP_{B}(\Id)+ \PP_{B}(\Rone) - (\PP_{B}(\F) + \PP_{B}(\Rone)) \\ \nonumber
&= \PP_{B}(\Id) -  \PP_{B}(\F) . 
\end{align}
This leads us to the first precise connection between ACE and $\RTKnow$:
\begin{prop}\label{ACEmono}
Consider a causal structure with two binary observed  variables $X$ and $Y$ and no confounders, where $X$ is a parent of $Y$. $ {\rm ACE}_{\text{dist}}(\PP_B)$ is a monotone in $\RTKnow$.
\end{prop}

\begin{proof}
Given the parametrization of $\PP_B$ as in Eq.~\eqref{generalparametrization}, we know that
\begin{align}
   \PP_B (\Id)= \beta\frac{1-\alpha}{2}, \quad \PP_B(\F)=\beta\frac{1+\alpha}{2}.
\end{align}
From Eq.~\eqref{eq:ACEIF}, we know that $ {\rm ACE}_{\text{dist}}(\PP_B) = \PP_B (\Id) - \PP_B (\F)$. Thus, 
\begin{align} \label{eq:ACEba}
    | {\rm ACE}_{\text{dist}}(\PP_B) |=\beta|\alpha|.
\end{align}
Both $\beta$ and $|\alpha|$ are monotones in $\RTKnow$ and both are nonnegative valued, thus their product is also a monotone. Hence $| {\rm ACE}_{\text{dist}}(\PP_B) |$ is a monotone in $\RTKnow$.
\end{proof}
Proposition~\ref{ACEmono} tells us that even though ACE does not quantify the strength of causal influence (which is given by $M_{\beta}=\beta$ in $\RTKnow$ for binary variables), the function induced by it on $\PP_B$ is related to $M_{\beta}$. More precisely, it provides some coarse-grained information about $M_{\beta}$, where the coarse-graining comes from one's uncertainty on $\Id$ versus $\F$ (quantified by $|\alpha|$).

From  Eq.~\eqref{eq:ACEba}, we can derive an interpretation of $|\mathrm{ACE}(\PP_{Y|X})|$ for the case of binary variables: it is the lower bound on the strength of causal influence needed to realize the given stochastic map $\PP_{Y|X}$.
\begin{prop}\label{prop:ACE}
Consider a causal structure with two binary observed variables $X$ and $Y$ and no confounders, where $X$ is a parent of $Y$. Then, 
\begin{align}
  |\mathrm{ACE}(\PP_{Y|X})|=  \operatorname*{min}\limits_{\bar{\Gamma}[\PP_{Y|X}]}
      \Big\{ M_\beta(\mathbb{P}_B) \Big\}.
\end{align}
That is, 
\begin{align}
 \forall \PP_B\in \bar{\Gamma}(\PP_{Y|X}), \, M_{\beta}(\PP_B) \geq|\mathrm{ACE}(\PP_{Y|X})|, \text{ and the bound is tight}.
\end{align}
\end{prop}

\begin{proof}

The ACE depends only on the stochastic map, and is therefore invariant under the choice of $\PP_B$ within a given preimage $\bar\Gamma(\PP_{Y|X})$. Since ${\rm ACE}_{\text{dist}}(\PP_{B}) =\PP_{B}(\Id) -  \PP_{B}(\F)$, when we vary $\PP_B$ within a given preimage $\bar\Gamma(\PP_{Y|X})$, any decrease in $\PP_B(\Id)$ must be accompanied by an equal decrease in $\PP_B(\F)$. Therefore, the minimum $M_\beta$ among probability distributions over functions that belong to $\bar{\Gamma}[\PP_{Y|X}]$ is achieved when  $\PP_B(\Id)$ and  $\PP_B(\F)$ are decreased as far as possible, i.e.,  when one of them is zero. Geometrically, this minimum is achieved by the bottom point of the line segment defined by $\bar\Gamma(\PP_{Y|X})$, such as the bottom blue and green points of the tetrahedron of Fig.~\ref{fig:simplex_projection}. 

Denote this bottom point by $\PP_{B{\rm min}}$, then $\forall \PP_B\in \bar{\Gamma}(\PP_{Y|X}), \, M_{\beta}(\PP_B)\geq M_\beta(\PP_{B{\rm min}})$.
Since the bottom point $\PP_{B{\rm min}}$ lies on one of the bottom two faces of the tetrahedron, corresponding to the case where $|\alpha|=1$, Eq.~\eqref{eq:ACEba} becomes $|{\rm ACE}_{\text{dist}}(\PP_{B{\rm min}})|=\beta=M_\beta(\PP_{B{\rm min}})$. 
Furthermore, since  $\Gamma(\PP_{B{\rm min}})=\PP_{Y|X}$, we also have that $|{\rm ACE}(\PP_{Y|X})|=|{\rm ACE}_{\text{dist}}(\PP_{B{\rm min}})|$. Thus, $\forall \PP_B\in \bar{\Gamma}(\PP_{Y|X}), \, M_{\beta}(\PP_B)\geq|{\rm ACE}(\PP_{Y|X})|$ and the equal sign is achieved by $\PP_{B{\rm min}}$.
\end{proof}

The above shows how, at least in the context of the DAG $X\to Y$ with binary variables, the ACE can be derived from a rigorous resource-theoretic perspective and related rigorously to strength of causal influence, rather than merely being viewed as a measure of observational differences that arise from different possible values for the cause. We stress however that the ACE lacks information when compared to quantifiers of causal influence based on probability distributions over functions, that are derived from our resource theories, such as $\beta$ and $|\alpha|$. In particular, in the binary case the absolute value of the ACE only gives the \emph{least} amount of causal influence necessary to realize a certain conditional probability distribution, and provides coarse-grained information for $\beta$ and $|\alpha|$.  This deficiency is found not only for the ACE, but also for any other purported quantifiers of causal influence based on conditional probability distributions, such as the ones discussed in Ref.~\cite{Janzing2013}.

\subsection{The relation to Shannon theory}\label{sec_Shannon}

Shannon theory~\cite{Shannon}, herein denoted by $\Shannon$,  can be formulated as a resource theory of one-way classical communication channels (see Example 2.6 in Ref.~\cite{coecke2016mathematical}). The basic scenario involves a sender who wishes to transmit a message $X$ to a receiver, who in turn receives an output $Y$ from the channel. 
In the language of Section~\ref{sec_flag}, here Alice and Bob do not have access to any information about the flag variable that encodes the function applied by the channel,  and consequently, the channel is modeled by the conditional probability distribution $\PP_{Y|X}$ one obtains by marginalizing over the flag variable.  Note that we here consider the version of Shannon theory wherein the sender and receiver are assumed to have access to shared randomness (which can help in pre- and post- processing the channel, but that does not influence the input $X$).  The free resources are taken to be the conditionals $\PP_{Y|X}$ such that $\PP_{Y|X}= \PP_{Y}$, that is, channels that do not allow any communication. 
The free operations on channels  are given by encodings and decodings of the message, represented as stochastic maps that may be correlated by shared randomness. The aim of this theory is to find the ordering over channels, since the higher a channel is in the order, the more communication can be achieved with it.

When Alice and Bob {\em do} have access to partial information about the flag variable — i.e., they can characterize the channel with a probability distribution over functions rather then just by a stochastic map — this information cannot be captured within $\Shannon$. Consider again the bit-to-bit case we explored in Sec.~\ref{sec_flag} and studied in detail in Sec.~\ref{sec:b2bc}. In $\RTKnow$, resources may be nonfree even when they do not enable communication. For instance, $\PP_{B}^{1}=\frac{1}{2}[\Id]+\frac{1}{2}[\F]$ is a nonfree resource because it does transmit causal influence. It can  only   be used for communication,   however,  if the value of the flag variable is known. In contrast, $\PP_{B}^{2}=\frac{1}{2}[\Rzero]+\frac{1}{2}[\Rone]$ is a free resource as it erases all information about the input variable. Because $\PP_{B}^{1}$ and $\PP_{B}^{2}$ induce the same stochastic map,  they correspond to the same resource in $\Shannon$, indeed a {\em free} resource. 

This distinction highlights how  $\RTKnow$  concerns knowledge of causal influence  
rather than simply the possibility of 
reliable message transmission.   The difference is that there are many distributions over causal dependences that imply the inability to transmit information. Consider the situation wherein  the distribution over functions is $\PP_{B}^{1}$ (the equal mixture of $\Id$ and $\F$).  In Shannon theory, where one does not incorporate any knowledge about the flag variable into the resource description, there is no possibility of transmitting information in this situation and so the channel is considered to be resourceless,  By contrast, in  $\RTKnow$, the same situation is modeled as a nonfree resource.  This can be seen as encoding the fact that if one were able to learn the value of the flag variable, it would become  possible to transmit some information through the channel.

\subsection{The categorical perspective}\label{sec_categorical}

Resource theories are often defined via partitioned process theories \cite{coecke2016mathematical}. A process theory can be formally viewed as a symmetric monoidal category $\mathcal{C}$, and a partitioned process theory is a symmetric monoidal category with a specified free symmetric monoidal subcategory $\mathcal{C}_{\text{free}}\subseteq \mathcal{C}$. 

For us, we are not interested in arbitrary SMCs $\mathcal{C}$, but specifically SMCs where the morphisms are combs in some other category, that is, $\mathcal{C}=\mathbf{Comb}[\mathcal{D}]$ for some $\mathcal{D}$. Specifically, objects in $\mathbf{Comb}[\mathcal{D}]$ are given by pairs of objects in $\mathcal{D}$, e.g., $(X,Y)$ and then morphisms from $(X,Y)\to (X',Y')$ are defined as triples $(Z,f:X'\to X\otimes Z, g:Y\otimes Z \to Y')$ where $Z$ is an object in $\mathcal{D}$ and $f$ and $g$ are morphisms in $\mathcal{D}$\footnote{Formally it is better to think about combs as a suitable quotient of these triples, but this isn't necessary for us here.}. That is, morphisms in $\mathcal{C}$ are combs in $\mathcal{D}$ such as 
\begin{align}\label{fig:cat1}
\tikzset{every picture/.style={line width=0.75pt}}
\begin{tikzpicture}[x=0.75pt,y=0.75pt,yscale=-1,xscale=1]

\draw [line width=0.75]    (25,125) -- (25,150) ;
\draw [line width=0.75]    (36,170) -- (36,193) ;
\draw  [line width=0.75]  (11,150) -- (61,150) -- (61,170) -- (11,170) -- cycle ;
\draw [line width=0.75]    (49,80) -- (49,150) ;
\draw [line width=0.75]    (25,80) -- (25,105) ;
\draw [line width=0.75]    (35,35) -- (35,60) ;
\draw  [line width=0.75]  (11,60) -- (61,60) -- (61,80) -- (11,80) -- cycle ;

\draw (55,110) node [anchor=north west][inner sep=0.75pt]    {$Z$};
\draw (8,131) node [anchor=north west][inner sep=0.75pt]    {$X$};
\draw (7,87) node [anchor=north west][inner sep=0.75pt]    {$Y$};
\draw (13,40) node [anchor=north west][inner sep=0.75pt]    {$Y'$};
\draw (14,176) node [anchor=north west][inner sep=0.75pt]    {$X'$};
\draw (29,153) node [anchor=north west][inner sep=0.75pt]    {$f$};
\draw (31,65) node [anchor=north west][inner sep=0.75pt]    {$g$};

\end{tikzpicture},
\end{align} 
where sequential and parallel composition of the morphisms in $\mathcal{D}$ are given by composition of the underlying combs as illustrated in Fig.~\ref{fig:cat2}.
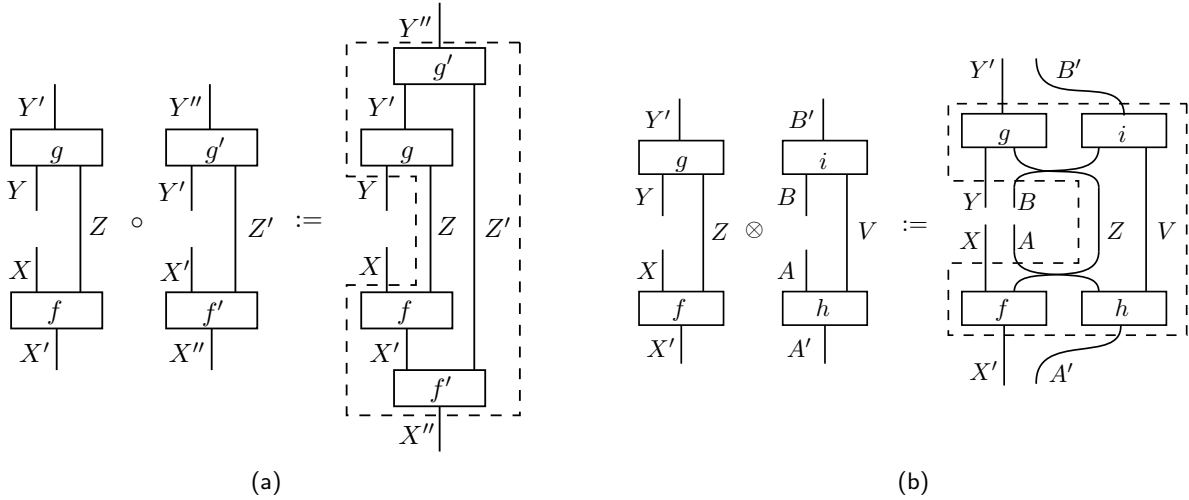
\begin{figure}[ht]
\centering
\begin{subfigure}[t]{0.45\textwidth}
  \centering
  \resizebox{\linewidth}{!}{%
    \tikzset{every picture/.style={line width=0.75pt}}

\begin{tikzpicture}[x=0.75pt,y=0.75pt,yscale=-1,xscale=1]

\draw [line width=0.75]    (19,129) -- (19,154) ;
\draw [line width=0.75]    (30,174) -- (30,197) ;
\draw  [line width=0.75]  (5,154) -- (55,154) -- (55,174) -- (5,174) -- cycle ;
\draw [line width=0.75]    (43,84) -- (43,154) ;
\draw [line width=0.75]    (19,84) -- (19,109) ;
\draw [line width=0.75]    (29,39) -- (29,64) ;
\draw  [line width=0.75]  (5,64) -- (55,64) -- (55,84) -- (5,84) -- cycle ;
\draw [line width=0.75]    (104,129) -- (104,154) ;
\draw [line width=0.75]    (115,174) -- (115,197) ;
\draw  [line width=0.75]  (90,154) -- (140,154) -- (140,174) -- (90,174) -- cycle ;
\draw [line width=0.75]    (128,84) -- (128,154) ;
\draw [line width=0.75]    (104,84) -- (104,109) ;
\draw [line width=0.75]    (114,39) -- (114,64) ;
\draw  [line width=0.75]  (90,64) -- (140,64) -- (140,84) -- (90,84) -- cycle ;
\draw [line width=0.75]    (211,129) -- (211,154) ;
\draw [line width=0.75]    (222,174) -- (222,197) ;
\draw  [line width=0.75]  (197,154) -- (247,154) -- (247,174) -- (197,174) -- cycle ;
\draw [line width=0.75]    (235,84) -- (235,154) ;
\draw [line width=0.75]    (211,84) -- (211,109) ;
\draw [line width=0.75]    (221,39) -- (221,64) ;
\draw  [line width=0.75]  (197,64) -- (247,64) -- (247,84) -- (197,84) -- cycle ;
\draw  [line width=0.75]  (215,19) -- (265,19) -- (265,39) -- (215,39) -- cycle ;
\draw  [line width=0.75]  (215,197) -- (265,197) -- (265,217) -- (215,217) -- cycle ;
\draw [line width=0.75]    (240,217) -- (240,242) ;
\draw [line width=0.75]    (240,-6) -- (240,19) ;
\draw [line width=0.75]    (259,39) -- (259,197) ;
\draw  [dash pattern={on 4.5pt off 4.5pt}][line width=0.75]  (284,16) -- (284,222) -- (189,222) -- (189,150) -- (227,150) -- (227,90) -- (189,90) -- (189,16) -- (284,16) -- cycle ;

\draw (2,135) node [anchor=north west][inner sep=0.75pt]    {$X$};
\draw (1,91) node [anchor=north west][inner sep=0.75pt]    {$Y$};
\draw (7,44) node [anchor=north west][inner sep=0.75pt]    {$Y'$};
\draw (8,180) node [anchor=north west][inner sep=0.75pt]    {$X'$};
\draw (23,157) node [anchor=north west][inner sep=0.75pt]    {$f$};
\draw (25,72) node [anchor=north west][inner sep=0.75pt]    {$g$};
\draw (85,135) node [anchor=north west][inner sep=0.75pt]    {$X'$};
\draw (84,91) node [anchor=north west][inner sep=0.75pt]    {$Y'$};
\draw (90,44) node [anchor=north west][inner sep=0.75pt]    {$Y''$};
\draw (90,180) node [anchor=north west][inner sep=0.75pt]    {$X''$};
\draw (108,157) node [anchor=north west][inner sep=0.75pt]    {$f'$};
\draw (110,67) node [anchor=north west][inner sep=0.75pt]    {$g'$};
\draw (194,133) node [anchor=north west][inner sep=0.75pt]    {$X$};
\draw (193,91) node [anchor=north west][inner sep=0.75pt]    {$Y$};
\draw (199,44) node [anchor=north west][inner sep=0.75pt]    {$Y'$};
\draw (200,180) node [anchor=north west][inner sep=0.75pt]    {$X'$};
\draw (215,157) node [anchor=north west][inner sep=0.75pt]    {$f$};
\draw (217,72) node [anchor=north west][inner sep=0.75pt]    {$g$};
\draw (235,21) node [anchor=north west][inner sep=0.75pt]    {$g'$};
\draw (233,199) node [anchor=north west][inner sep=0.75pt]    {$f'$};
\draw (46,111) node [anchor=north west][inner sep=0.75pt]    {$Z$};
\draw (133,111) node [anchor=north west][inner sep=0.75pt]    {$Z'$};
\draw (263,111) node [anchor=north west][inner sep=0.75pt]    {$Z'$};
\draw (237,111) node [anchor=north west][inner sep=0.75pt]    {$Z$};
\draw (216,225) node [anchor=north west][inner sep=0.75pt]    {$X''$};
\draw (216,0) node [anchor=north west][inner sep=0.75pt]    {$Y''$};
\draw (69,111) node [anchor=north west][inner sep=0.75pt]    {$\circ$};
\draw (159,110) node [anchor=north west][inner sep=0.75pt]    {$:=$};

\end{tikzpicture}
  }
  \caption{}
\end{subfigure}
\hfill
\begin{subfigure}[t]{0.48\textwidth}
  \centering
  \raisebox{.2\height}{%
  \resizebox{\linewidth}{!}{%
    \tikzset{every picture/.style={line width=0.75pt}}

\begin{tikzpicture}[x=0.75pt,y=0.75pt,yscale=-1,xscale=1]

\draw [line width=0.75]    (19,129) -- (19,154) ;
\draw [line width=0.75]    (30,174) -- (30,197) ;
\draw  [line width=0.75]  (5,154) -- (55,154) -- (55,174) -- (5,174) -- cycle ;
\draw [line width=0.75]    (43,84) -- (43,154) ;
\draw [line width=0.75]    (19,84) -- (19,109) ;
\draw [line width=0.75]    (29,39) -- (29,64) ;
\draw  [line width=0.75]  (5,64) -- (55,64) -- (55,84) -- (5,84) -- cycle ;
\draw [line width=0.75]    (104,129) -- (104,154) ;
\draw [line width=0.75]    (115,174) -- (115,197) ;
\draw  [line width=0.75]  (90,154) -- (140,154) -- (140,174) -- (90,174) -- cycle ;
\draw [line width=0.75]    (128,84) -- (128,154) ;
\draw [line width=0.75]    (104,84) -- (104,109) ;
\draw [line width=0.75]    (114,39) -- (114,64) ;
\draw  [line width=0.75]  (90,64) -- (140,64) -- (140,84) -- (90,84) -- cycle ;
\draw [line width=0.75]    (210,114) -- (210,154) ;
\draw [line width=0.75]    (221,174) -- (221,210) ;
\draw  [line width=0.75]  (196,154) -- (246,154) -- (246,174) -- (196,174) -- cycle ;
\draw [line width=0.75]    (210,68) -- (210,104) ;
\draw [line width=0.75]    (220,15) -- (220,48) ;
\draw  [line width=0.75]  (196,48) -- (246,48) -- (246,68) -- (196,68) -- cycle ;
\draw  [line width=0.75]  (267,154) -- (317,154) -- (317,174) -- (267,174) -- cycle ;
\draw [line width=0.75]    (305,68) -- (305,154) ;
\draw  [line width=0.75]  (267,48) -- (317,48) -- (317,68) -- (267,68) -- cycle ;
\draw    (227,130) .. controls (227,154) and (277,135) .. (277,154) ;
\draw    (227,68) .. controls (227,78) and (235,80) .. (244,82) .. controls (259,83) and (277,81) .. (277,92) ;
\draw    (277,130) .. controls (277,154) and (227,135) .. (227,154) ;

\draw [line width=0.75]    (277,92) -- (277,130) ;

\draw    (277,68) .. controls (277,92) and (227,73) .. (227,92) ;
\draw    (290,174) .. controls (291,198) and (239,177) .. (240,209) ;
\draw    (292,49) .. controls (292,21) and (240,49) .. (240,15) ;

\draw [line width=0.75]    (227,92) -- (227,104) ;
\draw [line width=0.75]    (227,114) -- (227,130) ;

\draw  [dash pattern={on 4.5pt off 4.5pt}][line width=0.75]  (329,42) -- (329,180) -- (188,180) -- (188,137) -- (265,137) -- (265,88) -- (188,88) -- (188,42) -- (329,42) -- cycle ;

\draw (2,135) node [anchor=north west][inner sep=0.75pt]    {$X$};
\draw (1,91) node [anchor=north west][inner sep=0.75pt]    {$Y$};
\draw (7,44) node [anchor=north west][inner sep=0.75pt]    {$Y'$};
\draw (8,180) node [anchor=north west][inner sep=0.75pt]    {$X'$};
\draw (23,157) node [anchor=north west][inner sep=0.75pt]    {$f$};
\draw (25,72) node [anchor=north west][inner sep=0.75pt]    {$g$};
\draw (85,135) node [anchor=north west][inner sep=0.75pt]    {$A$};
\draw (85,91) node [anchor=north west][inner sep=0.75pt]    {$B$};
\draw (92,44) node [anchor=north west][inner sep=0.75pt]    {$B'$};
\draw (91,180) node [anchor=north west][inner sep=0.75pt]    {$A'$};
\draw (108,157) node [anchor=north west][inner sep=0.75pt]    {$h$};
\draw (110,70) node [anchor=north west][inner sep=0.75pt]    {$i$};
\draw (46,112) node [anchor=north west][inner sep=0.75pt]    {$Z$};
\draw (133,111) node [anchor=north west][inner sep=0.75pt]    {$V$};
\draw (66,111) node [anchor=north west][inner sep=0.75pt]    {$\otimes$};
\draw (159,110) node [anchor=north west][inner sep=0.75pt]    {$:=$};
\draw (193,117) node [anchor=north west][inner sep=0.75pt]    {$X$};
\draw (195,94) node [anchor=north west][inner sep=0.75pt]    {$Y$};
\draw (199,14) node [anchor=north west][inner sep=0.75pt]    {$Y'$};
\draw (199,194) node [anchor=north west][inner sep=0.75pt]    {$X'$};
\draw (214,157) node [anchor=north west][inner sep=0.75pt]    {$f$};
\draw (216,56) node [anchor=north west][inner sep=0.75pt]    {$g$};
\draw (250,14) node [anchor=north west][inner sep=0.75pt]    {$B'$};
\draw (246,195) node [anchor=north west][inner sep=0.75pt]    {$A'$};
\draw (285,157) node [anchor=north west][inner sep=0.75pt]    {$h$};
\draw (287,54) node [anchor=north west][inner sep=0.75pt]    {$i$};
\draw (310,111) node [anchor=north west][inner sep=0.75pt]    {$V$};
\draw (280,111) node [anchor=north west][inner sep=0.75pt]    {$Z$};
\draw (227,117) node [anchor=north west][inner sep=0.75pt]    {$A$};
\draw (227,94) node [anchor=north west][inner sep=0.75pt]    {$B$};

\end{tikzpicture}
  }
    }
  \caption{}
\end{subfigure}
\caption{(a) Sequential and (b) parallel composition of the morphisms in $\mathcal{D}$.}
\label{fig:cat2}
\end{figure}

We can then define the free subtheory $\mathbf{Comb}[\mathcal{D}]_{\text{c.c.}}\subseteq \mathbf{Comb}[\mathcal{D}]$ of common cause combs in which the morphisms $(Z,f,g)$ are such that $f$  factorises as $f = (h\otimes \mathds{1}_Z) \circ (\mathds{1}_X \otimes s)$ where $s:I\to A\otimes Z$ and $h:X\otimes A\to X'$, that is, where those in which there is no direct cause from $f$ to $g$ but only a common cause given by the state $s$. Diagramamtically, this means that we restrict to combs of the form
\begin{align}\label{fig:cat3}
\tikzset{every picture/.style={line width=0.75pt}}

\begin{tikzpicture}[x=0.75pt,y=0.75pt,yscale=-1,xscale=1]

\draw [line width=0.75]    (25,125) -- (25,150) ;
\draw [line width=0.75]    (36,170) -- (36,193) ;
\draw  [line width=0.75]  (11,150) -- (61,150) -- (61,170) -- (11,170) -- cycle ;
\draw [line width=0.75]    (49,80) -- (49,150) ;
\draw [line width=0.75]    (25,80) -- (25,105) ;
\draw [line width=0.75]    (35,35) -- (35,60) ;
\draw  [line width=0.75]  (11,60) -- (61,60) -- (61,80) -- (11,80) -- cycle ;
\draw [line width=0.75]    (118,125) -- (118,150) ;
\draw [line width=0.75]    (113,170) -- (113,236) ;
\draw  [line width=0.75]  (106,150) -- (128,150) -- (128,170) -- (106,170) -- cycle ;
\draw [line width=0.75]    (150,80) -- (150,189) ;
\draw [line width=0.75]    (125,80) -- (125,105) ;
\draw [line width=0.75]    (135,35) -- (135,60) ;
\draw  [line width=0.75]  (111,60) -- (161,60) -- (161,80) -- (111,80) -- cycle ;
\draw   (137,211) -- (116,190) -- (157,189) -- cycle ;
\draw [line width=0.75]    (123,171) -- (123,191) ;
\draw  [dash pattern={on 4.5pt off 4.5pt}][line width=0.75]  (101,148) -- (165,148) -- (165,216) -- (101,216) -- cycle ;

\draw (8,131) node [anchor=north west][inner sep=0.75pt]    {$X$};
\draw (7,87) node [anchor=north west][inner sep=0.75pt]    {$Y$};
\draw (13,40) node [anchor=north west][inner sep=0.75pt]    {$Y'$};
\draw (14,176) node [anchor=north west][inner sep=0.75pt]    {$X'$};
\draw (29,153) node [anchor=north west][inner sep=0.75pt]    {$f$};
\draw (31,68) node [anchor=north west][inner sep=0.75pt]    {$g$};
\draw (54,104) node [anchor=north west][inner sep=0.75pt]    {$Z$};
\draw (101,131) node [anchor=north west][inner sep=0.75pt]    {$X$};
\draw (107,87) node [anchor=north west][inner sep=0.75pt]    {$Y$};
\draw (113,40) node [anchor=north west][inner sep=0.75pt]    {$Y'$};
\draw (94,220) node [anchor=north west][inner sep=0.75pt]    {$X'$};
\draw (112,153) node [anchor=north west][inner sep=0.75pt]    {$h$};
\draw (131,68) node [anchor=north west][inner sep=0.75pt]    {$g$};
\draw (154,104) node [anchor=north west][inner sep=0.75pt]    {$Z$};
\draw (131,193) node [anchor=north west][inner sep=0.75pt]    {$s$};
\draw (125,171) node [anchor=north west][inner sep=0.75pt]    {$A$};
\draw (78,105) node [anchor=north west][inner sep=0.75pt]    {$=$};

\end{tikzpicture},
\end{align} 
and it is straightforward to show that such combs are closed under sequential and parallel composition. 

The three resource theories, $\RT$, $\RTKnow$, and $\Shannon$ can all be then viewed as different instances of this construction, all that varies is the starting category $\mathcal{D}$. In the case of $\RT$ we need to start from the SMC $\mathbf{FinSet}$ of finite sets and functions, for $\Shannon$ we start with the SMC $\mathbf{FinStoch}$ of finite sets and stochastic maps, and for $\RTKnow$ we start from an SMC we call $\mathbf{ProbFunc}$ which has finite sets as objects and probability distributions over functions as morphisms\footnote{Formally this can be viewed as the SMC that we obtain by changing the base of enrichment for $\mathbf{FinStoch}$ from being self enriched to being enriched in $\mathbf{FinStoch}$.}.

\subsection{The quantum generalization of our resource theories}

It remains somewhat unclear how, precisely, to model causal influence and knowledge of causal influence in a quantum world.  Recent progress on how to do so is reported in Refs.~\cite{allen2017quantum,Barrett2019}. For a description of the particular challenge of disentangling influence and inference in quantum theory, see Ref.~\cite{omlet}.

This lack of clarity constitutes an obstacle to  working out the quantum generalizations of our results.  Nonetheless, {\em attempting} to uncover such a generalization---that is, quantum resource theories of causal influence and knowledge of causal influence---may provide a novel and insightful angle on the conceptual problem and may ultimately lead to greater clarity. 

In the  case of knowledge of causal influence, one of the central challenges, we believe, is to determine what is the quantum analogue of a probability distribution over functions.  Note that a completely positive trace-preserving (CPTP) map is the quantum analogue of a stochastic map, not the quantum analogue of a distribution over functions. Furthermore, the quantum analogue cannot simply be a probability distribution over unitary maps (even though this may be a special case of said analogue) because only {\em unital} CPTP maps admit of a convex decomposition into unitaries.  Finally, note that finding this analogue is also important for the research program described  in Ref.~\cite[Sec. X]{omlet}, which aims to develop a quantum realist causal-inferential theory.

In this vein, it is worth noting that the quantum analogue of $\Shannon$ is the theory of quantum channels where we allow shared entanglement between the sender and receiver~\cite{hsieh2010entanglement}.  Here, CPTP maps play the role of stochastic maps and consequently quantum Shannon theory is not a resource theory of causal influence or of knowledge of causal influence in the quantum world for the same reasons that $\Shannon$ is not such a resource theory in the classical world, i.e., the reasons presented in Sec.~\ref{sec_Shannon}. 

It is also worth nothing that extensions of quantum Shannon theory wherein one allows of processings of quantum channels by supermaps describing { superpositions of causal orders}\cite{Kristjnsson2020, Milz2022}  also do not provide a resource-theoretic account of causal influence or knowledge of causal influence.  This is because in such approaches, the resources are still taken to be CPTP maps, which are not the quantum analogue of a distribution over functions. In particular, the free resources are defined as those that do not permit signaling between the parties, but as noted above, lack of signaling does not imply lack of causal influence. Thus, the resource theories we are seeking will necessarily differ from the resource theory presented in Ref.~\cite{Milz2022} in its assessment of what are the free resources.  

The categorical perspective described in Sec.~\ref{sec_categorical} also offers some opportunities for making progress on this problem.  Indeed, the distinctions in the categorical presentations of the three classical resource theories we have described above may suggest formal analogues of these distinctions in the quantum sphere.

\section{Conclusions}

We introduced two classical resource theories---one characterizing  causal influence in the presence of perfect information about the functional dependence of $Y$ on $X$, and another characterizing causal influence in the case when one has some uncertainty about the functional dependence.

In the latter resource theory, we ultimately see that resourcefulness consists {\em both} of having knowledge of the function relating one's causal relata {\em and} for the actual function relating the causal relata to be one with large causal connectivity. Thus this resource theory can aptly be summarized as follows:
\begin{center}
\setlength{\fboxsep}{10pt}
  \begin{minipage}{0.8\textwidth} \centering 
    \itshape
    If you have causal influence and you know it, clap your hands.
  \end{minipage}
\end{center}

In both cases, we developed the resource theory extensively for the simplest DAG (with two observed variables with a direct causal influence from one to the other, and no confounders), for arbitrary finite cardinalities. Moreover, we derived a complete set of monotones for variables with arbitrary finite cardinalities in the first resource theory, and for binary variables in the second resource theory.  Extending the second resource theory to the case of arbitrary cardinalities (including infinite cardinalities) is an obvious next direction, and first steps towards it were presented in Sec.~\ref{se:n2n}. After that, the next (and more difficult) direction would be to deal with confounders and ultimately to extend our resource theories to arbitrary DAGs. 

One motivation for understanding these resource theories (even in simple contexts) is to ultimately develop analogous resource theories in the context of quantum causal influence. A first lesson of our work (continuing on from the arguments of  Ref.~\cite{omlet}) is that one must be careful even about what kind of thing is taken to constitute a resource in the context of understanding causal influence. For example, we argued here that it is useful to study probability distributions over functions rather than stochastic maps in the classical case, which suggests that the most insightful object to study in the quantum case is not merely a CPTP map, but rather some other object characterizing uncertainty about more fine-grained causal relations.

\section*{Acknowledgments}

MMA, DS, YY, and RWS 
were supported by Perimeter Institute for Theoretical Physics. Research at Perimeter Institute is supported in part by the Government of Canada through the Department of Innovation, Science and Economic Development and by the Province of Ontario through the Ministry of Colleges and Universities. 
YY and MMA are also supported by the Natural Sciences and Engineering Research Council of Canada (Grant No. RGPIN-2024-04419).
This work is partially carried out under IRA Programme, project no.~FENG.02.01-IP.05-0006/23, financed by the FENG program 2021-2027, Priority FENG.02, Measure FENG.02.01., with the support of the FNP. BZ is supported by the Foundation for Polish Science.
JHS was funded by the European Commission by the QuantERA project ResourceQ under the grant agreement UMO2023/05/Y/ST2/00143.
Some figures were prepared using Mathcha.

\bibliographystyle{quantum}
\bibliography{bib}

\appendix

\section{Guessing probability for $\PP_B^4$ and $\PP_B^5$}\label{app:guess}

Consider a communication channel $X \to Y$ from Alice to Bob. In this appendix, we calculate the probability of Bob guessing the value of $X$ correctly when the channel is characterized by one of the following two probability distributions over functions:
\begin{align}\nonumber
    &\PP_B^4 = \frac{1}{3} \, [\F] + \frac{2}{3} \, [\Rzero],     \\
     &\PP_B^5 = \frac{1}{3} \, [\F] + \frac{1}{3} \, [\Rzero] + \frac{1}{3} \, [\Rone],
\end{align}
which correspond to Eq.~\eqref{eq_Flip_Rzero} and Eq.~\eqref{eq_Flip_Rzero_Rone} in the main text.
Suppose Alice and Bob describe the channel by $\PP_B^4$. This leads to the conditional probability distribution
\begin{equation}\label{eq:sigma}
    \PP_{Y|X=0}=\frac{2}{3}[0]+\frac{1}{3}[1]\quad,\quad
 \PP_{Y|X=1}=[0].
\end{equation}
Since Bob does not have any prior knowledge about $X$, we assume he assigns a uniform prior over $X$, that is, $\PP_X(0)=\PP_X(1)=1/2$. In this case, the conditional probability distribution above leads to $\PP_Y(0)=5/6$. With this, we can calculate
\begin{align}
    &\PP_{Y|X}(0|0)=\frac{\PP_{Y|X}(0|0)\PP_X(0)}{\PP_Y(0)}=\frac{2}{5}, \\
    &\PP_{Y|X}(0|1)=\frac{\PP_{Y|X}(1|0)\PP_X(0)}{\PP_Y(1)}=1.
\end{align}
Therefore, when Bob receives $Y=0$ he should guess $X=1$, and he will have a probability of $3/5$ of being correct. When Bob receives $Y=1$ he should guess $X=0$, and he will have a $100\%$ probability of being correct. Denote $C$ the variable representing is his guess is correct and let $C=\checkmark$ be ``success''.
Overall, his probability of guessing $X$ correctly is 
\begin{equation}
    \PP_C(\checkmark)=\PP_{C|Y}(\checkmark|0)\PP_Y(0)+\PP_{C|Y}(\checkmark|1)\PP_Y(1)=\frac{3}{5}\cdot\frac{5}{6}+1\cdot\frac{1}{6}=\frac{2}{3}.
\end{equation}

Suppose now that Alice and Bob describe the channel by $\PP_B^5$. This leads to the conditional probability distribution
\begin{equation}
     \PP_{Y|X=0}=\frac{2}{3}[0]+\frac{1}{3}[1],\qquad 
    \PP_{Y|X=1}=\frac{1}{3}[0]+\frac{2}{3}[1].
\end{equation}
Recall that we assume that he has a uniform prior over $X$. In this case, the conditional probability distribution above leads to $\PP_Y(0)=\PP_Y(1)=1/2$. With this, we can calculate
\begin{align}
    &\PP_{Y|X}(0|0)=\frac{\PP_{Y|X}(0|0)\PP_X(0)}{\PP_Y(0)}=\frac{1}{3}, \\
    &\PP_{Y|X}(0|1)=\frac{\PP_{Y|X}(1|0)\PP_X(0)}{\PP_Y(1)}=\frac{2}{3}.
\end{align}
Therefore, when Bob receives $Y=0$ he should guess $X=1$, and he will have a probability of $2/3$ of being correct. When Bob receives $Y=1$ he should guess $X=0$, and he will have a probability of $2/3$ of being correct.  Overall, his probability of guessing $X$ correctly is 
\begin{equation}
    \PP_C(\checkmark)=\PP_{C|Y}(\checkmark|0)\PP_Y(0)+\PP_{C|Y}(\checkmark|1)\PP_Y(1)=\frac{2}{3}\cdot\frac{1}{2}+\frac{2}{3}\cdot\frac{1}{2}=\frac{2}{3}.
\end{equation}

Therefore, for both $\PP_B^4$ and $\PP_B^5$, the overall probability of Bob obtaining the correct value of $X$ is $2/3$.  Even though this overall probability is the same, note that  $\PP_B^4$ gives Bob certainty that there is causal connection in at least some specific 
rounds (those where he receives $Y=1$), while  $\PP_B^5$ does not.

\section{Proof of Proposition~\ref{prop:complete}}\label{app:proof-monotones}
Below we recall Proposition~\ref{prop:complete} from the main text and provide its proof.
\Prop*
\begin{proof} \quad \\

\noindent \textbf{$\bullet$ Proof that $M_\beta(\PP_B)$ is a monotone} 

To see that this is a monotone, notice from Table \ref{table1} that any extremal free operation can only decrease $\beta$ or keep it invariant. Because $\beta$ is a linear function of $\PP_B$, it also cannot increase under any mixture of extremal free oeprations. Hence, $M_\beta(\PP_B) \geq M_\beta(\Gamma_{\textrm{free}}[\PP_{B}])$ for any free operation $\Gamma_{\textrm{free}}$. 

\noindent \textbf{$\bullet$ Proof that $M_{|\alpha|}(\PP_{B})$ is a monotone} 

Let us show this for the set of resources where $M_{|\alpha|}(\PP_{B})$ is well defined (i.e., all nonfree resources).
Notice from Table \ref{table1} that all the extremal free operations that take $\PP_{B}$ to a nonfree resource leave the value of $|\alpha|$ unchanged (i.e., they only possibly change the sign of $\alpha$). Since $M_{|\alpha|}(\PP_{B})$ is not a linear function of $\PP_{B}$, it remains to check whether a convex mixture of free operations could generate a resource $\PP_{B}^{*}$ with parameters $(\alpha^*,\beta^*,\gamma^*)$ such that $|\alpha^*|>|\alpha|$. Let us take a fixed but arbitrary free operation 
\begin{align} 
\Gamma' \coloneqq \PP_{\Gamma_\text{pre},\Gamma_\text{post}} = c_{(i)} [\Id,\Id]+c_{(ii)^0} [\Id,\Rzero] + c_{(ii)^1} [\Id,\Rone] + c_{(v)} [\F,\Id] + c_{(vi)}[\Id,\F] + c_{(vii)} [\F,\F]\,,
\end{align}
where the $\{c_{(j)}\}$ are nonnegative coefficients that sum up to $1$, and $[f,g]$ denotes the point distribution where the pre-processing is $\Gamma_\text{pre}=f$ and the post-processing is $\Gamma_\text{post}=g$. 

The probability of $\Id$ given by $\Gamma'[\PP_{B}]$ is 
\begin{align}
  \beta^*\left( \frac{1-\alpha^*}{2} \right) &= \beta \left(\frac{(1-\alpha)(c_{(i)}+c_{(vii)})+(1+\alpha)(c_{(v)}+c_{(vi)})}{2} \right) \nonumber \\
  &=  \beta \left(\frac{(c_{(i)}+c_{(v)}+c_{(vi)}+c_{(vii)})-\alpha(c_{(i)}+c_{(vii)}-c_{(v)}-c_{(vi)})}{2} \right) \,. \label{eq_proofalpha_1}
\end{align}

And  the probability of $\F$ given by  $\Gamma'[\PP_{B}]$  is
\begin{align}
  \beta^*\left( \frac{1+\alpha^*}{2} \right) &= \beta \left(\frac{(1-\alpha)(c_{(v)}+c_{(vi)})+(1+\alpha)(c_{(i)}+c_{(vii)})}{2} \right) \nonumber \\
  &=  \beta \left(\frac{(c_{(i)}+c_{(v)}+c_{(vi)}+c_{(vii)})+\alpha(c_{(i)}+c_{(vii)}-c_{(v)}-c_{(vi)})}{2} \right) \,. \label{eq_proofalpha_2}
\end{align}

Together, Eqs.~\eqref{eq_proofalpha_1} and ~\eqref{eq_proofalpha_2} imply that $\beta^*=(c_{(i)}+c_{(v)}+c_{(vi)}+c_{(vii)})\beta$. Using this back in Eq.~\eqref{eq_proofalpha_1}, we obtain 

\begin{align}
    \alpha^* = \alpha \, \frac{c_{(i)}+c_{(vii)}-c_{(v)}-c_{(vi)}}{c_{(i)}+c_{(vii)}+c_{(v)}+c_{(vi)}}\,.
\end{align}
Now we can apply the triangle inequality 
\begin{align}
    |(c_{(i)}+c_{(vii)})+(-c_{(v)}-c_{(vi)})| \leq |c_{(i)}+c_{(vii)}|+|-c_{(v)}-c_{(vi)}| = c_{(i)} + c_{(vii)} + c_{(v)} + c_{(vi)}
\end{align}
to conclude that $|\alpha^*| \leq |\alpha|$. 
Hence, any convex mixture of such free operations cannot increase the value of $|\alpha|$.

Therefore, $M_{|\alpha|}(\PP_{B}) \geq M_{|\alpha|}(\Gamma_{\textrm{free}}[\PP_{B}])$ for any free operation $\Gamma_{\textrm{free}}$ such that $\Gamma_{\textrm{free}}[\PP_{B}]$ is nonfree, which implies that $M_{|\alpha|}$ is a resource monotone.

\noindent \textbf{$\bullet$ Proof that $M_{|\gamma|,\beta}(\PP_{B})$ is a monotone} 

To see that $M_{|\gamma|,\beta}(\PP_{B})$ is a monotone, we will show that it is the analytic form of the cost-construction of Eq.~\eqref{eq:cost} applied to the set 
\begin{align}
    \mathbf{S} := \left\{\PP_{B}^{\epsilon} \coloneqq \epsilon \, [\Id] + (1 - \epsilon) \, \left( \frac12[\Rzero]+\frac12[\Rone] \right)\right\}_{\epsilon \in [0,1]}
\end{align}
and the function $f=M_\beta$. The set $\mathbf{S} $ is the line that connects $[\Id]$ to $\frac12[\Rzero]+\frac12[\Rone]$ in the tetrahedron, i.e., the line defined by $\gamma=0$, $\alpha=-1$ and varying $\beta$. 
Note that for $\PP_{B}^{\epsilon} \in \mathbf{S}$, $M_\beta(\PP_{B}^{\epsilon}) = \epsilon$. To obtain the cost-construction for a resource $\PP_B$, we will look at the resource $\PP_{B}^{\epsilon}\in \mathbf{S}$ with minimal value of $M_\beta(\PP_{B}^{\epsilon})=\epsilon$ such that $\PP_{B}^{\epsilon} \to \PP_{B}$.

First, note that the extremal free operations map $\PP_{B}^{\epsilon} \in \mathbf{S}$ to the following resources:
\begin{align}\label{eq:down-eps1}
    &\PP_{B}^{\epsilon} \to \PP_{B}^{\epsilon} = \epsilon \, [\Id] + (1 - \epsilon) \, \left( \frac12[\Rzero]+\frac12[\Rone] \right) &(i), (vii) \\
    &\PP_{B}^{\epsilon} \to \PP_{B}^{'\epsilon} \coloneqq \epsilon \, [\F] + (1 - \epsilon) \, \left( \frac12[\Rzero]+\frac12[\Rone] \right) &(v), (vi) \\
    &\PP_{B}^{\epsilon} \to [\Rzero]  & (ii) \\
    &\PP_{B}^{\epsilon} \to [\Rone]  & (ii)
    \label{eq:down-eps2}
\end{align}

Hence, the downward closure polytope $\mathcal{P}_{\downarrow}[\PP_{B}^{\epsilon}]$ is defined by the convex hull of those four distinct points\footnote{These are distinct whenever $\epsilon>0$. When $\epsilon=0$, $\mathcal{P}_{\downarrow}[\PP_{B}^{\epsilon}]$ collapses to the line segment defined by $\Rzero$ and $\Rone$.}, as shown in Fig.~\ref{fig_third_mon_proof}. The condition that $\PP_{B}^{\epsilon} \to \PP_{B}$ means that $\epsilon$ should be chosen such that $\PP_{B} \in \mathcal{P}_{\downarrow}[\PP_{B}^{\epsilon}]$. The condition that $\epsilon$ be minimized imposes that $\PP_{B}$ lie on one of the four facets of $\mathcal{P}_{\downarrow}[\PP_{B}^{\epsilon}]$. This is because, if $\PP_{B}$ is in the interior of $\mathcal{P}_{\downarrow}[\PP_{B}^{\epsilon}]$, then it is possible to construct another tetrahedron $\mathcal{P}_{\downarrow}[\PP_{B}^{\epsilon'}]$ that includes $\PP_B$ and has $\epsilon'<\epsilon$, i.e., such that the points $\PP_B^\epsilon$ and $\PP_B^{'\epsilon}$ of Fig.~\ref{fig_third_mon_proof} are further away from the points $[\Id]$ and $[\F]$ respectively.  

In what follows, we show that $\PP_{B}$ necessarily lies either on the facet $(\PP_{B}^{\epsilon},\PP_{B}^{'\epsilon},[\Rzero])$ or on the facet $(\PP_{B}^{\epsilon},\PP_{B}^{'\epsilon},[\Rone])$ of the tetrahedron $\mathcal{P}_{\downarrow}[\PP_{B}^{\epsilon}]$.

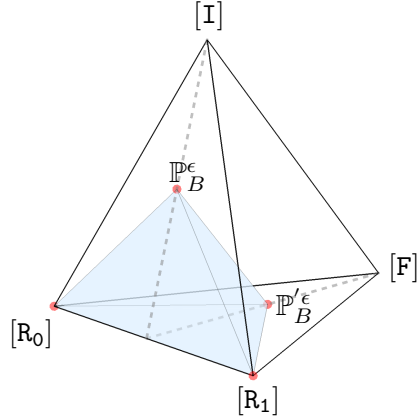
\begin{figure}[htb!]
     \centering
     \includegraphics[width=0.4\linewidth]{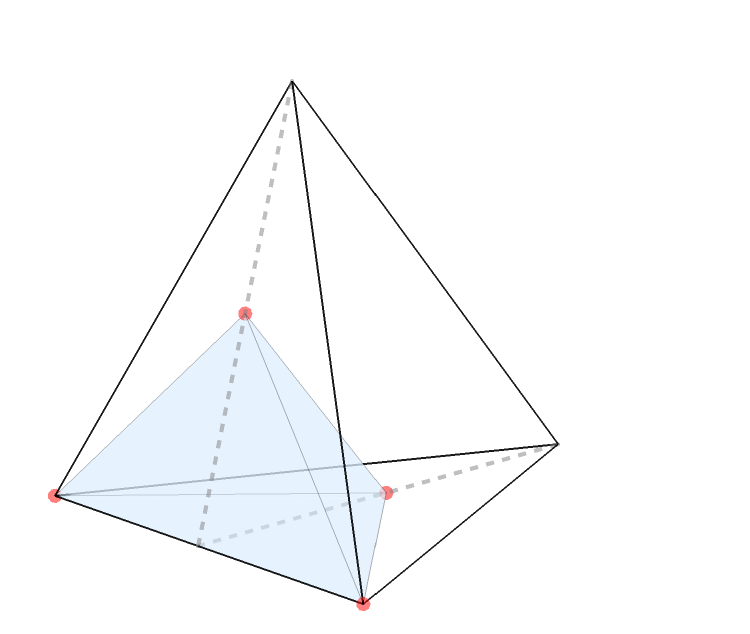}
             \put(-185,20){$[\Rzero]$}
        \put(-100,-5){$[\Rone]$}
        \put(-116,140){$[\Id]$}
        \put(-43,45){$[\F]$}
        \put(-125,80){$\PP^{\epsilon}_B$}
    \put(-85,30){$\PP_{B}^{'\epsilon}$}
     \caption{The downward closure polytope $\mathcal{P}_{\downarrow}[\PP_{B}^{\epsilon}]$. The red dots represent the resources given in Eqs.~\eqref{eq:down-eps1}-~\eqref{eq:down-eps2} and their convex hull is indicated in blue.  }
    \label{fig_third_mon_proof}
 \end{figure}

If $\PP_{B}$ is an interior point of the $([\Id],[\F],[\Rzero],[\Rone])$ tetrahedron, then it clearly lies either on the $(\PP_{B}^{\epsilon},\PP_{B}^{'\epsilon},[\Rzero])$ facet or the $(\PP_{B}^{\epsilon},\PP_{B}^{'\epsilon},[\Rone])$ facet, because the other two facets of $\mathcal{P}_{\downarrow}[\PP_{B}^{\epsilon}]$ are entirely contained on facets of the $([\Id],[\F],[\Rzero],[\Rone])$ tetrahedron. 

Now, consider the case where $\PP_{B}$ belongs to a facet of  the $([\Id],[\F],[\Rzero],[\Rone])$ tetrahedron. If  $\PP_{B}$ belongs to the facet $([\Id], [\F], [\Rzero])$ or to the facet $([\Id], [\F], [\Rone])$, then the only $\mathcal{P}_{\downarrow}[\PP_{B}^{\epsilon}]$ that can possibly include $\PP_B$ is the one where $\epsilon=1$, and thus $\PP_{B}$ lies on either the facet $(\PP_{B}^{\epsilon},\PP_{B}^{'\epsilon},[\Rzero])=([\Id],[\F],[\Rzero])$ or the facet $(\PP_{B}^{\epsilon},\PP_{B}^{'\epsilon},[\Rone])=([\Id],[\F],[\Rone])$.  If  $\PP_{B}$ belongs to the facet $([\Id], [\Rzero], [\Rone])$ or to the facet $([\F], [\Rzero], [\Rone])$, then the condition that $\epsilon$ be minimized imposes that $\PP_B$ lies on one of the line segments defined by $(\PP_B^\epsilon,[\Rzero])$, $(\PP_B^\epsilon,[\Rone])$, $(\PP_B^\epsilon,[\Rzero])$ or $(\PP_B^\epsilon,[\Rone])$. Therefore, $\PP_B$ also lies either on the facet $(\PP_{B}^{\epsilon},\PP_{B}^{'\epsilon},[\Rzero])$ or on the facet $(\PP_{B}^{\epsilon},\PP_{B}^{'\epsilon},[\Rone])$.

Let us now look at each of the cases of $\PP_B$ belonging to the facet $(\PP_{B}^{\epsilon},\PP_{B}^{'\epsilon},[\Rzero])$ or to the facet $(\PP_{B}^{\epsilon},\PP_{B}^{'\epsilon},[\Rone])$, and for each case relate  $\epsilon$ with the parameters $(\alpha, \beta,\gamma)$ that parametrize our resource $\PP_{B}$. Let us start with the case where $\beta \neq 1$, ($\gamma$ is not unspecified) and $\PP_{B}$ is a nonfree resource.
\begin{compactitem}
    \item Let us consider the case where $\PP_{B}$ belongs to the facet defined by ($\PP_{B}^{\epsilon}$, $\PP_{B}^{'\epsilon}$, $[\Rone]$). A generic resource on this facet is parametrized by
    \begin{align}
        c_1 \epsilon [\Id] + c_2 \epsilon [\F] + \frac{1-\epsilon}{2}(c_1+c_2) [\Rzero] + \left(\frac{1-\epsilon}{2}(c_1+c_2) + c_3 \right) [\Rone]\,,
    \end{align}
    where $c_1$, $c_2$, and $c_3$ denote the coefficients in the convex mixture of the resources $\PP_{B}^{\epsilon}$, $\PP_{B}^{'\epsilon}$, and $[\Rone]$, respectively. If $\PP_{B}$ is on that facet, then it must happen that
    \begin{align}
        \beta = \epsilon (c_1 + c_2)\,,\quad \frac{1-\alpha}{2} = \frac{c_1}{c_1+c_2}\,,\quad 1-\gamma = \frac{(1-\epsilon) (c_1+c_2)}{1-\beta}\,.
    \end{align}
    The last equality, together with the normalization condition $c_1+c_2+c_3=1$, yields $c_3 = \gamma (1-\beta)$. Furthermore, note that all points of the facet ($\PP_{B}^{\epsilon}$, $\PP_{B}^{'\epsilon}$, $[\Rone]$) are such that $\gamma \geq 0$. Therefore, we have that $c_3 = |\gamma| (1-\beta)$. From this it follows that 
    \begin{align}
        \frac{\beta}{1 - |\gamma| (1-\beta)} = \frac{\epsilon (c_1+c_2)}{1 - c_3} = \epsilon\,.
    \end{align}
    \item  Let us consider the case where $\PP_{B}$ belongs to the facet defined by ($\PP_{B}^{\epsilon}$, $\PP_{F}^{' \epsilon}$, $[\Rzero]$). A generic resource on this facet is parametrised by
    \begin{align}
        c_1 \epsilon [\Id] + c_2 \epsilon [\F] + \frac{1-\epsilon}{2}(c_1+c_2) [\Rone] + \left(\frac{1-\epsilon}{2}(c_1+c_2) + c_3 \right) [\Rzero]\,,
    \end{align}
    where $c_1$, $c_2$, and $c_3$ denote the coefficients in the convex mixture of the resources $\PP_{B}^{\epsilon}$, $\PP_{B}^{'\epsilon}$, and $[\Rzero]$, respectively. If $\PP_{B}$ is on that facet, then it must happen that
    \begin{align}
        \beta = \epsilon (c_1 + c_2)\,,\quad \frac{1-\alpha}{2} = \frac{c_1}{c_1+c_2}\,,\quad 1+\gamma = \frac{(1-\epsilon) (c_1+c_2)}{1-\beta}\,.
    \end{align}
    The last equality, together with the normalization condition $c_1+c_2+c_3=1$, yields $c_3 = -\gamma (1-\beta)$. Furthermore, note that all points of the facet ($\PP_{B}^{\epsilon}$, $\PP_{B}^{'\epsilon}$, $[\Rzero]$) are such that $\gamma\leq 0$. Therefore, we have that $c_3 = |\gamma| (1-\beta)$. From this it follows that 
    \begin{align}
        \frac{\beta}{1 - |\gamma| (1-\beta)} = \frac{\epsilon (c_1+c_2)}{1 - c_3} = \epsilon\,.
    \end{align}
\end{compactitem}
We see in both cases that 
\begin{align}
    \epsilon = \frac{\beta}{1 - |\gamma| (1-\beta)} =M_{|\gamma|,\beta}(\PP_{B}) \, ,
\end{align}
as we intended to show.

Now, if $\PP_{B}$ is a free resource (i.e., $\alpha$ is unspecified) then $\epsilon=0$ is the solution to the optimization problem in the definition of the cost-monotone, and since in this case $\beta=0$, we can also write
\begin{align}
    0 = \epsilon = \frac{\beta}{1 - |\gamma| (1-\beta)} \,.
\end{align}
Finally, let us discuss the case where $\beta=1$. Here, $\PP_{B}$ lies in the line segment defined by $([\Id], [\F])$, and $\gamma$ is unspecified. For $\PP_{B}$ to belong to $\mathcal{P}_{\downarrow}[\PP_{B}^{\epsilon}]$ then it must happen that $\mathcal{P}_{\downarrow}[\PP_{B}^{\epsilon}]$ is the full $([\Id], [\F], [\Rzero], [\Rone])$ tetrahedron, for which it must be that $\epsilon = 1$. Therefore, the solution to the optimization problem in the definition of the cost-monotone is $\epsilon=1$, which means 
\begin{align}
    \epsilon = 1 \quad \text{if} \quad \beta=1\,.
\end{align}

Putting all these cases together, we see that the cost-monotone $M^{\text{cost}}_{M_\beta,\mathbf{S}}(\PP_{B}) = \epsilon$ from Eq.~\eqref{eq:cost} has an analytical expression given by
\begin{align}
    M^{\text{cost}}_{M_\beta,\mathbf{S}}(\PP_{B}) =  \begin{cases} \frac{\beta}{1-|\gamma|(1-\beta)} \quad &\text{if} \quad \beta \neq 1\,,\\
\qquad 1 \quad &\text{if} \quad \beta = 1\,.
\end{cases}
\end{align}
Hence, the function $M_{|\gamma|,\beta}(\PP_{B}) := M^{\text{cost}}_{M_\beta,\mathbf{S}}(\PP_{B})$ is a resource monotone. 

\noindent \textbf{$\bullet$ Proof that the three monotones form a complete set}

The last thing to show is that $M_\beta$, $M_{|\alpha|}$ and $M_{|\gamma|,\beta}$  form a complete set of monotones for $\RTKnow$ when $X$ and $Y$ are binary.
 Let us begin by considering the case of resources for which all the values of $(M_\beta,M_{|\alpha|},M_{|\gamma|,\beta})$ are defined. In this case, the values of the monotones allow one to fully specify the values of the parameters $(|\alpha|,\beta,|\gamma|)$, or of just $(|\alpha|,\beta)$ when $\gamma$ is undefined. These  in turn are sufficient to specify the equivalence class of the resource, as argued in Sec.~\ref{se:b2brc}. Hence, for resources where all the monotones have defined values, these monotones fully specify the pre-order. 

The resources for which not all the values of $(M_\beta,M_{|\alpha|},M_{|\gamma|,\beta})$ are defined are the free resources. For them, $M_{|\alpha|}$ is not defined and $(M_\beta,M_{|\gamma|,\beta})= (0,0)$, because $\beta=0$. The free resources are the only ones for which $(M_\beta,M_{|\gamma|,\beta})= (0,0)$, so their equivalence class is fully specified by these two monotones.

We can therefore conclude that the values of the monotones $(M_\beta,M_{|\alpha|},M_{|\gamma|,\beta})$ allow one to deduce the location of a resource in the partial order, and so they form a complete set of monotones.

\end{proof}

\section{Proof of Lemma~\ref{lem:canform}}\label{app:proof-lemma}
\Lemma*

\begin{proof}
    In Sec.~\ref{se:b2brc} we showed that the parameters $(|\alpha|,\beta,|\gamma|)$, and hence $(M_\beta,M_{|\alpha|},M_{|\gamma|,\beta})$, are sufficient for characterizing the pre-order of resources, even in the cases when some of those parameters have unspecified values. So we only need to show that these values are necessary. To see this, let us focus on the three monotones $(M_\beta,M_{|\alpha|},M_{|\gamma|,\beta})$.
\begin{compactitem}
    \item Consider two resources $\PP_{B}^{1}$ and $\PP_{B}^{2}$ defined by $(\alpha,\beta,\gamma)$ as follows: 
    $\PP_{B}^{1} \leftrightarrow (\alpha,\beta,\gamma)=(0,0.5,0.3)$ and $\PP_{B}^{2} \leftrightarrow (\alpha,\beta,\gamma)=(0,0.5,0.7)$. These two resources have the same values for $M_{|\alpha|}$ and $M_\beta$. However, $M_{|\gamma|,\beta}(\PP_{B}^{1}) \simeq 0.588 < 0.769 \simeq M_{|\gamma|,\beta}(\PP_{B}^{2})$, and since $M_{|\gamma|,\beta}$ is a monotone, it follows that $\PP_{B}^{1} \not\to \PP_{B}^{2}$. We therefore have two resources that have the same value for $(M_{|\alpha|},M_\beta)$ but do not belong to the same equivalence class. Thus, the information provided by $M_{|\gamma|,\beta}$ is necessary to determine the equivalence class of a resource.
    \item Consider two resources $\PP_{B}^{1}$ and $\PP_{B}^{2}$ defined by $(\alpha,\beta,\gamma)$ as follows: 
    $\PP_{B}^{1} \leftrightarrow (\alpha,\beta,\gamma)=(0,0.5,0.3)$ and $\PP_{B}^{2} \leftrightarrow (\alpha,\beta,\gamma)=(0.5,0.5,-0.3)$. These two resources have the same values for $M_\beta$ and $M_{|\gamma|,\beta}$. However, $M_{|\alpha|}(\PP_{B}^{1}) = 0 < 0.5 = M_{|\alpha|}(\PP_{B}^{2})$, and since $M_{|\alpha|}$ is a monotone, it follows that $\PP_{B}^{1} \not\to \PP_{B}^{2}$. We therefore have two resources that have the same value for $(M_\beta,M_{|\gamma|,\beta})$ but do not belong to the same equivalence class. Thus, the information provided by $M_{|\alpha|}$ is necessary to determine the equivalence class of a resource.
    \item Finally, consider two resources $\PP_{B}^{1}$ and $\PP_{B}^{2}$ defined by $(\alpha,\beta,\gamma)$ as follows: 
    $\PP_{B}^{1}\leftrightarrow (\alpha,\beta,\gamma)=(0,0.5,0)$ and $\PP_{B}^{2} \leftrightarrow (\alpha,\beta,\gamma)=(0,0.25,2/3)$. These two resources have the same values for $M_{|\alpha|}$ and $M_{|\gamma|,\beta}$. However, $M_\beta(\PP_{B}^{1}) \simeq 0.5 > 0.25 \simeq M_\beta(\PP_{B}^{2})$, and since $M_{\beta}$ is a monotone, it follows that $\PP_{B}^{2} \not\to \PP_{B}^{1}$. We therefore have two resources that have the same value for $(M_{|\alpha|},M_{|\gamma|,\beta})$ but do not belong to the same equivalence class. Thus, the information provided by $M_{\beta}$ is necessary to determine the equivalence class of a resource.
\end{compactitem}
We have shown that in some cases, the values of the three monotones $(M_\beta,M_{|\alpha|},M_{|\gamma|,\beta})$ are necessary to determine the pre-order of resources. Since these three values fully specify the values of $(|\alpha|,\beta,|\gamma|)$, this implies that there are cases where the latter are necessary to specify the equivalence classes of resources. 
\end{proof}

\section{Alternative proposal for free operations} 

In this appendix, we define an alternative resource theory where the resources are also probability distributions over functions. The free operations of this alternative resource theory allow pre- and post-processings to depend on what the function being processed actually is. 
You can picture this diagrammatically in the framework of causal-inferential theories~\cite{omlet}, by copying the so-called inferential wire on the given resource and sending that copy into the common cause that may correlate the pre- and post-processings. In this case,  one can faithfully represent the equivalence classes of resources by 
the beta vector $\{\beta_k\}_k$ introduced in Sec.~\ref{se:n2n}, and the set of monotones introduced in Thm.~\ref{thm:betaM} forms a complete set, as we show next.

\begin{lem} \quad \\
In this alternative resource theory of knowledge of causal influence,  the  equivalence classes of resources are faithfully represented by 
the beta vector $\{\beta_k\}_k$ introduced in Sec.~\ref{se:n2n}, where for each $k$, $\beta_k$ represents the total probability assigned to functions with image size $k$ in the decomposition of $\PP_{F}$. 
\end{lem}
\begin{proof}
Let us begin by proving that with this new set of free operations one can freely inter-convert any pair of resources with the same image size. First, notice that the free operations here include the free operations we defined for the resource theory of functions in Sec.~\ref{se:rtf}. Hence, any resource conversion for functions we can do in that resource theory, we can also achieve here. From Prop.~\ref{prop_increase_imagesize} it follows that in the current case we also have that functions with the same image size can be freely inter-converted. Now assume that we have two resources $\PP_{F}$ and $\PP_{F}'$, which have the same image size. Let us decompose them as $\PP_{F} = \sum_{i} c_i [f_i]$ and $\PP_{F}' = \sum_{i} d_i [f_i]$, where both sum run over all the functions in the set  $\mathcal{F}:=\{f_i \, | \, f_i \quad  \text{has the same image size as} \quad \PP_{F} \}$. To see that one can freely transform between $\PP_{F}$ and $\PP_{F}'$, pick $f_0 \in \mathcal{F}$. From Prop.~\ref{prop_increase_imagesize}, $\exists \, \tau_i$ s.t.~$f_i = \tau_i[f_0]$.  Also, notice that the transformation ``reset to $f_0$'' ($\mathtt{R_{f_0}}$) is a free operation. Hence, the free operation $\Gamma_F := \left( \sum_i c_i \tau_i \right) \circ \mathtt{R_{f_0}}$ achieves $\PP_{F} = \Gamma_F[\PP_{F}']$, and the free operation $\Gamma_R := \left( \sum_i d_i \tau_i \right) \circ \mathtt{R_{f_0}}$ achieves $\PP_{F}' = \Gamma_R[\PP_{F}]$. Notice that here is it crucial to know what resources we are interconverting, so we can tailor the free operation to them. This is something we hence couldn't achieve in the Resource theory of Sec.\ref{se:rtfd}. 

Now we need to show that the new set of free operations we have does not allow us to freely increase the image size of a function. For this, recall that the pre- and post-processings that we do on the given function (call it $f$) are functions themselves, since the outcome resource needs to be a function itself. The difference with Sec.~\ref{se:skft} is that here which pre- and post-processing functions we use may depend on which $f$ we are transforming. But this doesn't change the fact that the processings will be functions at the end of the day. Following the same argument as in Sec.~\ref{sec_functionRT_resource_conversion} a post-processing of $f$ by a function cannot produce a new function with a larger image size than $f$, which proves the claim. 

The next step of the proof is to show that a resource $\PP_{F}$ cannot be transformed into one (call it $\PP^{*}_{F^{*}}$) that has a larger image size. This follows from the last claim, noticing that $\PP^{*}_{F^{*}}$ is a probability distribution over functions, and hence we transform a collection of functions (the ones that correspond to the resource $\PP_{F}$) to another collection of functions (the ones that correspond to the resource $\PP^{*}_{F^{*}}$). The pre- and post-processings that we employ for any transformation are functions themselves, and hence we cannot increase the image size of any of the functions in the decomposition of $\PP_{F}$. This means that all the functions pertaining to $\PP^{*}_{F^{*}}$ have image size smaller than or equal to those that pertain to $\PP_{F}$, and therefore $\PP^{*}_{F^{*}}$ cannot have an image size larger than that of $\PP_{F}$.

The last step is to notice that, since free operations cannot increase the image size of a resource, then the sets of equivalent resources that we found at the beginning of the proof (resources with the same image size) are full equivalence classes that do not contain any extra resources beyond those already in each set. 

We see then that to specify the equivalence class of a resource, all we need to know is its image size, for which it's necessary and sufficient to know the beta vector $\{\beta_k\}_k$ introduced in Sec.~\ref{se:n2n}. This vector is fully specified by the set of monotones introduced in Thm.~\ref{thm:betaM}, which are necessary and sufficient to reconstruct $\{\beta_k\}_k$, showing this is a complete set of monotones in this new resource theory.
    
\end{proof}

A resource theory hinged on these free operations has a somewhat contradictory feel. On the one hand, the resources themselves are probability distributions -- representing some agent's ignorance -- over the function describing causal influence, but, on the other hand, the free operations themselves do ``know'' exactly what the function is. It isn't particularly clear in which circumstances (if any) it would make sense for the free operation to have more knowledge than the agent! 

Relative to these free operations,
the resources $[\Id]$ and $\frac{1}{2}[\Id]+\frac{1}{2}[\F]$ are considered equivalent, 
i.e., knowledge of whether the actual functional relation is one or the other is not resourceful anymore. In this sense, this resource theory quantifies the actual causal influence which is happening irrespective of what we know about it, while our original motivation (discussed in Sec.~\ref{sec:examples}) was aimed at also quantifying how much we know about the actual causal influence. Hence, why we primarily focused on the resource theory $\RTKnow$ in the main text.

\end{document}